% Commenced on July 14, 2007
% Revised on 10 February, 2008

%============== Plain TeX ========================

\magnification=\magstep 1

%-----Paper Format---------
\vsize=18.6cm 
\hsize=14.6cm
\hoffset=-0.7cm
\voffset=1cm
\vglue 1cm

%----Headers-------
\def\rightheadline{\hfil{\small Zeta-function and 
Hecke congruence subgroups.\ II}
\hfil\tenrm\folio}
\def\leftheadline{\tenrm\folio\hfil{\scsc Y. Motohashi}\hfil}
\def\firstheadline{\hfil}
\headline{\ifnum\pageno=1 \firstheadline\else
\ifodd\pageno \rightheadline \else \leftheadline\fi\fi}

%-----Footers------
\def\firstpage{\hss{\vbox to 0cm{\vfil\hbox{\rm\folio}}}\hss}
\def\emptyfootline{\hfil}
\footline{\ifnum\pageno=1\firstpage\else
\emptyfootline\fi}

%-----AMS Symbols-------
\input amssym.def
\input amssym.tex

%-----Special Fonts-----
\font\small=cmr8

\font\csc=cmcsc10
\font\scsc=cmcsc10 at 8pt
\font\title=cmbx10 at 12pt

\font\teneusm=eusm10
\font\seveneusm=eusm7
\font\fiveeusm=eusm5
\newfam\eusmfam
\def\eusm{\fam\eusmfam\teneusm}
\textfont\eusmfam=\teneusm
\scriptfont\eusmfam=\seveneusm
\scriptscriptfont\eusmfam=\fiveeusm

\font\teneufm=eufm10
\font\seveneufm=eufm7
\font\fiveeufm=eufm5
\newfam\eufmfam
\def\eufm{\fam\eufmfam\teneufm}
\textfont\eufmfam=\teneufm
\scriptfont\eufmfam=\seveneufm
\scriptscriptfont\eufmfam=\fiveeufm

%----Small Macros----
\def\varGamma{{\mit \Gamma}}
\def\Re{{\rm Re}\,}
\def\Im{{\rm Im}\,}
\def\txt#1{{\textstyle{#1}}}
\def\scr#1{{\scriptstyle{#1}}}
\def\r#1{{\rm #1}}
\def\B#1{{\Bbb #1}}
\def\e#1{{\eusm #1}}

\def\sgn{{\rm sgn}}
\def\f#1{{\eufm #1}}

%----Text-------

\centerline{\title The Riemann Zeta-Function and Hecke Congruence
Subgroups.\ II}
\vskip 1cm
\centerline{\csc Yoichi Motohashi}
\vskip 1cm
\noindent
This is a rework of our old file 
on an explicit spectral decomposition of the mean value
$$
M_2(g;A)=\int_{-\infty}^\infty
\left|\zeta\!\left(\txt{1\over2}+it\right)\right|^4
\left|A\!\left(\txt{1\over2}+it\right)\right|^2g(t)dt
$$
that has been left unpublished since September 1994, though its
summary account is given in [9] (see also [11, Section 4.6]);
here 
$$
A(s)=\sum_n\alpha_nn^{-s}
$$
is a finite Dirichlet series and $g$ is assumed to be
even, regular, real-valued on $\B{R}$, and of fast decay on a sufficiently
wide horizontal strip. At this occasion we shall add greater
details  as well as a rigorous treatment of the Mellin transform
$$
Z_2(s;A)=\int_1^\infty
\left|\zeta\!\left(\txt{1\over2}+it\right)\right|^4
\left|A\!\left(\txt{1\over2}+it\right)\right|^2t^{-s}dt
$$
which was scantly touched on in [9].
\medskip
We shall proceed with an arbitrary $A$ to a
considerable extent but later restrict ourselves to the situation
where $\alpha_n$ is supported by the set of square-free integers.
This is solely to avoid certain technical complexities pertaining
to Kloosterman sums associated with Hecke congruence subgroups which do
not appear particularly worth dealing with thoroughly, for our
present principal purpose is to look into the nature of $Z_2(s;A)$.
\medskip
Our result on $Z_2(s;A)$  seems to allow us to have a glimpse of 
the nature of the plain sixth power moment
$$
M_3(g;1)=
\int_{-\infty}^\infty\left|\zeta\!\left(\txt{1\over2}+it\right)\right|^6
g(t)dt,
$$
although we shall set out only certain ensuing problems which are
to be solved before stating anything precisely. In
fact, this motivation which was implicit in our original file was similar
to that expressed in [4]. Our approach was, however, more explicit, being
a natural extension of our treatment of the plain fourth moment $M_2(g;1)$
that was later published in [11].
\medskip
As we noted at a few occasions, the reason of the success with 
$M_2(g;1)$ lies probably in the fact that the Eisenstein
series in the framework of $\r{SL}(2,\B{R})$ is closely related to the
product of two zeta-values and and in that the group is of real rank one,
with the observation that the later is reflected in that the integral for
$M_2(g;1)$ is single (as is inferred from the arguments developed
in e.g.\ [2][12]). Extrapolating
this, we surmise that a proper formulation of the sixth moment of the
zeta-function might be expressed instead in terms of a double integral,
since the group
$\r{SL}(3,\B{R})$ appears to be closely related to the product of three
zeta-values and it is of real rank 2.
Nevertheless, we shall consider $M_2(g;A)$, as it stands between the
pure fourth and sixth moments and requires less machineries than the
plausible direct approach to  the sixth moment via the spectral theory
of $L^2(\r{PSL}(3,\B{Z})\backslash\r{PSL}(3,\B{R}))$
such as proposed in [11, Section 5.4]. 
\medskip
There are at least three ways for us to proceed along. The first is the
argument that we took in [7][11], the second is a representation
theoretic approach developed in [2], and the third is the one in
[12] which is more representation theoretic and in fact
generalizes to quite a wide extent. We shall take the
first way, as we have indicated above,
for it appears to be the most explicit and allow us to exploit best
the peculiarity of our problem, i.e., the presence of the square of the
zeta-function in place of the first power of an automorphic
$L$-function. However, it should be stressed
that the methods in [2] and [12]
have a definite advantage over that in [7][11]; see {\csc Remark 3} in
Section 15 below. 
\bigskip
\noindent
{\bf Convention.} {\it We shall assume throughout our discussion that
there exist no exceptional eigenvalues for any Hecke congruence 
subgroup $\varGamma_0(q)$\/}.
\medskip
\noindent
Thus all spectral data
$\kappa_j$ should be understood to be real and non-negative. With this,
we might not appear prudent enough, but actually our discussion of
$Z_2(s;A)$ is not essentially affected by the assumption, though we are
aware of the possible existence of poles in the interval
$\left({1\over2},1\right)$.   
\medskip
\noindent
{\csc Remark 1}. Readers are warned of a number of notational conflicts,
none of which should, however, cause any serious misunderstanding. We
remark also that our discussion contains details which must be often
excessive for experts; nevertheless, we do this because our old file had
been  prepared for an abortive series of lectures to be given to
beginners,  and we want to keep the original style. 
By the way, there exists as well an abridged version of the file that was
to be included in [11] as its sixth chapter, but the plan was put
away because of a reason which we can no longer remember.
\medskip
\noindent
{\csc Remark 2}. We do not mention any of works
on mean values of automorphic $L$-functions done in recent years, 
notably by D. Goldfeld and his colleagues, some of which
in fact come close to our interest on $Z_2(s;A)$. This is solely
due to our wish to keep ourselves within the framework of the unpublished
file  of ours; the necessary updating will be made in our relevant
forthcoming works. 
\medskip
In passing, we stress that our work [8] (see also
[11, Section 5.3]) on $Z_2(s;1)$ was done without any knowledge of the
existence of A. Good's work [5] on the Mellin transform of the square of
an arbitrary automorphic $L$-function. His argument depends on a
clever choice of a Poincar\'e series, whereas ours exploits fully the
peculiarity of the Riemann zeta-function as indicated above and
produces results more explicit than his. We add that our  
reasoning extends beyond Good's situation. This is 
a consequence of our latest work [12]  lying on the lines
developed in [2], [7], and [11]. 

\vskip 1cm
\noindent
{\bf 1.} To begin with, we have
$$
M_2(g;A)=\sum_{\scr{a,b,c}\atop \scr{(a,b)=1}}
{\alpha_{ac}\overline{\alpha_{bc}}\over c\sqrt{ab}}I_2(g;b/a),\eqno(1.1)
$$
where
$$
I_2(g;b/a)=\int_{-\infty}^\infty\left|\zeta\!\left(\txt{1\over2}
+it\right)\right|^4(b/a)^{it}g(t)dt.\eqno(1.2)
$$
\medskip
To study the latter we introduce
$$
I(u,v,w,z;g;b/a)=-i\int_{(0)}
\zeta(u+t)\zeta(v-t)\zeta(w+t)\zeta(z-t)(b/a)^{t}g(-it)dt\eqno(1.3)
$$
with $(a,b)=1$ and $\Re u,\ldots,\Re z>1$.
Shifting the contour to $(\alpha)$ lying in the far right, we have
$$
\eqalignno{
I(u,v,w,z&;g;b/a)=-i\int_{(\alpha)}\cdots dt\cr
+2\pi\Big\{&\zeta(u+v-1)
\zeta(v+w-1)\zeta(z-v+1)(b/a)^{v-1}g(i(1-v))\cr
+\;&\zeta(u+z-1)
\zeta(w+z-1)\zeta(v-z+1)(b/a)^{z-1}g(i(1-z))\Big\};&(1.4)
}
$$
Thus $I(u,v,w,z;g;b/a)$ is meromorphic throughout $\B{C}^4$.
With this, we assume that $\Re u,\ldots,\break\Re z<1$ and shift the last
contour back to the original, getting
$$
\eqalignno{
I(u,v,w,z&;g;b/a)=-i\int_{(0)}\cdots dt\cr
+2\pi\Big\{&\zeta(u+v-1)
\zeta(v+w-1)\zeta(z-v+1)(b/a)^{v-1}g(i(1-v))\cr
+\;&\zeta(u+z-1)
\zeta(w+z-1)\zeta(v-z+1)(b/a)^{z-1}g(i(1-z))\cr
+\;&\zeta(w-u+1)\zeta(u+v-1)\zeta(u+z-1)(b/a)^{1-u}g(i(u-1))\cr
+\;&\zeta(u-w+1)
\zeta(v+w-1)\zeta(w+z-1)(b/a)^{1-w}g(i(w-1))\Big\}.&(1.5)
}
$$
In the vicinity of $p_{1\over2}=\left({1\over2},{1\over2},
{1\over2},{1\over2}\right)$, the part in the braces is
equal to
$$
\eqalignno{
&\zeta(u+v-1)\zeta(v+w-1){1\over z-v}(1+c_E(z-v)+\cdots)(b/a)^{v-1}
g(i(1-v))\cr
&+\zeta(u+v-1)\left(1+{\zeta'\over\zeta}(u+v-1)(z-v)+\cdots\right)\cr
&\times\zeta(v+w-1)\left(1+{\zeta'\over\zeta}(v+w-1)(z-v)+\cdots\right)
{1\over v-z}(1+c_E(v-z)+\cdots)\cr
&\times(b/a)^{v-1}\left(1+(\log b/a)(z-v)+\cdots\right)g(i(1-v))
\left(1-i{g'\over g}(i(1-v))(z-v)+\cdots\right)\cr
&+{1\over w-u}(1+c_E(w-u)+\cdots)\zeta(u+v-1)\zeta(u+z-1)(b/a)^{1-u}
g(i(u-1))\cr
&+{1\over u-w}(1+c_E(u-w)+\cdots)\zeta(u+v-1)
\left(1+{\zeta'\over\zeta}(u+v-1)(w-u)+\cdots\right)\cr
&\times\zeta(u+z-1)\left(1+{\zeta'\over\zeta}(u+z-1)(v-u)+\cdots\right)\cr
&\times(b/a)^{1-u}\left(1+(\log b/a)(u-w)+\cdots\right)g(i(u-1))
\left(1+i{g'\over g}(i(u-1))(w-u)+\cdots\right)\cr
=&\,\zeta(u+v-1)\zeta(v+w-1)(b/a)^{v-1}g(i(1-v))\cr
&\times\Big\{2c_E-{\zeta'\over\zeta}(u+v-1)-{\zeta'\over\zeta}(v+w-1)
-\log b/a +{g'\over g}(i(1-v))\Big\}+O(|z-v|)\cr
&+\zeta(u+v-1)\zeta(u+z-1)(b/a)^{1-u}g(i(u-1))\cr
&\times\Big\{2c_E-{\zeta'\over\zeta}(u+v-1)-{\zeta'\over\zeta}(u+z-1)
+\log b/a +{g'\over g}(i(u-1))\Big\}+O(|u-w|),&(1.6)
}
$$
where $c_E$ is the Euler constant.
Hence, in particular, $I(u,v,w,z;g;b/a)$ is regular in a neighbourhood
of $p_{1\over2}$, and  we get
$$
\eqalignno{
I_2(g;b/a)&=I\!\left(p_{1\over2};g;b/a\right)\cr
&-{\pi\over2}(b/a)^{-1/2}g\left(\txt{1\over2}i\right)
\left\{2c_E-2\log 2\pi -\log(b/a)+i{g'\over
g}\!\left(\txt{1\over2}i\right)\right\} \cr
&-{\pi\over2}(a/b)^{-1/2}g\left(\txt{1\over2}i\right)
\left\{2c_E-2\log 2\pi -\log(a/b)-i{g'\over
g}\!\left(\txt{1\over2}i\right)\right\}.&(1.7)
}
$$
The last two terms can be regarded as practically negligible.
\bigskip\noindent
{\bf 2.} On the other hand, we have, in the region of absolute
convergence,
$$
\eqalignno{
I(u,v,w,z;g;b/a)&=\sum_{k,l,m,n}{1\over k^u l^v m^w n^z}
\hat{g}\left(\log{bln\over akm}\right)\cr
&=\zeta(u+v)\sum_{\scr{k,l}\atop\scr{(k,l)=1}}{1\over k^ul^v}\sum_{m,n}
{1\over m^wn^z}\hat{g}\left(\log{bln\over akm}\right)
,&(2.1)
}
$$
where $\hat{g}$ is the Fourier transform of $g$; and 
$(ak,bl)=(a,l)\cdot(b,k)=c\cdot d$, say; note that $(a,b)=1$. We have
$$
\eqalignno{
I(u,v,w,z;g;b/a)&=\zeta(u+v)\sum_{c|a,d|b}{1\over c^vd^u}
\sum_{\scr{k,l, (k,l)=1}\atop \scr{(a/c,l)=1, (b/d,k)=1}}
{1\over k^ul^v}\sum_{m,n}{1\over m^wn^z}
\hat{g}\left(\log{bln/d\over akm/c}\right)\cr
&=\zeta(u+v){1\over a^vb^u}\sum_{c|a,d|b}c^vd^u
\sum_{\scr{k,l, (k,l)=1}\atop \scr{(c,l)=1, (d,k)=1}}
{1\over k^ul^v}\sum_{m,n}{1\over m^wn^z}
\hat{g}\left(\log{dln\over ckm}\right)\cr
&=\zeta(u+v){1\over a^vb^u}\sum_{c|a,d|b}c^vd^u
\sum_{\scr{k,l}\atop \scr{(ck,dl)=1}}
{1\over k^ul^v}\sum_{m,n}{1\over m^wn^z}
\hat{g}\left(\log{dln\over ckm}\right)\cr
&=\zeta(u+v){1\over a^vb^u}\sum_{c|a,d|b}c^vd^u
J(u,v,w,z;g;d/c) ,&(2.2)
}
$$
say.
\medskip
Then we apply the dissection:
$$
J(u,v,w,z;g;d/c)=\left\{J_0+J_++J_-\right\}(u,v,w,z;g;d/c),\eqno(2.3)
$$
where $J_-(u,v,w,z;g;d/c)=J_+(v,u,z,w;g;c/d)$ and
$$
\eqalignno{
J_0(u,v,w,z;g;d/c)&=\hat{g}(0)\sum_{\scr{k,l}\atop\scr{(ck,dl)=1}}
{1\over k^ul^v}\sum_{\scr{m,n}\atop\scr{ckm=dln}}{1\over m^wn^z},
&(2.4)\cr
J_+(u,v,w,z;g;d/c)&=\sum_{\scr{k,l}\atop\scr{(ck,dl)=1}}
{1\over k^ul^v}\sum_{\scr{m,n}\atop\scr{ckm>dln}}{1\over m^wn^z}
\hat{g}\left(\log{dln\over ckm}\right).
&(2.5)
}
$$
We have
$$
\eqalignno{
&J_0(u,v,w,z;g;d/c)=\hat{g}(0)\sum_{\scr{k,l}\atop\scr{(ck,dl)=1}}
{1\over k^ul^v}\sum_n{1\over(dln)^w(ckn)^z}\cr
&=\hat{g}(0)c^{-z}d^{-w}\zeta(w+z)
\sum_{\scr{k,l}\atop\scr{(ck,dl)=1}}
{1\over k^{u+z}l^{v+w}}\cr
&=\hat{g}(0)c^{-z}d^{-w}\zeta(w+z)
\sum_{k,l}{1\over k^{u+z}l^{v+w}}\sum_{r|(ck,dl)}\mu(r)\cr
&=\hat{g}(0)c^{-z}d^{-w}\zeta(w+z)\sum_r\mu(r)
\sum_{r/(c,r)|k}{1\over k^{u+z}}
\sum_{r/(d,r)|l}{1\over l^{v+w}}\cr
&=\hat{g}(0)c^{-z}d^{-w}\zeta(w+z)\zeta(u+z)\zeta(v+w)
\sum_r\mu(r)((c,r)/r)^{u+z}((d,r)/r)^{v+w}\cr
&=\hat{g}(0)c^{-z}d^{-w}\zeta(w+z)\zeta(u+z)\zeta(v+w)\cr
&\hskip 1cm\times\prod_{p\nmid cd}\left(1-{1\over
p^{u+v+w+z}}\right)\prod_{p|c}\left(1-{1\over p^{v+w}}\right)
\prod_{p|d}\left(1-{1\over p^{u+z}}\right),&(2.6)
}
$$
where $p$ denotes a generic prime and the condition $(c,d)=1$ has been
used. The contribution of $J_0(u,v,w,z;g;d/c)$ to
$I(u,v,w,z;g;b/a)$ is thus equal to
$$
\eqalignno{
\hat{g}(0)&a^{-v}b^{-u}
{\zeta(u+v)\zeta(u+z)\zeta(w+v)\zeta(w+z)\over\zeta(u+v+w+z)}\cr
&\times\left\{\sum_{c|a}c^{v-z}\prod_{p|c}{1-\displaystyle{1\over
p^{v+w}}\over 1-\displaystyle{1\over p^{u+v+w+z}}}\right\} 
\left\{\sum_{d|b}d^{u-w}\prod_{p|d}{1-\displaystyle{1\over
p^{u+z}}\over 1-\displaystyle{1\over p^{u+v+w+z}}}\right\}. &(2.7)
}
$$
\bigskip\noindent{\bf 3.}
Next, we shall consider the non-diagonal part $J_+$. We have
$$
\eqalignno{
&J_+(u,v,w,z;g;d/c)\cr
=&
\sum_{\scr{k,l}\atop\scr{(ck,dl)=1}}
{1\over k^ul^v}\sum_f\sum_{\scr{m,n}\atop\scr{ckm=dln+f}}{1\over m^wn^z}
\hat{g}\left(\log{dln\over ckm}\right)\cr
=&\sum_{\scr{k,l}\atop\scr{(ck,dl)=1}}
{1\over k^ul^v}\sum_f\sum_{\scr{m,n}\atop\scr{ckm=dln+f}}
{(ck)^w\over (ckm)^wn^z}
\hat{g}\left(\log{dln\over ckm}\right)\cr
=&\sum_{\scr{k,l}\atop\scr{(ck,dl)=1}}
{(ck)^w\over k^ul^v}\sum_f\sum_{\scr{n}\atop\scr{dln+f\equiv 0 \bmod\,
ck}} {1\over (dln+f)^wn^z}
\hat{g}\left(\log{dln\over dln+f}\right)\cr
=&(c/d)^w\sum_{\scr{k,l}\atop\scr{(ck,dl)=1}}
{1\over k^{u-w}l^{v+w}}\sum_f\sum_{n\equiv -\overline{dl}f \bmod\,
ck}{1\over n^{w+z}}\left(1+{f\over dln}\right)^{-w}
\hat{g}\left(\log\left(1+{f\over dln}\right)
\right).\quad&(3.1)
}
$$
\medskip
We introduce the Mellin transform
$$
\eqalignno{
g^*\!(s,w)&=\int_0^\infty \hat{g}(\log(1+x)){x^{s-1}
\over(1+x)^w}dx\cr
&=\Gamma(s)\int_{-\infty}^\infty 
{\Gamma(w-s+it)\over\Gamma(w+it)}g(t)dt, &(3.2)
}
$$
provided $\Re w>\Re s>0$. Shifting the last contour downward
appropriately, we see that $g^*(s,w)/\Gamma(s)$ is entire in
$s,w$; and an upward shift gives that $g^*(s,w)$ is of rapid decay in
$s$ as far as $w$ and $\Re s$ are bounded (see [11, Lemma 4.1]). In
particular, we have
$$
\eqalignno{
J_+(u,v,w,z;g;d/c)&={(c/d)^w\over2\pi
i}\sum_{\scr{k,l}\atop\scr{(ck,dl)=1}} {1\over
k^{u-w}l^{v+w}}\cr
&\times\sum_f\sum_{n\equiv -\overline{dl}f \bmod\, ck}{1\over
n^{w+z}}\int_{(\eta)}g^*(s,w)\left(f\over dln\right)^{-s}ds, &(3.3)
}
$$
with $\eta>0$, which converges absolutely if
$$
\eta>1,\; \Re u>\Re w+1,\; \Re(v+w)>\eta+1,\; \Re(w+z)>\eta+1.\eqno(3.4)
$$
On this condition, we have
$$
\eqalignno{
J_+&(u,v,w,z;g;d/c)\cr
=&{(c/d)^w\over2\pi
i}\int_{(\eta)}g^*(s,w)d^s\sum_{\scr{k,l}\atop\scr{(ck,dl)=1}} {1\over
k^{u-w}l^{v+w-s}}\sum_f{1\over f^s}\sum_{n\equiv -\overline{dl}f \bmod\,
ck}{1\over n^{w+z-s}}\,ds\cr
=&{c^{-z}d^{-w}\over2\pi
i}\int_{(\eta)}g^*(s,w)(cd)^s
\sum_{\scr{k,l}\atop\scr{(ck,dl)=1}} {1\over
k^{u+z-s}l^{v+w-s}}\sum_f{1\over f^s}
\zeta\left(w+z-s,-{\overline{dl}f\over ck}\right)ds, &(3.5)
}
$$
where $\zeta(s,\omega)$ is the Hurwitz zeta-function. Classifying
$l$ into residue classes $\bmod\,\,ck$, we have
$$
\eqalignno{
&J_+(u,v,w,z;g;d/c)={c^{-z}d^{-w}\over2\pi
i}\int_{(\eta)}g^*(s,w)(cd)^s\sum_{(k,d)=1}
{1\over k^{u+z-s}}\cr
&\hskip 1cm\times
\sum_f{1\over f^s}
\sum_{\scr{h=1}\atop\scr{(h,ck)=1}}^{ck}\sum_{l\equiv
h\bmod\, ck} {1\over l^{v+w-s}}\zeta\left(w+z-s,-{\overline{dh}f\over
ck}\right)ds\cr
&={c^{-v-w-z}d^{-w}\over2\pi
i}\int_{(\eta)}g^*(s,w)(c^2d)^s\sum_{(k,d)=1}
{1\over k^{u+v+w+z-2s}}\cr
&\hskip
1cm\times\sum_f{1\over f^s}\sum_{\scr{h=1}\atop\scr{(h,ck)=1}}^{ck}
\zeta\left(v+w-s,{h\over ck}\right)
\zeta\left(w+z-s,-{\overline{dh}f\over
ck}\right)ds.&(3.6)
}
$$
\bigskip\noindent
{\bf 4.} We are going to shift the last contour. To this end we
assume that there exists a large $\eta_1$ such that $\eta_1>
\eta+1$ and
$$
\Re(v+w)<\eta_1,\; \Re(w+z)<\eta_1,\; \Re(u+v+w+z)>2(\eta_1+1).
\eqno(4.1)
$$
On this and $1<\Re s<\eta_1+\varepsilon$ with a small $\varepsilon>0$,
the sum
$$
\sum_{(k,d)=1}
{1\over k^{u+v+w+z-2s}}\sum_f{1\over f^s}
\sum_{\scr{h=1}\atop\scr{(h,ck)=1}}^{ck}
\zeta\left(v+w-s,{h\over ck}\right)
\zeta\left(w+z-s,-{\overline{dh}f\over ck}\right)\eqno(4.2)
$$
is a meromorphic function of the five complex variables. To see this
we note that for any finite $s$
$$
\zeta(s,\omega)\ll |s-1|^{-1}+\omega^{-\Re s} \qquad
(0<\omega\le1),
\eqno(4.3)
$$
as it follows via an application of partial summation to the
Dirichlet series defining $\zeta(s,\omega)$.
Thus $(4.2)$ is, provided neither $v+w-s$ nor $w+z-s$ is
too close to $1$,
$$
\eqalignno{
\ll&\sum_k k^{\Re(2s-u-v-w-z)+1}\left(1+k^{\Re(v+w-s)}\right)
\left(1+k^{\Re(w+z-s)}\right)\cr
=&\sum_k\left\{k^{\Re(2s-u-v-w-z)+1}+k^{\Re(s-u-z)+1}
+k^{\Re(s-u-v)+1}+k^{\Re(w-u)+1}\right\},&(4.4)
}
$$
in which we have
$$
\eqalignno{
&\Re(2s-u-v-w-z)+1<\Re(2s)-2(\eta_1+1)+1,\cr
&\Re(s-u-z)+1=\Re(s)-\Re(u+v+w+z)+\Re(v+w)+1\cr
&\hskip 2cm<\Re(s)-2(\eta_1+1)+\eta_1+1,\cr
&\Re(s-u-v)+1=\Re(s)-\Re(u+v+w+z)+\Re(w+z)+1\cr
&\hskip 2cm < \Re(s)-2(\eta_1+1)+\eta_1+1,\cr
&\Re(w-u)+1=\Re(v+w)+\Re(w+z)-\Re(u+v+w+z)+1\cr
&\hskip 2cm <\eta_1+\eta_1-2(\eta_1+1)+1;&(4.5)
}
$$
and the assertion follows.
\medskip
With this, we shift the contour in $(3.6)$ to $(\eta_1)$. We encounter
poles at $s=w+z-1,v+w-1$; we may assume without loss of generality that
they do not coincide. Before computing the residues,
we note that
$$
\sum_{\scr{h=1}\atop\scr{(h,q)=1}}\zeta(s,hm/q)=
\zeta(s)\sum_{\delta|q}\delta\mu(q/\delta)(\delta/(\delta,m))^{s-1}.
\eqno(4.6)
$$
To show this we use the functional
equation
$$
\zeta(s,\omega)=2(2\pi)^{s-1}\Gamma(1-s)\sum_n\sin\left(\txt{1\over2}\pi
s+2\pi n\omega\right)n^{s-1}\qquad(\Re s<0).\eqno(4.7)
$$
Thus, for $\Re s<0$,
$$
\eqalignno{
\sum_{\scr{h=1}\atop\scr{(h,q)=1}}^q\zeta(s,hm/q)&=
2\Gamma(1-s)(2\pi)^{s-1}\sum_n\sum_{\scr{h=1}\atop\scr{(h,q)=1}}^q
\sin\left(\txt{1\over2}\pi
s+2\pi{h\over q} mn\right)n^{s-1}\cr
&=2\Gamma(1-s)(2\pi)^{s-1}\sin\left(\txt{1\over2}\pi s\right)
\sum_n c_q(mn)n^{s-1}\cr
&=2\Gamma(1-s)(2\pi)^{s-1}\sin\left(\txt{1\over2}\pi s\right)
\sum_n n^{s-1}\sum_{\delta|(q,mn)}\delta\mu(q/\delta)\cr
&=2\Gamma(1-s)(2\pi)^{s-1}\sin\left(\txt{1\over2}\pi s\right)
\zeta(1-s)\sum_{\delta|q}\delta\mu(q/\delta)(\delta/(\delta,m))^{s-1},
\qquad&(4.8)
}
$$
with the Ramanujan sum $c_q \bmod\, q$; and $(4.6)$ follows via the
functional equation for $\zeta$.
\medskip
Let us compute the residue at $s=w+z-1$. This is equal to
$$
\eqalignno{
&2\pi i(c^2d)^{w+z-1}g^*(w+z-1,w)\sum_{(k,d)=1}
{1\over k^{u+v-w-z+2}}\sum_f{1\over f^{w+z-1}}\cr
&\hskip2cm\times\sum^{ck}_{\scr{h=1}\atop\scr{(h,ck)=1}}
\zeta\left(v-z+1,{h\over ck}\right)\cr
=&2\pi i(c^2d)^{w+z-1}g^*(w+z-1,w)\zeta(w+z-1)\zeta(v-z+1)\cr
&\hskip2cm\times\sum_{(k,d)=1}{1\over k^{u+v-w-z+2}}
\sum_{\delta|ck}\delta^{v-z+1}\mu(ck/\delta)\cr
=&2\pi ic^{2w+v+z-1}d^{w+z-1}g^*(w+z-1,w)\zeta(w+z-1)\zeta(v-z+1)\cr
&\hskip2cm\times\sum_{(k,d)=1}{1\over k^{u-w+1}}
\prod_{p|ck}\left(1-{1\over p^{v-z+1}}\right)\cr
=&2\pi ic^{2w+v+z-1}d^{w+z-1}g^*(w+z-1,w)\zeta(w+z-1)\zeta(v-z+1)\cr
&\times\prod_{p\nmid cd}\left(1+{1\over p^{u-w+1}}
\left(1-{1\over p^{v-z+1}}\right)\left(1-{1\over
p^{u-w+1}}\right)^{-1}\right)\cr
&\hskip2cm\times
\prod_{p|c}\left(1-{1\over p^{v-z+1}}\right)
\left(1-{1\over p^{u-w+1}}\right)^{-1}.&(4.9)
}
$$
Returning to $(3.6)$, we see that the contribution of the residue to
$J_+(u,v,w,z,;g;d/c)$ is
$$
\eqalignno{
&c^{w-1}d^{z-1}g^*(w+z-1,w){\zeta(v-z+1)\zeta(w+z-1)\zeta(u-w+1)\over
\zeta(u+v-w-z+2)}\cr
&\times\prod_{p|c}\left({1-\displaystyle{1\over p^{v-z+1}}
\over 1-\displaystyle{1\over p^{u+v-w-z+2}}}\right)
\prod_{p|d}\left({1-\displaystyle{1\over p^{u-w+1}}
\over 1-\displaystyle{1\over p^{u+v-w-z+2}}}\right).&(4.10)
}
$$
\bigskip\noindent
{\bf 5.} The residue at $s=v+w-1$ is equal to
$$
\eqalignno{
&2\pi i(c^2d)^{v+w-1}g^*(v+w-1,w)\sum_{(k,d)=1}
{1\over k^{u-v-w+z+2}}\cr
&\hskip2cm\times\sum_f{1\over
f^{v+w-1}}\sum^{ck}_{\scr{h=1}\atop\scr{(h,ck)=1}}
\zeta\left(z-v+1,-{\overline{dh}f\over ck}\right)\cr
=&2\pi i(c^2d)^{v+w-1}g^*(v+w-1,w)\zeta(z-v+1)\sum_{(k,d)=1}
{1\over k^{u-v-w+z+2}}\cr
&\hskip2cm\times\sum_{\delta|ck}\delta^{z-v+1}\mu(ck/\delta)
\sum_f{(f,\delta)^{v-z}\over f^{v+w-1}},&(5.1)
}
$$
as before. Here 
$$
\eqalignno{
\sum_f{(f,\delta)^{v-z}\over f^{v+w-1}}&=\sum_f{1\over f^{v+w-1}}
\sum_{\lambda|(f,\delta)}\lambda^{v-z}\prod_{p|\lambda}
\left(1-{1\over p^{v-z}}\right)\cr
&=\zeta(v+w-1)\sum_{\lambda|\delta}{1\over \lambda^{w+z-1}}
\prod_{p|\lambda}\left(1-{1\over p^{v-z}}\right),&(5.2)
}
$$
and
$$
\eqalignno{
&\sum_{\delta|ck}\delta^{z-v+1}\mu(ck/\delta)
\sum_f{(f,\delta)^{v-z}\over f^{v+w-1}}\cr
=&\zeta(v+w-1)\sum_{\lambda|ck}{1\over\lambda^{w+z-1}}
\prod_{p|\lambda}\left(1-{1\over p^{v-z}}\right)
\sum_{\delta|(ck)/\lambda}(\delta\lambda)^{z-v+1}
\mu(ck/\delta\lambda)\cr
=&(ck)^{z-v+1}\zeta(v+w-1)
\sum_{\lambda|ck}{1\over\lambda^{w+z-1}}
\prod_{p|\lambda}\left(1-{1\over p^{v-z}}\right)\prod_{p|(ck)/\lambda}
\left(1-{1\over p^{z-v+1}}\right).
&(5.3)
}
$$
Thus
$$
\eqalignno{
&\sum_{(k,d)=1}
{1\over k^{u-v-w+z+2}}
\sum_{\delta|ck}\delta^{z-v+1}\mu(ck/\delta)
\sum_f{(f,\delta)^{v-z}\over f^{v+w-1}}\cr
=&c^{z-v+1}\zeta(v+w-1)\sum_{(k,d)=1}{1\over k^{u-w+1}}
\sum_{\lambda|ck}{1\over\lambda^{z+w-1}}
\prod_{p|\lambda}\left(1-{1\over
p^{v-z}}\right)\prod_{p|(ck)/\lambda}
\left(1-{1\over p^{z-v+1}}\right)\cr
=&c^{z-v+1}\zeta(v+w-1)\sum_{(\lambda,d)=1}{(c,\lambda)^{u-w+1}
\over\lambda^{u+z}}\prod_{p|\lambda}\left(1-{1\over p^{v-z}}\right)\cr
&\hskip2cm\times\sum_{(k,d)=1}{1\over k^{u-w+1}}
\prod_{p|(ck)/(c,\lambda)}\left(1-{1\over p^{z-v+1}}\right)\cr
=&c^{z-v+1}\zeta(v+w-1)\sum_{(\lambda,d)=1}{(c,\lambda)^{u-w+1}
\over\lambda^{u+z}}\prod_{p|\lambda}\left(1-{1\over p^{v-z}}\right)\cr
&\hskip2cm\times\prod_{p\nmid (cd)/(c,\lambda)}
\left(1+{1\over p^{u-w+1}}\left(1-{1\over p^{z-v+1}}\right)
\left(1-{1\over p^{u-w+1}}\right)^{-1}\right)\cr
&\hskip2cm\times\prod_{p|c/(c,\lambda)}
\left(1-{1\over p^{z-v+1}}\right)\left(1-{1\over p^{u-w+1}}\right)^{-1}\cr
=&c^{z-v+1}{\zeta(v+w-1)\zeta(u-w+1)\over\zeta(u-v-w+z+2)}\prod_{p|d}
\left(1-\displaystyle{1\over p^{u-w+1}}
\over 1-\displaystyle{1\over p^{u-v-w+z+2}}\right)\cr
&\hskip2cm\times\sum_{(\lambda,d)=1}{(c,\lambda)^{u-w+1}
\over\lambda^{u+z}}\prod_{p|\lambda}\left(1-{1\over p^{v-z}}\right)
\prod_{p|c/(c,\lambda)}
\left(1-\displaystyle{1\over p^{z-v+1}}
\over 1-\displaystyle{1\over p^{u-v-w+z+2}}\right).
&(5.4)
}
$$
In the last sum we write $\lambda=\lambda_1\lambda_2$ with
$(\lambda_1,c)=1$ and $\lambda_2|c^\infty$; and we see that
the sum is equal to
$$
\eqalignno{
&{\zeta(u+z)\over\zeta(u+v)}\prod_{p|cd}
\left({1-\displaystyle{1\over p^{u+z}}\over
1-\displaystyle{1\over p^{u+v}}}\right)\cr
&\times\prod_{p^\beta\Vert c}
\left(\sum_{j=0}^\infty{p^{(u-w+1)\min(\beta,j)}\over p^{j(u+z)}}
\left(1-{1\over p^{v-z}}\right)^{\xi(j)}
\left({1-\displaystyle{1\over p^{z-v+1}}\over
1-\displaystyle{1\over
p^{u-v-w+z+2}}}\right)^{\xi(\beta-\min(\beta,j))}\right),\qquad&(5.5) 
}
$$
where $1-\xi$ is the unit measure placed at the origin. 
One could compute the last sum into a finite expression.
\medskip
The contribution of the residue at $s=v+w-1$ to $J_+(u,v,w,z;g;d/c)$
is equal to
$$
\eqalignno{
c^{w-1}d^{v-1}&g^*(v+w-1,w){\zeta(z-v+1)\zeta(v+w-1)\zeta(u-w+1)
\zeta(u+z)\over\zeta(u+v)\zeta(u-v-w+z+2)}\cr
&\times \prod_{p|cd}
\left({1-\displaystyle{1\over p^{u+z}}\over
1-\displaystyle{1\over p^{u+v}}}\right)
\prod_{p|d}
\left(1-\displaystyle{1\over p^{u-w+1}}
\over 1-\displaystyle{1\over p^{u-v-w+z+2}}\right)
\prod_{p^\beta\Vert c}(\cdots),&(5.6)
}
$$
where the last product is as in $(5.5)$.
\bigskip\noindent
{\bf 6.} Now let us turn to
$$
\eqalignno{
J_+^*(u,&v,w,z;g;d/c)={c^{-v-w-z}d^{-w}\over2\pi
i}\int_{(\eta_1)}g^*(s,w)(c^2d)^s
\sum_{(k,d)=1}{1\over k^{u+v+w+z-2s}}\cr
&\times\sum_f{1\over f^s}\sum_{\scr{h=1}\atop\scr{(h,ck)=1}}^{ck}
\zeta\left(v+w-s,{h\over ck}\right)
\zeta\left(w+z-s,-{\overline{dh}f\over
ck}\right)ds,&(6.1)
}
$$
where $(4.1)$ holds. On noting that $\Re(v+w-s)<0$, $\Re(w+z-s)<0$, we
appeal to the functional equation $(4.7)$. Then the last 
double sum is equal to
$$
\eqalignno{
&4{\Gamma(1+s-v-w)\Gamma(1+s-w-z)\over(2\pi)^{2+2s-v-2w-z}}
\sum_{f,m,n}m^{v+w-s-1}n^{w+z-1}(fn)^{-s}\cr
\times& \sum_{\scr{h=1}\atop\scr{(h,ck)=1}}^{ck}
\sin\left(\txt{1\over2}\pi(v+w-s)+2\pi {m\over ck}h\right)
\sin\left(\txt{1\over2}\pi(w+z-s)-2\pi{fn\over ck}\overline{dh}\right)\cr
=&2{\Gamma(1+s-v-w)\Gamma(1+s-w-z)\over(2\pi)^{2+2s-v-2w-z}}
\sum_{f,m,n}m^{v+w-s-1}n^{w+z-1}(fn)^{-s}\cr
\times&\left\{\cos\left(\txt{1\over2}\pi(v-z)\right)
S\!\left(m,\overline{d}fn;ck\right)-
\cos\left(\txt{1\over2}\pi(v+2w+z-2s)\right)
S\!\left(m,-\overline{d}fn;ck\right)\right\}\cr
=&2{\Gamma(1+s-v-w)\Gamma(1+s-w-z)\over(2\pi)^{2+2s-v-2w-z}}
\sum_{m,n}m^{v+w-s-1}n^{-s}\sigma_{w+z-1}(n)\cr
\times&\left[\cos\left(\txt{1\over2}\pi(v-z)\right)
S\!\left(m,\overline{d}n;ck\right)-
\cos\left(\txt{1\over2}\pi(v+2w+z-2s)\right)
S\!\left(m,-\overline{d}n;ck\right)\right],\quad&(6.2)
}
$$
where $S$ is the ordinary Kloosterman sum, and
$\sigma_\tau(n)=\sum_{\lambda|n}
\lambda^\tau$.
\medskip
Thus
$$
\eqalignno{
J_+^*&(u,v,w,z;g;d/c)={c^ud^{{1\over2}(u+v-w+z)}\over\pi i
(2\pi)^{u-w+1}}\cr
&\times\sum_{m,n}m^{{1\over2}(v+w-u-z-1)}
n^{-{1\over2}(u+v+w+z-1)}\sigma_{w+z-1}(n)
\sum_{(k,d)=1}{1\over ck\sqrt{d}}\cr
&\times\int_{(\eta_1)}
\Big[\cos\left(\txt{1\over2}\pi(v-z)\right)
S\!\left(m,\overline{d}n;ck\right)
-\cos\left(\txt{1\over2}\pi(v+2w+z-2s)\right)
S\!\left(m,-\overline{d}n;ck\right)\Big]\cr
&\times\Gamma(1+s-v-w)\Gamma(1+s-w-z)
g^*(s,w)
\left({2\pi\sqrt{mn}\over ck\sqrt{d}}\right)^{u+v+w+z-2s-1}
ds.&(6.3)
}
$$
We put
$$
\eqalignno{
&\tilde{g}_+(u,v,w,z;x)={1\over2\pi i}
\cos\left(\txt{1\over2}\pi(v-z)\right)\cr
&\hskip1cm\times\int_{(\eta_1)}
\Gamma(1+s-v-w)\Gamma(1+s-w-z)g^*(s,w)(x/2)^{u+v+w+z-2s-1}ds,\cr
&\tilde{g}_-(u,v,w,z;x)\cr
&=-{1\over2\pi i}\int_{(\eta_1)}
\cos\left(\txt{1\over2}\pi(v+2w+z-2s)\right)\cr
&\hskip1cm\times\Gamma(1+s-v-w)\Gamma(1+s-w-z)g^*(s,w)
(x/2)^{u+v+w+z-2s-1}ds,&(6.4)
}
$$
and
$$
\eqalignno{
Y_\pm(u,v,w,z;g;d/c;m,n)&=\sum_{(k,d)=1}{1\over ck\sqrt{d}}
S\!\left(m,\pm\overline{d}n;ck\right)
\tilde{g}_\pm\!\left(u,v,w,z;{4\pi\sqrt{mn}\over ck\sqrt{d}}\right)\cr
&=\sum_{(k,d)=1}{1\over ck\sqrt{d}}
S\!\left(n,\pm\overline{d}m;ck\right)
\tilde{g}_\pm\!\left(u,v,w,z;{4\pi\sqrt{mn}\over ck\sqrt{d}}\right).
&(6.5)
}
$$
We have
$$
J_+^*(u,v,w,z;g;d/c)=\left[K_++K_-\right](u,v,w,z;g;d/c),\eqno(6.6)
$$
with
$$
\eqalignno{
&K_{\pm}(u,v,w,z;g;d/c)=2{c^ud^{{1\over2}(u+v-w+z)}\over
(2\pi)^{u-w+1}}\cr
&\times\sum_{m,n}m^{{1\over2}(v+w-u-z-1)}
n^{-{1\over2}(u+v+w+z-1)}\sigma_{w+z-1}(n)
Y_\pm(u,v,w,z;g;d/c;m,n).&(6.7)
}
$$
\bigskip\noindent
{\bf 7.} We need to spectrally decompose the sums $Y_\pm$. To this end
we shall begin with some basic facts about a generic discrete subgroup
$\varGamma$ of $\r{PSL}(2,\B{R})$ and later proceed to the Kuznetsov sum
formula for the Hecke congruence subgroup $\varGamma_0(q)$. 
\medskip
Thus, let $\varGamma$ be a discrete subgroup of $\r{PSL}(2,\B{R})$
which has a fundamental domain of finite volume. We call $\f{a}$ a cusp
of $\varGamma$ if and only if there exists a 
$\sigma\in\varGamma$ such that
$\sigma$ is parabolic, i.e., $\r{Tr}(\sigma)=\pm 2$ and
$\sigma(\f{a})=\f{a}
\in \B{R}\cup{\infty}$. 
Let $\varGamma_{\f{a}}$ be $\{\sigma\in\varGamma:\sigma(\f{a})=\f{a}\}$,
i.e., the stabilizer of $\f{a}$. Then $\varGamma_{\f{a}}$ is cyclic,
so all elements in it are parabolic. Hence, there exits a $\sigma_\f{a}$
such that $\sigma_\f{a}(\infty)=\f{a}$ and
$\sigma_{\f{a}}^{-1}\varGamma_\f{a}\sigma_\f{a}=\varGamma_\infty
=[S]$ with $S=\left({1\atop}{1\atop1}\right)$. 
\medskip
The discussion below
depends on the choice of $\sigma_\f{a}$ which is not unique. 
If $\sigma_\f{a}'$ is another choice, then there exists a $b$
such that $\sigma_\f{a}'=\sigma_\f{a}S^b$. In fact, since
$\sigma_\f{a}^{-1}\varGamma_\f{a}\sigma_\f{a}=
{\sigma_\f{a}'}^{-1}\varGamma_\f{a}\sigma_\f{a}'$, we have $
\sigma_\f{a}S\sigma_\f{a}^{-1}=
\sigma_\f{a}'S^{\pm1}{\sigma_\f{a}'}^{-1}$ or $\sigma_\f{a}^{-1}
\sigma_\f{a}'S^{\pm1}=S\sigma_\f{a}^{-1}\sigma_\f{a}'$. On the other
hand $\sigma_\f{a}^{-1}\sigma_\f{a}'(\infty)=\infty$ implies that
$\sigma_\f{a}^{-1}\sigma_\f{a}'=\left({a\atop }{b\atop c}\right)$; and
$\left({a\atop }{b\atop c}\right)\left({1\atop }{\pm1\atop 1}\right)
=\left({1\atop }{1\atop 1}\right)\left({a\atop }{b\atop c}\right)$ yields
that $a=\pm c$, that is, $a=c=1$ and the assertion follows.
\medskip
Let $f$ be a $\varGamma$-automorphic form of weigh $2k$, with a positive
integer $k$; namely, for $\sigma=\left({a\atop c}{b\atop d}\right)
\in\varGamma$,
$$
\eqalignno{
f\left(\sigma(z)\right)&=(cz+d)^{2k}f(z)\cr
&=\jmath(\sigma,z)^{2k}f(z).&(7.1)
}
$$
The function $f(\sigma_\f{a}(z))
(\jmath(\sigma_\f{a},z))^{-2k}$ is of
period
$1$. In fact,
$$
\eqalignno{
f(\sigma_\f{a}S(z))(\jmath(\sigma_\f{a},S(z)))^{-2k}
&=f(\sigma_\f{a}S\sigma_\f{a}^{-1}\sigma_\f{a}(z))
(\jmath(\sigma_\f{a},S(z)))^{-2k}\cr
&=f(\sigma_\f{a}(z))(\jmath(\sigma_\f{a}S\sigma_\f{a}^{-1},
\sigma_\f{a}(z)))^{2k}(\jmath(\sigma_\f{a},S(z)))^{-2k}\cr
&=f(\sigma_\f{a}(z))\left[\jmath(\sigma_\f{a}S,z)
/\jmath(\sigma_\f{a},z)\right]^{2k}(\jmath(\sigma_\f{a},S(z)))^{-2k}\cr
&=f(\sigma_\f{a}(z))(\jmath(\sigma_\f{a},z))^{-2k}.&(7.2)
}
$$
Thus, if $f(\sigma_\f{a}(z))$ is regular near $\infty$, then the function
$f(\sigma_\f{a}(\log z/2\pi i))
(\jmath(\sigma_\f{a},\log z/2\pi i))^{-2k}$
is single valued and regular on a small disk centered and punctured at the
origin. Hence
$$
f(\sigma_\f{a}(z))(\jmath(\sigma_\f{a},z))^{-2k}=
\sum_n \varrho(n,\f{a})\exp(2\pi i nz),\eqno(7.3)
$$
which is called the Fourier expansion of $f$ around the cusp $\f{a}$.
\medskip
Note that this expansion 
depends on the choice of $\sigma_\f{a}$. In fact
If $\sigma_\f{a}'$ is another choice , then
$\sigma_\f{a}'=\sigma_\f{a}S^b$ with a $b$. We have 
$f(\sigma_\f{a}'(z))(\jmath(\sigma_{\f{a}'},z))^{-2k}
=f(\sigma_\f{a}(z+b))
(\jmath(\sigma_\f{a},z+b))^{-2k}$. That is, $\varrho(n,\f{a})$ is
multiplied by $\exp(2\pi inb)$.
\medskip
If $f$ is regular on the upper half plane ${\cal H}=\{z=x+iy:
-\infty<x<\infty, y>0\}$ and 
$\varrho(n,\f{a})=0$ for any
$n\le0$ and any $\f{a}$, then $f$ is termed a holomorphic cusp-form. Let 
${\cal S}_k(\varGamma)$ be the space of all cusp-forms of weight $2k$.
Then ${\cal S}_k(\varGamma)$ is a finite dimensional Hermitian space
with the Petersson inner product
$$
\langle f,g\rangle_k=\int_{\varGamma\backslash\cal H}
f(z)\overline{g(z)}y^{2k}d\mu(z),\quad d\mu(z)=dxdy/y^2.\eqno(7.4)
$$
We let $\left\{\psi_{j,k}(z),1\le j\le\vartheta(k)\right\}$ stand for an
orthonormal base of ${\cal S}_k(\varGamma)$.
\medskip

\medskip
\noindent
{\bf 8.} Let $k\ge2$. We introduce the Poincar\'e series 
$$
P_m(z,\f{a};k)=\sum_{\gamma\in\varGamma_\f{a}\backslash\varGamma}
(\jmath(\sigma_\f{a}^{-1}\gamma,z))^{-2k}\exp(2\pi i m\sigma_\f{a}^{-1}
\gamma(z)).\eqno(8.1)
$$
This is a holomorphic cusp form of weight $2k$ for any integer $m>0$. 
We shall confirm this claim, though we skip
the convergence issue, which causes no difficulty when $k\ge2$.
\medskip
First, each summand is a function
over
$\varGamma_\f{a}\backslash\varGamma$. In fact, if
$\varGamma_\f{a}\gamma=\varGamma_\f{a}\gamma'$, then $\varGamma_\infty
\sigma_\f{a}^{-1}\gamma=\varGamma_\infty\sigma_\f{a}^{-1}\gamma'$ and
$\sigma_\f{a}^{-1}\gamma(z)\equiv \sigma_\f{a}^{-1}\gamma'(z)\bmod\, 1$
as well as $\jmath(\sigma_\f{a}^{-1}\gamma,z)
=\jmath(\sigma_\f{a}^{-1}\gamma',z)$. Also the relation
$P_m(\gamma(z),\f{a};k)=(\jmath(\gamma,z))^{2k}P_m(z,\f{a};k)$ is obvious;
and $P_m(z,\f{a};k)$ is regular over $\cal H$. Thus, it remains to
consider the Fourier expansion at a given cusp $\f{b}$. 
We have
$$
\eqalignno{
P_m&(\sigma_\f{b}(z),\f{a};k)(\jmath(\sigma_\f{b},z))^{-2k}\cr
&=\sum_{\gamma\in\varGamma_\f{a}\backslash\varGamma}
(\jmath(\sigma_\f{b},z))^{-2k}(\jmath(\sigma_\f{a}^{-1}
\gamma,\sigma_\f{b}(z)))^{-2k}\exp(2\pi
im\sigma_\f{a}^{-1}\gamma\sigma_\f{b}(z)) \cr
&=\sum_{\gamma\in\varGamma_\f{a}\backslash\varGamma}
(\jmath(\sigma_\f{a}^{-1}\gamma\sigma_\f{b},z))^{-2k}\exp(2\pi
im\sigma_\f{a}^{-1}\gamma\sigma_\f{b}(z))\cr
&=\sum_{\gamma\in\varGamma_\f{a}\backslash\varGamma}
{1\over(cz+d)^{2k}}\exp\left(2\pi i m{az+b\over cz+d}\right),\qquad
\sigma_\f{a}^{-1}\gamma\sigma_\f{b}=\left({a\atop c}{b\atop d}\right).
&(8.2)
}
$$
If $c=0$, then $\sigma_\f{a}^{-1}\gamma\sigma_\f{b}(\infty)=\infty$ or
$\gamma(\f{b})=\f{a}$, that is, $\f{a}\equiv\f{b}\bmod\,\varGamma$ as well
as $\gamma\sigma_\f{b}=\sigma_\f{a}S^b$. Moreover, if
$\gamma'\sigma_\f{b}=\sigma_\f{a}S^{b'}$, then $\gamma'(\f{b})=\f{a}
=\gamma(\f{b})$, that is, $\gamma'\gamma^{-1}\in\varGamma_\f{a}$
or $\varGamma_\f{a}\gamma=\varGamma_\f{a}\gamma'$. Hence
$$
\sum_{\scr{\gamma\in\varGamma_\f{a}\backslash\varGamma}
\atop\scr{c=0}}=\delta_{\f{a},\f{b}}\exp(2\pi im(z+b)).\eqno(8.3)
$$
As to the remaining part, we have
$$
\sum_{\scr{\gamma\in\varGamma_\f{a}\backslash\varGamma}
\atop\scr{c\ne0}}=\sum_{\scr{\gamma\in\varGamma_\f{a}\backslash\varGamma}
\atop\scr{c\ne0}}{1\over(cz+d)^{2k}}
\exp\left(2\pi im{a\over c}-2\pi im{1\over c(cz+d)}\right).\eqno(8.4)
$$
We observe that if $\sigma_\f{a}^{-1}\gamma\sigma_\f{b}
=\left({a\atop c}{b\atop d}\right)$ appears in the right side, then
$\sigma_\f{a}^{-1}\gamma\sigma_\f{b}S^n
\sigma_\f{b}^{-1}\sigma_\f{b}=\left({a\atop c}{b+an\atop d+cn}\right)$
does for all $n\in\B{Z}$. In fact, $\gamma\sigma_\f{b}S^n\sigma_\f{b}^{-1}
\in\varGamma$ and thus 
$\varGamma_\f{a}\gamma\sigma_\f{b}S^n\sigma_\f{b}^{-1}$ is an
element of $\varGamma_\f{a}\backslash\varGamma$. Moreover, if
$\varGamma_\f{a}\gamma\sigma_\f{b}S^m\sigma_\f{b}^{-1}=
\varGamma_\f{a}\gamma\sigma_\f{b}S^n\sigma_\f{b}^{-1}$, then
$\sigma_\f{a}\varGamma_\infty\sigma_\f{a}^{-1}
\gamma\sigma_\f{b}S^m\sigma_\f{b}^{-1}=
\sigma_\f{a}\varGamma_\infty\sigma_\f{a}^{-1}
\gamma\sigma_\f{b}S^n\sigma_\f{b}^{-1}$ or 
$\sigma_\f{a}^{-1}\gamma\sigma_\f{b}S^m=S^l
\sigma_\f{a}^{-1}\gamma\sigma_\f{b}S^n$. This means that
$$
\left({a\atop c}{b+am\atop d+cm}\right)=
\left({a+cl\atop c}{b+an+(d+cn)l\atop d+cn}\right);\eqno(8.5)
$$
and we get $l=0, m=n$, which confirms our claim. On the other hand, since
$\{\gamma\sigma_\f{b}S^n\sigma_\f{b}^{-1}:n\in\B{Z}\}=
\gamma\varGamma_\f{b}$, we
should classify the summands in $(8.4)$ according to the double coset
decomposition $\varGamma_\f{a}\backslash\varGamma/\varGamma_\f{b}$,
which naturally we could have introduced already at $(8.2)$.
\medskip
We have  thus
$$
\sum_{\scr{\gamma\in\varGamma_\f{a}\backslash\varGamma}
\atop\scr{c\ne0}}=
\sum_{\scr{\gamma\in\varGamma_\f{a}\backslash\varGamma/\varGamma_\f{b}}
\atop\scr{c\ne0}}\sum_n{1\over(c(z+n)+d))^{2k}}
\exp\left(2\pi im{a\over c}-2\pi im{1\over c(c(z+n)+d)}\right).\eqno(8.6)
$$
More explicitly, we have the relation 
$\gamma\in\varGamma_\f{a}\backslash\varGamma/\varGamma_\f{b}$ is
equivalent to $\left({a\atop c}{b\atop d}\right)\in\varGamma_\infty
\backslash\sigma_\f{a}^{-1}\varGamma\sigma_\f{b}/\varGamma_\infty$.
With this one may proceed just in the same way as the case
of the full modular group and get
$$
\eqalignno{
&P_m(\sigma_\f{b}(z),\f{a};k)(j(\sigma_\f{b},z))^{-2k}=\delta_{\f{a},
\f{b}}\exp(2\pi im(z+b))\cr
+&2\pi(-1)^k\sum_{n>0}\left\{\sum_{c>0}{1\over c}S(m,n;c;\f{a},\f{b})
\left({n\over m}\right)^{k-{1\over2}}J_{2k-1}
\left(4\pi{\sqrt{mn}\over c}\right)\right\}\exp(2\pi inz).&(8.7)
}
$$
Here
$$
S(m,n;c;\f{a},\f{b})=\sum_\gamma\exp\left(2\pi i(am+dn)/c\right)
\eqno(8.8)
$$
is a Kloosterman sum associated with $\varGamma$,
where $\gamma$ runs over the representatives of
$\varGamma_\f{a}
\backslash\varGamma/\varGamma_\f{b}$
with the same $c$ in the sense remarked after $(8.6)$. The expression
$(8.8)$ and the constant $b$ in $(8.7)$ depend of course on the choice of
$\sigma_\f{a}$,
$\sigma_\f{b}$. 
\medskip
The last summands are functions
on $\varGamma_\f{a}\backslash\varGamma/\varGamma_\f{b}$. In fact, let 
$\varGamma_\f{a}\gamma\varGamma_\f{b}
=\varGamma_\f{a}\gamma'\varGamma_\f{b}$. Then $\sigma_\f{a}
\varGamma_\infty\sigma_\f{a}^{-1}\gamma\sigma_\f{b}\varGamma_\infty
\sigma_\f{b}^{-1}\ni\gamma'$ or $\varGamma_\infty\sigma_\f{a}^{-1}
\gamma\sigma_\f{b}=\sigma_\f{a}^{-1}\gamma'\sigma_\f{b}\varGamma_\infty$,
which means that there exist two integers $l,l'$ such that
$S^l\left({a\atop c}{b\atop d}\right)
=\left({a'\atop c'}{b'\atop d'}\right)S^{l'}$. Hence $c=c'$ and
$a\equiv a'$, $d\equiv d' \bmod\, c$. Also, for each $c>0$ there are
at most finitely many double cosets having $c$ as the lower-left 
element; otherwise the convergence would be violated. 
\medskip
On the assumption that there exists a $c_0>0$ such that for
any non-zero integers $m,n$ and any pair of cusps $\f{a},\f{b}$
$$
\sum_c {1\over c^{2k}}|S(m,n;c;\f{a},\f{b})|\ll (mn)^{c_0},\eqno(8.9)
$$
we have
$$
P_m(\sigma_\f{b}(\infty),\f{a};k)=0,\eqno(8.10)
$$
implying that $P_m$ is a holomorphic cusp form of weight $2k$.
\bigskip\noindent
{\bf 9.} We consider the spectral decomposition
$$
\langle P_m(\cdot,\f{a};k),P_n(\cdot,\f{b};k)\rangle_k
=\sum_{j=1}^{\vartheta(k)}
\langle P_m(\cdot,\f{a};k),\psi_{j,k}\rangle_k
\overline{\langle P_n(\cdot,\f{b};k),\psi_{j,k}\rangle_k}\,.\eqno(9.1)
$$
The left side is
$$
\eqalignno{
&\sum_{\gamma\in\varGamma_\f{b}\backslash\varGamma}\int_{\varGamma
\backslash\cal H}P_m(z,\f{a};k)
\overline{(\jmath(\sigma_\f{b}^{-1}\gamma,z))^{-2k}
\exp(2\pi in\sigma_\f{b}^{-1}\gamma(z))}y^{2k}d\mu(z)\cr
=&\sum_{\gamma\in\varGamma_\f{b}\backslash\varGamma}
\int_{\sigma_\f{b}^{-1}\gamma(\varGamma
\backslash\cal H)}P_m(\gamma^{-1}\sigma_\f{b}(z),\f{a};k)
\overline{(\jmath(\sigma_\f{b}^{-1}
\gamma,\gamma^{-1}\sigma_\f{b}(z)))^{-2k}
\exp(2\pi inz)}{y^{2k}d\mu(z)
\over|\jmath(\gamma^{-1}\sigma_\f{b},z)|^{4k}}\cr
=&\sum_{\gamma\in\varGamma_\f{b}\backslash\varGamma}
\int_{\sigma_\f{b}^{-1}\gamma(\varGamma
\backslash\cal H)}P_m(\gamma^{-1}\sigma_\f{b}(z),\f{a};k)
(\jmath(\gamma^{-1}\sigma_\f{b},z)))^{-2k}
\exp(-2\pi in\bar{z})y^{2k}d\mu(z)\cr
=&\sum_{\gamma\in\varGamma_\f{b}\backslash\varGamma}
\int_{\sigma_\f{b}^{-1}\gamma(\varGamma
\backslash\cal H)}P_m(\sigma_\f{b}(z),\f{a};k)
(\jmath(\gamma^{-1},\sigma_\f{b}(z))^{2k}
(\jmath(\gamma^{-1}\sigma_\f{b},z)))^{-2k}
\exp(-2\pi in\bar{z})y^{2k}d\mu(z)\cr
=&\sum_{\gamma\in\varGamma_\f{b}\backslash\varGamma}
\int_{\sigma_\f{b}^{-1}\gamma(\varGamma
\backslash\cal H)}P_m(\sigma_\f{b}(z),\f{a};k)
(\jmath(\sigma_\f{b},z)))^{-2k}
\exp(-2\pi in\bar{z})y^{2k}d\mu(z)\cr
=&\int_{\sigma_\f{b}^{-1}
\bigcup_{\gamma\in\varGamma_\f{b}\backslash\varGamma}
\gamma(\varGamma
\backslash\cal H)}P_m(\sigma_\f{b}(z),\f{a};k)
(\jmath(\sigma_\f{b},z)))^{-2k}
\exp(-2\pi in\bar{z})y^{2k}d\mu(z)\cr
=&\int_0^\infty\int_0^1
P_m(\sigma_\f{b},\f{a};k)
(\jmath(\sigma_\f{b},z)))^{-2k}
\exp(-2\pi in\bar{z})y^{2k-2}dxdy\cr
=&2\pi\Gamma(2k-1)(4\pi\sqrt{mn})^{1-2k}
\Big\{{1\over2\pi}\delta_{\f{a},\f{b}}\delta_{m,n}\exp(2\pi inb)\cr
&\hskip2cm+(-1)^k\sum_c{1\over
c}S(m,n;c;\f{a},\f{b})J_{2k-1}\left(4\pi{\sqrt{mn}\over c}\right)
\Big\},
&(9.2)
}
$$
where we have used that $\sigma_\f{b}^{-1}
\bigcup_{\gamma\in\varGamma_\f{b}\backslash\varGamma}
\gamma(\varGamma\backslash{\cal H})=\sigma_\f{b}^{-1}(\varGamma_\f{b}
\backslash\cal H)=\varGamma_\infty\backslash{\cal H}$; in fact,
since\break
$\sigma_\f{b}\varGamma_\infty
\sigma_\f{b}^{-1}(\varGamma_\f{b}\backslash {\cal H})={\cal H}$, we have
$\varGamma_\infty
\sigma_\f{b}^{-1}(\varGamma_\f{b}\backslash {\cal H})={\cal H}$.
\medskip
On the other hand, we have in much the same way
$$
\eqalignno{
\langle P_m(\cdot,\f{a};k),\psi_{j,k}\rangle_k
=&\int_0^\infty\int_0^1 \exp(2\pi imz)
\overline{\psi_{j,k}(\sigma_\f{a}(z))
\jmath(\sigma_\f{a},z)^{-2k}}d\mu(z)\cr
=&\Gamma(2k-1)(4\pi m)^{1-2k}\overline{\varrho_{j,k}(m,\f{a})},&(9.3)
}
$$
where we have put, following $(7.3)$,
$$
\psi_{j,k}(\sigma_\f{a}(z))\jmath(\sigma_\f{a},z)^{-2k}=\sum_{n>0}
\varrho_{j,k}(n,\f{a})\exp(2\pi inz).\eqno(9.4)
$$
\medskip
Hence we have obtained the Petersson Formula: 
\medskip
\noindent
{\bf Lemma 1.} {\it
For $k\ge2$
$$
\eqalignno{
&{1\over2\pi}{\Gamma(2k-1)\over(4\pi\sqrt{mn})^{2k-1}}
\sum_{j=1}^{\vartheta(k)}\overline{\varrho_{j,k}(m,\f{a})}
\varrho_{j,k}(n,\f{b})\cr
=&{1\over2\pi}\delta_{\f{a},\f{b}}\delta_{m,n}\exp(2\pi inb)+
(-1)^k\sum_c{1\over
c}S(m,n;c;\f{a},\f{b})J_{2k-1}\left(4\pi{\sqrt{mn}\over c}\right),&(9.5)
}
$$
provided $\varGamma$ satisfies $(8.9)$\/}.
\medskip
The case $k=1$ can also be treated in much the same way as is done with
the full modular group (see [11, pp.\ 52--54]), 
excepting that $(8.9)$ should be replaced by the assumption:
There exists a constant $\tau<2$ such that for any 
non-zero integers $m,n$ and for any pair of cusps $\f{a},\f{b}$
$$
\sum_c{1\over
c^\tau}|S(m,n;c;\f{a},\f{b})|\ll |mn|^{c_0}.\eqno(9.6)
$$
On this $(9.5)$ holds for all $k\ge1$. 
\bigskip\noindent
{\bf 10.} We turn to real analytic cusp forms. The procedure is
similar to the holomorphic case and also to the full modular situation,
and we can be brief.
\medskip
Let $f$ be a real analytic cusp form of weight zero with respect to
$\varGamma$ so that $f(\gamma(z))=f(z)$ for all $\gamma\in\varGamma$,
and $\Delta f=\nu f$ with $\Delta=-y^2(\partial_x^2+\partial_y^2)$.
Since $f(\sigma_\f{a}(z))$ is of period one, we have the Fourier
expansion
$$
f(\sigma_\f{a}(z))=\sum_n \varrho(n,\f{a};y)\exp(2\pi inx).\eqno(10.1)
$$
We require that 
$$
\hbox{$\displaystyle\lim_{z\to\infty}f(\sigma_\f{a}(z))=0$ for any 
$\f{a}$, and}\;
\int_{\varGamma\backslash {\cal H}}|f|^2d\mu(z)<\infty.\eqno(10.2)
$$
We have then
$$
f(\sigma_\f{a}(z))=y^{1\over2}\sum_{n\ne0}\varrho(n,\f{a})K_{i\kappa}
(2\pi|n|y)\exp(2\pi inx),\quad \nu=\kappa^2+\txt{1\over4}.\eqno(10.3)
$$
\medskip
One may consider more generally the decomposition of the
space $L^2(\varGamma\backslash G)$, $G=\r{PSL}(2,
\B{R})$ into irreducible subspaces and
appeal to the theory of representations of the Lie
group $G$. This will allow us to deal with all cusp forms of various
weights in a unified fashion. However, here we shall rather follow
the argument due to Kuznetsov and others.
\medskip
Thus, let us introduce the Poincar\'e series of the Selberg type
$$
U_m(z,\f{a};s)=\sum_{\gamma\in\varGamma_\f{a}\backslash\varGamma}
\left(\r{Im}\sigma_\f{a}^{-1}\gamma(z)\right)^s\exp(2\pi im
\sigma_\f{a}^{-1}\gamma(z)),\eqno(10.4)
$$
and the Eisenstein series $E(z,\f{a};s)=U_0(z;\f{a};s)$, associated
with the cusp $\f{a}$. Arguing
as in Section 8, we have the Fourier expansion
$$
\eqalignno{
U_m(\sigma_\f{b}(z),\f{a};s)=\delta_{\f{a},\f{b}}y^s\exp(2\pi im&(z+b))
+y^{1-s}\sum_n\exp(2\pi inx)\sum_c{1\over c^{2s}}S(m,n;c;\f{a},\f{b})\cr
&\times\int_{-\infty}^\infty\exp\left(-2\pi iny\xi-
{2\pi m\over c^2y(1-i\xi)}\right){d\xi\over(1+\xi^2)^s}.&(10.5)
}
$$
On the assumption $(9.6)$, $U_m(\sigma_\f{b}(z),\f{a};s)$ is regular
for $\Re s>\tau/2$, and also  $U_m(\sigma_\f{b}(z),\f{a};s)\ll
y^{1-\Re s}$ as $y\to\infty$. In particular, $U_m(z,\f{a};s)\in
L^2(\varGamma\backslash{\cal H})$ if $\Re s>\tau/2$. Also we have
$$
\eqalignno{
E(\sigma_\f{b}(z),\f{a};s)&=\delta_{\f{a},\f{b}}y^s+\sqrt{\pi}y^{1-s}
{\Gamma\left(s-{1\over2}\right)\over\Gamma(s)}e_0(s;\f{a},\f{b})\cr
&+{2\pi^s\over\Gamma(s)}y^{1\over2}\sum_{n\ne0}|n|^{s-{1\over2}}
e_n(s;\f{a},\f{b})K_{s-{1\over2}}(2\pi|n|y)\exp(2\pi inx),&(10.6)
}
$$ 
with 
$$
e_n(s;\f{a},\f{b})=\sum_c{1\over c^{2s}}S(0,n;\f{a},\f{b}).\eqno(10.7)
$$
It can be shown that $E(\sigma_\f{b}(z),\f{a};s)$ is meromorphic for
all $s$. Moreover, in the case of congruence subgroups,
$E(\sigma_\f{b}(z),\f{a};s)$ is regular for $\Re s\ge{1\over2}$ except
for a simple pole at $s=1$.
\medskip
Let $\{\psi_j:j\ge1\}$ be a complete orthonormal base of 
the cuspidal subspace of $L^2(\varGamma\backslash{\cal H})$ such
that $\Delta\psi_j=\nu_j\psi_j$ with $\nu_j=\kappa_j^2+{1\over4}$, and
$$
\psi_j(\sigma_\f{a}(z))=y^{1\over2}\sum_{n\ne0}\varrho_j(n,\f{a})
K_{i\kappa_j}(2\pi|n|y)\exp(2\pi inx).\eqno(10.8)
$$
We put also $\psi_0\equiv (\hbox{\rm volume of $\varGamma\backslash{\cal
H})$}^{-1/2}$. We suppose that $\varGamma$ is such that no $E(z,\f{a};s)$
has poles in the interval $\left({1\over2},1\right)$. Then we have the
spectral expansion: For any pair $f,g\in L^2(\varGamma\backslash
{\cal H})$, it holds that
$$
\langle f,g\rangle=\sum_{j=0}^\infty \langle f,\psi_j\rangle
\overline{\langle g,\psi_j\rangle}+{1\over4\pi}\sum_\f{c}
\int_{-\infty}^\infty {\cal E}(r,\f{c};f)
\overline{{\cal E}(r,\f{c};g)}dr,\eqno(10.9)
$$
where $\langle\cdot,\cdot\rangle=\langle\cdot,\cdot\rangle_0$ and
$$
{\cal E}(r,\f{c};f)=\int_{\varGamma\backslash{\cal H}}
f(z)\overline{E\left(z,\f{c};\txt{1\over2}+ir\right)}d\mu(z).\eqno(10.10)
$$
\bigskip\noindent
{\bf 11.} We collect here analogues of Bruggeman's and Kuznetsov's
formulas: On the basic assumption $(9.6)$ we have:
\medskip
\noindent
{\bf Lemma 2.} {\it Uniformly for any $n\ne0$ and $\f{a}$,
$$
\sum_{\kappa_j\le K}{|\varrho_j(n,\f{a})|^2\over\cosh\pi\kappa_j}
+\sum_\f{c}\int_{-K}^K\left|e_n
\left(\txt{1\over2}+ir;\f{c},\f{a}\right)\right|^2dr\ll K^2+
|n|^{c_1},\eqno(11.1)
$$
where $c_1$ depends on $\tau$, $c_0$ in $(9.6)$. In particular, we
have the bound
$$
\varrho_j(n,\f{a})\ll(\kappa_j+|n|^{{1\over2}c_1})
\exp\left(\txt{1\over2}\pi\kappa_j\right).\eqno(11.2)
$$
\/}
\medskip
\noindent
{\bf Lemma 3.} {\it
Let $h(r)$ be even, regular and of fast decay on the strip $|\r{Im}
r|<{1\over2}+\eta$ with an $\eta>0$.
Then it holds that for any $m,n>0$ and $\f{a},\f{b}$
$$
\eqalignno{
&\sum_{j=1}^\infty{\overline{\varrho_j(m,\f{a})}\varrho_j(\pm n,\f{b})
\over\cosh\pi\kappa_j}h(\kappa_j)+
{1\over\pi}\sum_\f{c}\int_{-\infty}^\infty (n/m)^{ir}
\overline{e_m\left(\txt{1\over2}+ir;\f{c},\f{a}\right)}
e_n\left(\txt{1\over2}+ir;\f{c},\f{b}\right)h(r)dr\cr
&={1\over\pi^2}\delta_{\f{a},\f{b}}\delta_{m,\pm n}\exp(2\pi imb)
\int_{-\infty}^\infty r\tanh(\pi r)h(r)dr
+\sum_c{1\over c}S(m,\pm n;c;\f{a},\f{b})h_\pm\!\left(4\pi{\sqrt{mn}\over
c}\right),\quad\qquad&(11.3) 
}
$$
where $\f{c}$ runs over all inequivalent cusps, and
$$
h_+(x)={2i\over\pi}\int_{-\infty}^\infty{rh(r)\over\cosh\pi r}
J_{2ir}(x)dr,\quad
h_-(x)={4\over\pi^2}\int_{-\infty}^\infty rh(r)\sinh(\pi r)
K_{2ir}(x)dr.\eqno(11.4)
$$
}
\medskip
\noindent
{\bf Lemma 4.} {\it
Let $\varphi$ be smooth and of fast decay over the positive real axis.
Then we have, for any $m,n>0$ and $\f{a},\f{b}$,
$$
\eqalignno{
\sum_c{1\over c}&S(m,\pm n;c;\f{a},\f{b})
\varphi\!\left(4\pi{\sqrt{mn}\over c}\right)=\sum_{j=1}^\infty
{\overline{\varrho_j(m,\f{a})}\varrho_j(\pm n,\f{b})
\over\cosh\pi\kappa_j}\hat{\varphi}_\pm(\kappa_j)\cr
&+{1\pm1\over4\pi(4\pi\sqrt{mn})^{2k-1}}\sum_{k=1}^\infty
\Gamma(2k)\hat\varphi_+\!\left(\left(\txt{1\over2}-2k\right)i\right)
\sum_{j=1}^{\vartheta(k)}\overline{\varrho_{j,k}(m,\f{a})}
\varrho_{j,k}(n,\f{b})\cr
&+{1\over\pi}\sum_\f{c}\int_{-\infty}^\infty (n/m)^{ir}
\overline{e_m\left(\txt{1\over2}+ir;\f{c},\f{a}\right)}
e_n\left(\txt{1\over2}+ir;\f{c},\f{b}\right)\hat{\varphi}_\pm(r)dr,&(11.5)
}
$$
where
$$
\eqalignno{
\hat{\varphi}_+(r)&={\pi i\over2\sinh\pi r}\int_0^\infty
\left\{J_{2ir}(x)-J_{-2ir}(x)\right\}\varphi(x){dx\over x},\cr
\hat{\varphi}_-(r)&=2\cosh(\pi r)\int_0^\infty K_{2ir}(x)\varphi(x)
{dx\over x}.&(11.6)
}
$$
}
\bigskip\noindent
{\bf 12.} With this, we shall consider the specialization $\varGamma
=\varGamma_0(q)$. Our discussion overlaps, to a certain extent, with that
developed in [3]; however, the present work can be read independently of
it. In this section we shall  fix a representative set of all cusps
inequivalent
$\bmod\,\varGamma_0(q)$.
\medskip
We introduce 
$V=\left\{\left({1\atop n}\,{\atop 1}\right):n\in\B{Z}\right\}$ the
stabilizer of the point $0$ in $\varGamma_0(1)$
and the double coset decomposition
$$
\varGamma_0(1)=\bigcup_\f{a}\varGamma_0(q)\gamma_\f{a}V,\eqno(12.1)
$$
where the symbol $\f{a}$ is to be regarded temporarily as to be just a
label. We begin with a particular $\gamma_\f{a}$, and transform it to
a matrix suitable for our purpose. We thus look into the product
$$
\left(\matrix{a&b\cr cq& d}\right)
\left(\matrix{e&f\cr g&h}\right)\left(\matrix{1&\cr n&1}\right)=
\left(\matrix{*&*\cr *&k}\right)\left(\matrix{1&\cr n&1}\right),
\eqno(12.2)
$$
where the middle matrix on the left side corresponds 
to $\gamma_\f{a}$. It is to be observed that $g$ is fixed $\bmod\, h$,
because of the action of $V$. We assume that $h\ne0$. We have
$k=cfq+dh$, and we claim  that this can be made equal to $(q,h)$. In fact
$c(fq/(q,h))+d(h/(q,h))=1$ is soluble in $c$ and $d$, 
for $(fq,h)=(q,h)$; then
$d\equiv\overline{h/(q,h)}\bmod\, fq/(q,h)$, and $d$ can be a prime large
enough so that $(d,q)=1$, and thus $(d,cq)=1$. 
With such a $d$ we may choose $a,b$ to satisfy
$ad-bcq=1$, which confirms our claim. On the other hand, if $h=0$, then
it suffices to put $c=\sgn(f), d=1$. Thus we may suppose that
$\gamma_\f{a}=\left({*\atop *}{*\atop w}\right)$ with $w|q$; that is,
each coset in $(12.1)$ contains elements of this property. 
\medskip
We then apply $(12.1)$ to the point $0$, getting
$$
\B{Q}\cup\{\infty\}=\bigcup_\f{a}\varGamma_0(q)\gamma_\f{a}(0).\eqno(12.3)
$$
This means that $\{\gamma_\f{a}(0):\f{a}\}$, with the current definition
of $\f{a}$, is the full set of inequivalent cusps $\bmod\,\varGamma_0(q)$.
In fact, that $\varGamma_0(q)\gamma_\f{a}(0)\ni\gamma_{\f{a}'}(0)$ 
implies readily that
$\varGamma_0(q)\gamma_\f{a}V=\varGamma_0(q)\gamma_{\f{a}'}V$; and
the stabilizer in $\varGamma_0(q)$ of 
$\gamma_\f{a}(0)$ is $\gamma_\f{a}V_{q/w}\gamma_\f{a}^{-1}$ with 
$V_d=\left\{\left({1\atop n}{\atop 1}\right): d|n\right\}$, provided
$\gamma_\f{a}=\left({*\atop *}{*\atop w}\right)$.
The labels $\{\f{a}\}$ indeed coincide with their former designation. 
Also, it should be noted that the element $w$ is unique to each double
coset, which can be proved by considering the relation
$\varGamma_0(q)\left({*\atop *}{*\atop
w}\right)V=\varGamma_0(q)\left({*\atop *}{*\atop w'}\!\right)V$
with respect to either $\bmod\, w$ or $\bmod\, w'$, getting $w|w'$ and
$w'|w$, respectively. Namely, if $w\ne w'$, then 
$\varGamma_0(q)\left({*\atop *}{*\atop
w}\right)V\cap\varGamma_0(q)\left({*\atop *}{*\atop w'}\!\right)V
=\emptyset$.
\medskip
Hence, it remains to see when the relation
$$
\varGamma_0(q)\left(\matrix{e&f\cr g&w}\right)V=
\varGamma_0(q)\left(\matrix{e'&f'\cr g'&w}\right)V\eqno(12.4)
$$
holds, where the two matrices are in $\varGamma_0(1)$ with $w|q$ and
$(gg',w)=1$. Thus, we have
$$
\eqalignno{
&\left(\matrix{a&b\cr c&d}\right)\left(\matrix{e&f\cr g&w}\right)=
\left(\matrix{e'&f'\cr g'&w}\right)\left(\matrix{1&\cr n&1}\right)
\quad \hbox{\rm with $q|c$}\cr\cr
\Longleftrightarrow\quad&\left(\matrix{a&b\cr c&d}\right)=
\left(\matrix{e'+nf'&f'\cr g'+nw&w}\right)\left(\matrix{w&-f\cr
-g&e}\right)\cr\cr
\Longleftrightarrow\quad& c=w(g'+nw)-gw=w(g'-g+nw)\cr\cr
\Longleftrightarrow\quad& w(g'-g+nw)\equiv 0\bmod\, q\cr
\Longleftrightarrow\quad& g'-g+nw\equiv0\bmod\, q/w\cr\cr
\Longleftrightarrow\quad& g'\equiv g \bmod\, (w,q/w).&(12.5)
}
$$
Hence
$$
(1.24)\quad\Longleftrightarrow\quad 
\hbox{$(gg',w)=1$ and $g\equiv g'\bmod\, (w,q/w)$}.\eqno(12.6)
$$
Namely, when $\gamma_\f{a}$ varies with $w$ fixed, then $g$ and thus $f$
runs over the complete residue classes $\bmod\, (w,q/w)$ while 
satisfying $(w,f)=1$. If $(u,(w,q/w))=1$, then obviously there exists an
$f$ such that $u\equiv f\bmod\, (w,q/w)$ and $(w,f)=1$. 
\medskip
Collecting the above, we have
\medskip
\noindent
{\bf Lemma 5.} {\it
A complete representative set of cusps inequivalent $\bmod\,
\varGamma_0(q)$ is given by
$$
\left\{{u\over w}:\, w|q,\,(u,w)=1,\, u\bmod\, (w,q/w)\right\},
\eqno(12.7)
$$
whose cardinality is
$$
\sum_{w|q}\varphi((w,q/w)). \eqno(12.8)
$$
}
\bigskip\noindent
{\bf 13.} Let us fix the stabilizers of those cusps given in $(12.7)$.
To this end we note first that if $\f{a}\ne\infty$ is a cusp of a discrete
group $\varGamma$, then
$$
\varGamma_\f{a}=\varGamma\cap\left\{
\left(\matrix{1+\nu\f{a}&-\nu\f{a}^2\cr \nu
&1-\nu\f{a}}\right):\,
\nu\in\B{R}\right\}.\eqno(13.1)
$$
In fact, since $(a\f{a}+b)/(c\f{a}+d)=\f{a}$, $a+d=2$, we see
that $\f{a}=(1-d)/c$, and the assertion follows with $c=\nu$.
If $\f{a}=u/w$ with $w|q,\,(u,w)=1$, then
$$
\varGamma_{u/w}=\varGamma_0(q)\cap\left\{
\left(\matrix{\displaystyle 1+\nu{u\over w}&
\displaystyle-
\nu {u^2\over w^2}\cr\cr 
\nu &\displaystyle 1-\nu
{u\over w}}\right):\nu\in\B{R}\right\},\eqno(13.2)
$$
and thus $\nu\in\B{Z}$, $\nu\equiv0\bmod\,q$,
$\nu\equiv0\bmod\,w^2$; namely
$$
\varGamma_{u/w}=\left\{
\left(\matrix{\displaystyle 1+\nu{u\over w}&
\displaystyle-
\nu {u^2\over w^2}\cr\cr 
\nu &\displaystyle 1-\nu
{u\over w}}\right)
: \nu\equiv 0\bmod\, [w^2,q]\right\}.\eqno(13.3)
$$
\medskip
We write
$$
q=vw=(v,w)^2v^*w^*,\quad v^*={v\over(v,w)},
\;w^*={w\over(v,w)}.\eqno(13.4)
$$
We put
$$
\varpi_{u/w}=
\left(\matrix{u&
\displaystyle{u\overline{u}-1\over
w}\cr\cr w&\overline{u}}\right),\quad u\overline{u}\equiv1\bmod
w,\eqno(13.5)
$$
and
$$
\tau_{v^*}=\left(\matrix{\sqrt{v^*}&\cr&
\displaystyle{1\over\sqrt{v^*}}}\right).
\eqno(13.6)
$$
Obviously we have $\varpi_{u/w}(\infty)=u/w$. Moreover,
we have
$$
\eqalignno{
\varpi_{u/w}^{-1}&\left(\matrix{\displaystyle 1+\nu{u\over w}&
\displaystyle-
\nu {u^2\over w^2}\cr\cr 
\nu &\displaystyle 1-\nu
{u\over w}}\right)\varpi_{u/w}\cr
=&\left(\matrix{\overline{u}&
-\displaystyle{u\overline{u}-1\over
w}\cr\cr -w& u}\right)\left(\matrix{\displaystyle 1+\nu{u\over w}&
\displaystyle-
\nu {u^2\over w^2}\cr\cr 
\nu &\displaystyle 1-\nu
{u\over w}}\right)\left(\matrix{u&
\displaystyle{u\overline{u}-1\over
w}\cr\cr w&\overline{u}}\right)\cr
=&\left(\matrix{1&-\displaystyle{\nu\over
w^2}\cr\cr&1}\right)=\tau_{v^*}
\left(\matrix{1&-\displaystyle{\nu\over v^*
w^2}\cr\cr&1}\right)\tau_{v^*}^{-1}. &(13.7)
}
$$
Hence, on noting that $[w^2,q]=v^*w^2$, we get
$$
\varGamma_{u/w}=\varpi_{u/w}\tau_{v^*}\varGamma_\infty
\tau_{v^*}^{-1}\varpi_{u/w}^{-1}\,,\eqno(13.8)
$$
which is equivalent to
$$
\varGamma_{u/w}=\varpi_{u/w}\left[S^{v^*}\right]
\varpi_{u/w}^{-1}\,.\eqno(13.9)
$$
\bigskip\noindent
{\bf 14.} In the 
the special instance where $q=v_iw_i$ with $(v_i,w_i)=1$, 
we shall consider the structure of the double coset decomposition
$\varGamma_{1/w_1}\backslash
\varGamma_0(q)/\varGamma_{1/w_2}$ and associated Kloosterman sums. 
\medskip
To this end we put 
$$
\eqalignno{
\sigma_{1/w_i}&=\varpi_{1/w_i}\tau_{v_i}S^{-\overline{w}_i/v_i}\cr
&=\varpi_{1/w_i}S^{-\overline{w}_i}\tau_{v_i},
&(14.1)
}
$$
where $w_i\overline{w}_i\equiv 1\bmod v_i$. The choice of a particular
value of $\overline{w}_i$ is irrelevant to our discussion of the
Kloosterman sums, as we shall show later. Note that
$$
\varGamma_{1/w_i}=\sigma_{1/w_i}\varGamma_\infty
\sigma_{1/w_i}^{-1}\,,\eqno(14.2)
$$
as is implied by $(13.8)$.
\medskip
We shall prove that
$$
S^{\overline{w}_1}
\varpi_{1/w_1}^{-1}\varGamma_0(q)\varpi_{1/w_2}
S^{-\overline{w}_2}=\left\{\left(
\matrix{(v_1,w_2)k&(v_1,v_2)l\cr\cr(w_1,w_2)r&(w_1,v_2)s}\right)
\in \r{SL}(2,\B{Z}),\; k,l,r,s\in\B{Z}\right\}\eqno(14.3)
$$
(cf.\ [6, p.\ 534]; note that there $q$ is square-free but here not
assumed to be so). In fact, we have, by
$(13.5)$,
$$
\varpi_{1/w_i}S^{-\overline{w}_i}=
\left(\matrix{1&-\overline{w}_i\cr w_i&1-w_i\overline{w}_i}\right);
\eqno(14.4)
$$
thus for $\left({a\atop c}{b\atop d}\right)\in\varGamma_0(q)$
$$
\eqalignno{
&S^{\overline{w}_1}
\varpi_{1/w_1}^{-1}
\left(\matrix{a&b\cr c&d}\right)
\varpi_{1/w_2}S^{-\overline{w}_2}\cr
&\equiv\left(\matrix{0&*\cr*&*}\right)\left(\matrix{a&b\cr0&d}\right)
\left(\matrix{*&*\cr0&*}\right)=\left(\matrix{0&*\cr*&*}\right)\bmod
(v_1,w_2),\cr
&\equiv\left(\matrix{0&*\cr*&*}\right)\left(\matrix{a&b\cr0&d}\right)
\left(\matrix{*&*\cr*&0}\right)=\left(\matrix{*&0\cr*&*}\right)\bmod
(v_1,v_2),\cr
&\equiv\left(\matrix{*&*\cr0&*}\right)\left(\matrix{a&b\cr0&d}\right)
\left(\matrix{*&*\cr0&*}\right)=\left(\matrix{*&*\cr0&*}\right)\bmod
(w_1,w_2),\cr
&\equiv\left(\matrix{*&*\cr0&*}\right)\left(\matrix{a&b\cr0&d}\right)
\left(\matrix{*&*\cr*&0}\right)=\left(\matrix{*&*\cr*&0}\right)\bmod
(w_1,v_2). &(14.5)
}
$$
On the other hand, we have that
$$
\eqalignno{
&\varpi_{1/w_1}S^{-\overline{w}_1}\left(
\matrix{(v_1,w_2)k&(v_1,v_2)l\cr\cr(w_1,w_2)r&(w_1,v_2)s}\right)
S^{\overline{w}_2}\varpi_{1/w_2}^{-1}\cr
&\equiv\left(\matrix{*&*\cr*&0}\right)\left(\matrix{0&*\cr*&*}\right)
\left(\matrix{*&*\cr0&*}\right)=\left(\matrix{*&*\cr0&*}\right)\bmod
(v_1,w_2),\cr
&\equiv\left(\matrix{*&*\cr*&0}\right)\left(\matrix{*&0\cr*&*}\right)
\left(\matrix{0&*\cr*&*}\right)=\left(\matrix{*&*\cr0&*}\right)\bmod
(v_1,v_2),\cr
&\equiv\left(\matrix{*&*\cr0&*}\right)\left(\matrix{*&*\cr0&*}\right)
\left(\matrix{*&*\cr0&*}\right)=\left(\matrix{*&*\cr0&*}\right)\bmod
(w_1,w_2),\cr
&\equiv\left(\matrix{*&*\cr0&*}\right)\left(\matrix{*&*\cr*&0}\right)
\left(\matrix{0&*\cr*&*}\right)=\left(\matrix{*&*\cr0&*}\right)\bmod
(w_1,v_2), &(14.6)
}
$$
and that $(v_1,w_2)(v_1,v_2)(w_1,w_2)(w_1,v_2)=q$. This proves $(14.3)$.
\medskip
Hence, we have
$$
\eqalignno{
&\varGamma_{1/w_1}\!\backslash\varGamma_0(q)/\varGamma_{1/w_2}
=\sigma_{1/w_1}\varGamma_\infty\sigma_{1/w_1}^{-1}
\backslash\varGamma_0(q)/\sigma_{1/w_2}\varGamma_\infty
\sigma_{1/w_2}^{-1}\cr\cr
\Longleftrightarrow\quad&\varGamma_\infty\backslash\tau_{v_1}^{-1}
S^{\overline{w}_1}
\varpi_{1/w_1}^{-1}\varGamma_0(q)\varpi_{1/w_2}
S^{-\overline{w}_2}
\tau_{v_2}/\varGamma_\infty\cr\cr
\Longleftrightarrow\quad&\varGamma_\infty\backslash\tau_{v_1}^{-1}
\left\{\left(
\matrix{(v_1,w_2)k&(v_1,v_2)l\cr\cr(w_1,w_2)r&(w_1,v_2)s}\right)
\in\r{SL}(2,\B{Z}),\;k,l,r,s\in\B{Z}\right\}
\tau_{v_2}/\varGamma_\infty\cr\cr\cr
\Longleftrightarrow\quad&\varGamma_\infty\backslash
\left\{\left(
\matrix{(v_1,w_2)k\sqrt{v_2/v_1}&(v_1,v_2)l/\sqrt{v_1v_2}\cr
\cr(w_1,w_2)r\sqrt{v_1v_2}&(w_1,v_2)s\sqrt{v_1/v_2}}\right)
\right\}/\varGamma_\infty\cr\cr\cr
\Longleftrightarrow\quad
&\hbox{\rm classifying the solutions of
$(v_1,w_2)(w_1,v_2)sk-(w_1,w_2)(v_1,v_2)rl=1$}\cr
&\hbox{\rm according to $(v_1,w_2)k\sqrt{v_2/v_1},
(w_1,v_2)s\sqrt{v_1/v_2}\bmod (w_1,w_2)r\sqrt{v_1v_2}$};\cr
&\hbox{\rm note the remark after $(8.6)$}\cr\cr
\Longleftrightarrow\quad&\hbox{\rm the moduli of the Kloosterman sums
have the form $(w_1,w_2)r\sqrt{v_1v_2}$}\cr 
&\hbox{\rm with
$((v_1,w_2)(w_1,v_2),r)=1$ and}\cr
&(v_1,w_2)(w_1,v_2)sk\equiv1\bmod(w_1,w_2)(v_1,v_2)r\cr
&(v_1,w_2)k\bmod v_1(w_1,w_2)r \longleftrightarrow k\bmod 
(v_1,v_2)(w_1,w_2)r\cr
&(w_1,v_2)s\bmod v_2(w_1,w_2)r \longleftrightarrow s\bmod 
(v_1,v_2)(w_1,w_2)r\cr\cr
\Longleftrightarrow\quad& c=(w_1,w_2)r\sqrt{v_1v_2},\; 
((v_1,w_2)(w_1,v_2),r)=1,\cr
&S(m,n;c;1/w_1,1/w_2)=
\sum_{\scr{s,k \bmod
(v_1,v_2)(w_1,w_2)r}\atop\scr{(v_1,w_2)(w_1,v_2)sk\equiv1\bmod
(v_1,v_2)(w_1,w_2)r}}\exp\left({2\pi
i(km+ns)\over(v_1,v_2)(w_1,w_2)r}\right) \cr\cr
&\phantom{S(m,n;c;1/w_1,1/w_2)}=
S(\overline{(v_1,w_2)}m,\overline{(w_1,v_2)} n; 
(v_1,v_2)(w_1,w_2)r),&(14.7)
}
$$
where the last member is an ordinary Kloosterman sum. 
\medskip
It remains to show the irrelevance of the choice of values
of $\overline{w}_j$. In fact, if we replace $\overline{w}_j$ by
$\overline{w}_j+nv_j$, $n\in\B{Z}$, then the first
equivalence assertion in $(14.7)$ does not change, for we have
$\tau_{v_j}^{-1}S^{nv_j}\tau_{v_j}=S^n\in\varGamma_\infty$.
\medskip
In particular, we find that if $q=cd$, $(c,d)=1$, and $(r,d)=1$, then
$$
S(m,n;cr\sqrt{d};1/q,1/c)=S(m,n;cr\sqrt{d};\infty,1/c)
=S(m,\overline{d} n; cr),\eqno(14.8)
$$
on the specification $(14.1)$ of $\sigma_{1/q}$ and $\sigma_{1/c}$.
\bigskip\noindent
{\bf 15.} We still need to see if $(9.6)$ is satisfied by the generic
$\varGamma_0(q)$. Until very recently we had been
unable to locate any rigorous treatment of those
generalized Kloosterman
sums over $\varGamma_0(q)$ in literature, excepting [9] and [10] where the
case with $q$ square-free is explicitly discussed on the basis of
$(14.7)$. With this situation, R. Bruggeman
kindly provided us with a treatment [1] of the sums using a partly
adelic reasoning; and it is assured that $(9.6)$ indeed holds with
any $\varGamma_0(q)$. Here we shall prove the same with an
alternative elementary method; this section can be read independently of
[1]. 
\medskip
We shall first redefine the Kloosterman sums associated with
the two cusps $u_i/w_i$, $i=1,2$, which are in the set $(12.7)$,
by introducing the convention
$$
\sigma_{u_i/w_i}=\varpi_{u_i/w_i}\tau_{v_i^*},\eqno(15.1)
$$
with $v_i^*$ as in Section 13, which is effective within
this section only. Note that when $u_i=1$ this does not coincide with
$(14.1)$; when discussing the absolute values of generalized Kloosterman
sums, obviously no difference is caused. Also, it is expedient to use the
Bruhat decomposition; that is,  in the big cell of
${\r{PSL}(2,\B{R})}$ we have
$$
\eqalignno{
\left(\matrix{a&b\cr c&d}\right)&=
\left(\matrix{1&a/c\cr&1}\right)
\left(\matrix{&-1/c\cr c&}\right)
\left(\matrix{1&d/c\cr&1}\right)\cr
&=B[a,d;c],&(15.2)
}
$$
say. 
\medskip
With this, let $\varkappa_q$ be the characteristic function of the set
$\varGamma_0(q)\subset\r{PSL}(2,\B{R})$. Then  Kloosterman sums
associated with the two cusps $u_i/w_i$, $i=1,2$, have moduli
$c\sqrt{v_1^*v_2^*}$,
$c\in\B{N}$; and under $(15.1)$ we have that
$$
\eqalignno{
&S(m,n;c\sqrt{v_1^*v_2^*};u_1/w_1,u_2/w_2)\cr
=&\sum_{\scr{ad\equiv 1\bmod c}
\atop{\scr{a\bmod v_1^*c}
\atop{\scr{d\bmod v_2^*c}}}}\varkappa_q
\left(\varpi_{u_1/w_1}B[a,d;c]
\varpi_{u_2/w_2}^{-1}\right)
\exp\left(2\pi i\left({ma\over v_1^*c}+
{nd\over v_2^*c}\right)\right),&(15.3)
}
$$
where $a,c,d\in\B{Z}$. In fact, by $(13.8)$ we need to consider
the double coset decomposition
$$
\eqalignno{
&\varGamma_\infty\backslash\tau_{v_1^*}^{-1}
\varpi_{u_1/w_1}^{-1}\varGamma_0(q)
\varpi_{u_2/w_2}\tau_{v_2^*}/\varGamma_\infty\cr
=&\varGamma_\infty\backslash\tau_{v_1^*}^{-1}
\left\{B[a,d;c]:\varkappa_q(\varpi_{u_1/w_1}B[a,d,c]
\varpi_{u_2/w_2}^{-1})=1\right\}\tau_{v_2^*}/\varGamma_\infty\cr
=&\varGamma_\infty\backslash
\left\{B\left[a\sqrt{v_2^*/v_1^*},d\sqrt{v_1^*/v_2^*};
c\sqrt{v_1^*v_2^*}\right]:
\varkappa_q(\varpi_{u_1/w_1}^{-1}B[a,d,c]
\varpi_{u_2/w_2})=1\right\}/\varGamma_\infty,\quad&(15.4)
}
$$
where $B[a,d;c]\in\varGamma_0(1)$, since $\varpi_{u_1/w_1}^{-1}
\varGamma_0(q)\varpi_{u_2/w_2}\subset\varGamma_0(1)$. 
The expression $(15.3)$ readily follows. 
In passing, we note that
$$
|S(m,n;c\sqrt{v_1^*v_2^*};u_1/w_1,u_2/w_2)|\le
v_1^*v_2^*\varphi(c),\eqno(15.5)
$$
for the number of summands on the right of $(15.3)$ is less than or
equal to $v_1^*v_2^*\varphi(c)$. In fact, a unique $d\bmod\,c$
corresponds to each $a$, $(a,c)=1$, or  $v_2^*$ 
classes  $d\bmod\,v_2^*c$ to each of $v_1^*\varphi(c)$ classes
$a\bmod v_1^*c$ with $(a,c)=1$.
\medskip
We remark that $\varkappa_q
\left(\varpi_{u_1/w_1}B[a,d;c]\varpi_{u_2/w_2}^{-1}\right)$
is a function over  $a\bmod v_1^*c$ and $d\bmod v_2^*c$.
To see this, we use the relation
$$
\eqalignno{
&\varpi_{u_1/w_1}B[a+a',d+d';c]\varpi_{u_2/w_2}^{-1}\cr
&=\varpi_{u_1/w_1}\left(\matrix{1& a'/c\cr&
1}\right)\varpi_{u_1/w_1}^{-1}\cdot
\varpi_{u_1/w_1}B[a,d;c]
\varpi_{u_2/w_2}^{-1}\cdot\varpi_{u_2/w_2}
\left(\matrix{1&d'/c\cr&
1}\right)\varpi_{u_2/w_2}^{-1}; \qquad&(15.6)
}
$$
and $(13.9)$ gives that
$$
\eqalignno{
&\varpi_{u_1/w_1}\left(\matrix{1& a'/c\cr&
1}\right)\varpi_{u_1/w_1}^{-1}\in
\varGamma_{u_1/w_1}\subset\varGamma_0(q),\cr
&\varpi_{u_2/w_2}\left(\matrix{1& d'/c\cr&
1}\right)\varpi_{u_2/w_2}^{-1}\in
\varGamma_{u_2/w_2}\subset\varGamma_0(q),
&(15.7)
}
$$
provided $v_1^*|(a'/c)\in\B{Z}$, 
$v_2^*|(d'/c)\in\B{Z}$, which proves the
assertion.
\medskip
Next, we shall show that if $ad\equiv1\bmod c$, then
$$
\varkappa_q\left(\varpi_{u_1/w_1}B[a,d;c]
\varpi_{u_2/w_2}^{-1}\right)=
\varkappa_q\left(\varpi_{\overline{c^*}u_1/w_1}
B[a,d;c_0]
\varpi_{\overline{c^*}u_2/w_2}^{-1}\right),\eqno(15.8)
$$
where $c=c_0c^*$ with
$c_0=(c,q^\infty)$, and
$\overline{c^*}c^*\equiv 1\bmod q$; note that
$\overline{c^*}u_i/w_i$ are cusps of $\varGamma_0(q)$. In fact,
computing the lower-left element of
$\varpi_{u_1/w_1}B[a,d;c]\varpi_{u_2/w_2}^{-1}$, we see that the value
of the left side of $(15.8)$ equals $1$ if and only if
$$
\overline{u}_2(aw_1+c\overline{u}_1)\equiv
w_2\left(w_1(ad-1)/c+d\overline{u}_1\right)
\bmod q;\eqno(15.9)
$$
and this is equivalent to the congruence
$$
\overline{\overline{c^*}u_2}(aw_1+
c_0\overline{\overline{c^*}u_1}\,)\equiv
w_2\left(w_1(ad-1)/c_0+d\overline{\overline{c^*}u_1}\,\right)
\bmod q,\;\eqno(15.10)
$$
which immediately implies $(15.8)$.
\medskip
Hence we have
$$
\eqalignno{
&S(m,n;c\sqrt{v_1^*v_2^*};u_1/w_1,u_2/w_2)\cr
=&\sum_{\scr{ad\equiv 1\bmod c}
\atop{\scr{a\bmod v_1^*c}
\atop{\scr{d\bmod v_2^*c}}}}\varkappa_q
\left(\varpi_{\overline{c^*}u_1/w_1}B[a,d;c_0]
\varpi_{\overline{c^*}u_2/w_2}^{-1}\right)
\exp\left(2\pi i\left({ma\over v_1^*c}+
{nd\over v_2^*c}\right)\right).&(15.11)
}
$$
Here we have
$$
{1\over v^*_ic}\equiv{\widetilde{c_i^*}\over
v_i^*c_0}+{\widetilde{v_i^*c_0}
\over c^*}\;\bmod 1,\eqno(15.12)
$$
with $\widetilde{c_i^*}c^*\equiv1\bmod v_ic_0$,
$\widetilde{v_i^*c_0}v_i^*c_0\equiv1\bmod c^*$. Inserting
this into $(15.11)$, putting $a\equiv a_0\bmod v_1^*c_0$,
$a\equiv a^*\bmod c^*$,
$d\equiv d_0\bmod v_1^*c_0$,
$d\equiv d^*\bmod c^*$, and further,
noting the congruence property of
$\varkappa_q$  proved in $(15.6)$--$(15.7)$, we may write $(15.11)$ as
$$
\eqalignno{
&S(m,n;c\sqrt{v_1^*v_2^*};u_1/w_1,u_2/w_2)=
\sum_{\scr{ad\,\equiv 1\bmod c}
\atop{\scr{a\bmod v_1^*c}
\atop{\scr{d\bmod v_2^*c}}}}\varkappa_q
\left(\varpi_{\overline{c^*}u_1/w_1}B[a_0,d_0;c_0]
\varpi_{\overline{c^*}u_2/w_2}^{-1}\right)\cr
&\times\exp\left(2\pi i\left({\widetilde{c^*_1}m a_0\over
v_1^*c_0}+ {\widetilde{c^*_2}nd_0\over
v_2^*c_0}\right)\right)\cdot
\exp\left(2\pi i\left({\widetilde{v^*_1c_0}ma^*\over
c^*}+ {\widetilde{v^*_2c_0}nd^*\over c^*}\right)\right).&(15.13) 
}
$$
We have thus obtained the factorization
$$
\eqalignno{
S&(m,n;c\sqrt{v_1^*v_2^*};u_1/w_1,u_2/w_2)\cr
&=S(\widetilde{c_1^*}m,\widetilde{c_2^*}n;c_0\sqrt{v_1^*v_2^*};
\overline{c^*}u_1/w_1,\overline{c^*}u_2/w_2)
S(\widetilde{v_1^*c_0}m,
\widetilde{v_2^*c_0}n;c^*),&(15.14)
}
$$
where the last $S$-factor is an ordinary Kloosterman sum. 
\medskip
In particular, applying $(15.5)$ and
the Weil bound, respectively, to the first and the second factors on
the right side of $(15.14)$, we get
$$
\eqalignno{
|S(m,n;c\sqrt{v_1^*v_2^*};u_1/w_1,u_2/w_2)|&\le v_1^*v_2^*\varphi(c_0)
|S(\widetilde{v_1^*c_0}m,\widetilde{v_2^*c_0}n;c^*)|\cr
&\ll v_1^*v_2^*c_0((m,n,c^*)c^*)^{{1\over2}+\varepsilon},&(15.15)
}
$$
with the implied constant depending only on $\varepsilon$. 
Thus we have, for any $\xi>{1\over2}$,
$$
\eqalignno{
&\sum_c{1\over(c\sqrt{v_1^*v_2^*})^\tau}
|S(m,n;c\sqrt{v_1^*v_2^*};u_1/w_1,u_2/w_2)|\cr
\ll&(v_1^*v_2^*)^{1-{1\over2}\tau}\left(\sum_{c|q^\infty}{1\over
c^{\tau-1}}\right)\left(
\sum_c{(m,n,c)^\xi\over c^{\tau-\xi}}\right),
&(15.16)
}
$$
which is finite if $\tau-\xi>1$.
Therefore, 
we have proved that any $\varGamma_0(q)$ satisfies $(9.6)$ with
$\tau>{3\over2}$.
\medskip
\noindent
{\csc Remark 3}. The methods in [2] and [12] extend to
$M_2(g;A)$ with an arbitrary $A$. Since they are independent of any
non-trivial treatment of generalized Kloosterman sums, 
the above confirmation of $(9.6)$ for generic $\varGamma_0(q)$ could be
regarded as redundant, as far as
the spectral decomposition of $M_2(g;A)$
is concerned. 
\bigskip\noindent
{\bf 16.} With this, we return to the second line of $(6.5)$. 
We stress that hereafter we shall again work with the definition
$(14.1)$. 
\medskip
In view of $(14.8)$ we have
$$
\eqalignno{
&Y_\pm(u,v,w,z;g;d/c;m,n)\cr
&=
\sum_{(k,d)=1}{1\over ck\sqrt{d}}S(n,\pm m;ck\sqrt{d};\infty,1/c)
\tilde{g}_\pm\left(u,v,w,z;4\pi{\sqrt{mn}\over ck\sqrt{d}}\right).&(16.1)
}
$$
Thus Lemma 4 gives the expansion
$$
\eqalignno{
&Y_\pm(u,v,w,z;g;d/c;m,n))\cr
&=\sum_{j=1}^\infty [g]_\pm(\kappa_j;u,v,w,z)
{\overline{\varrho_j(n,\infty)}\varrho_j(\pm m,1/c)
\over\cosh\pi\kappa_j}\cr
&+{1\pm1\over4\pi(4\pi\sqrt{mn})^{2k-1}}\sum_{k=1}^\infty
\Gamma(2k)[g]_+\!\left(\left(\txt{1\over2}-2k\right)i;u,v,w,z\right)
\sum_{j=1}^{\vartheta(k)}\overline{\varrho_{j,k}(n,\infty)}
\varrho_{j,k}(m,1/c)\cr
&+{1\over\pi}\sum_\f{c}\int_{-\infty}^\infty 
[g]_\pm(r;u,v,w,z)(m/n)^{ir}
\overline{e_n\left(\txt{1\over2}+ir;\f{c},\infty\right)}
e_m\left(\txt{1\over2}+ir;\f{c},1/c\right)
dr,&(16.2)
}
$$
where
$$
\eqalignno{
[g]_+(r;u,v,w,z)&={\pi i\over2\sinh\pi r}\int_0^\infty
\left\{J_{2ir}(x)-J_{-2ir}(x)\right\}\tilde{g}_+(u,v,w,z;x){dx\over x},\cr
[g]_-(r;u,v,w,z)&=2\cosh(\pi r)\int_0^\infty K_{2ir}(x)
\tilde{g}_-(u,v,w,z;x){dx\over x}.&(16.3)
}
$$
Further, by $(6.6)$--$(6.7)$ we have that
$$
\eqalignno{
&{(2\pi)^{u-w+1}\over2c^ud^{{1\over2}(u+v-w+z)}}J_+^*(u,v,w,z;g;d/c)\cr
=&\sum_\pm\sum_{j=1}^\infty{[g]_\pm(\kappa_j;u,v,w,z)\over\cosh\pi\kappa_j}
\left\{\sum_m{\overline{\varrho_j(n,\infty)}\sigma_{w+z-1}(n)
\over n^{{1\over2}(u+v+w+z-1)}}\right\}\cr
&\hskip 2cm\times\left\{
\sum_n{\varrho_j(\pm n,1/c)
\over n^{{1\over2}(u-v-w+z+1)}}\right\}\cr
+&{1\over2\pi}\sum_{k=1}^\infty{(2k-1)!\over (4\pi)^{2k-1}}
[g]_+\!\left(\left(\txt{1\over2}-2k\right)i;u,v,w,z\right)
\left\{\sum_m{\overline{\varrho_{j,k}(n,\infty)}\sigma_{w+z-1}(n)
\over n^{k-{1\over2}}n^{{1\over2}(u+v+w+z-1)}}\right\}\cr
&\hskip 2cm\times\left\{
\sum_n{\varrho_{j,k}(\pm n,1/c)
\over n^{k-{1\over2}}n^{{1\over2}(u-v-w+z+1)}}\right\}\cr
+&{1\over\pi}\sum_\pm\sum_\f{c}\int_{-\infty}^\infty
[g]_\pm(r;u,v,w,z)
\left\{\sum_m{\overline{e_n({1\over2}+ir;\f{c},\infty)}\sigma_{w+z-1}(n)
\over n^{{1\over2}(u+v+w+z-1)+ir}}\right\}\cr
&\hskip 2cm\times\left\{
\sum_n{e_n({1\over2}+ir;\f{c},1/c)
\over n^{{1\over2}(u-v-w+z+1)-ir}}\right\}dr,&(16.4)
}
$$
as Lemma 2 and the rapid decay of $[g]_\pm(r;u,v,w,z)$ yield absolute
convergence on the right side, provided $(4.1)$  (see
[11, Section 4.5]).
\bigskip\noindent
{\bf 17.} We need to continue $(16.4)$ to a neighbourhood of the
point $(u,v,w,z)=\left({1\over2},{1\over2},{1\over2},{1\over2}\right)$.
The continuation of $[g]_\pm$ is known already ([11, Section 4.6]), and we
are concerned with the nature of $L$-functions:
$$
\eqalignno{
L_j^\pm(s;1/c)&= \sum_n
\varrho_j(\pm n,1/c)n^{-s},\cr
D_j(s,\alpha)&=\sum_n\overline{\varrho_j(n,\infty)}
\sigma_\alpha(n)n^{-s},\cr
L_{j,k}(s;1/c)&= \sum_n\varrho_{j,k}
(n,1/c)n^{-s-k+{1\over2}},\cr
D_{j,k}(s,\alpha)&=\sum_n
\overline{\varrho_{j,k}(n,\infty)}\sigma_\alpha(n)n^{-s-k+{1\over2}},
&(17.1)
}
$$
where the sums converge absolutely if $\Re s$ is sufficiently
large, because of $(11.2)$. We shall especially require 
uniform bounds for these functions. The Dirichlet
series  involved in the last integral are to be discussed in detail
later, but under the restriction on $A$ mentioned in the introduction.
\medskip
In our continuation procedure of the right side of $(16.4)$, we exploit
the fact that above $L$-functions admit meromorphic continuation to
$\B{C}$ with respect to $s$, and with respect to $\alpha$ as well in
the second and the fourth $L$-functions. To reach $(16.4)$ we appealed
to Lemma 4, and hence the bound $(9.6)$ becomes crucial. Moreover,
the contribution of the continuous spectrum in $(16.4)$ makes it 
clear how important for us to have explicit representation of Fourier
coefficients of Eisenstein series at each cusp, and this is of course
closely related to the structure of generalized Kloosterman sums 
which is partly discussed in Section 15.
\medskip
We begin with relations
between $\sigma_\f{a}$ defined by $(14.1)$ and the two basic involutions
$J:z\mapsto -\overline{z}$, and
$F_q: z\mapsto -1/qz$, which satisfy
$$
J\varGamma_0(q)J^{-1}=\varGamma_0(q),\quad F_q\varGamma_0(q)F_q^{-1}
=\varGamma_0(q).
\eqno(17.2)
$$
We have
$$
J\sigma_\f{a}=\gamma_1\sigma_{\f{b}_1}S^{b_1},\quad F_q\sigma_\f{a}
=\gamma_2\sigma_{\f{b}_2}S^{b_2},\quad \gamma_1,\gamma_2\in\varGamma_0(q),
\;b_1,b_2\in\B{R},
\eqno(17.3)
$$
where $J(\f{a})$, $F_q(\f{a})$ are equivalent to $\f{b}_1$, $\f{b}_2$,
respectively. For instance, the
latter identity is due to the fact that the stabilizer of $\f{b}_2$ is
$$
(\gamma^{-1}_2F_q\sigma_\f{a})\varGamma_\infty(\gamma^{-1}_2
F_q\sigma_\f{a})^{-1}=
\gamma^{-1}_2F_q\varGamma_\f{a}F_q^{-1}\gamma_2\subset\varGamma_0(q)
\eqno(17.4)
$$  
(see the remark made prior to $(7.1)$).
\medskip
The reflection operator $J$ is isometric over $L^2(\varGamma
\backslash{\cal H})$, for
$J(\varGamma\backslash{\cal H})$ is a fundamental domain, 
and
$$
\Vert\psi J\Vert^2=\int_{\varGamma\backslash{\cal H}}|\psi J|^2d\mu
=\int_{J(\varGamma\backslash{\cal H})}|\psi|^2d\mu=
\int_{\varGamma\backslash{\cal
H}}|\psi|^2d\mu=\Vert\psi\Vert^2.\eqno(17.5)
$$
Besides, we have $J\Delta=\Delta J$ as well as the first relation in
$(17.3)$. Hence $\psi_jJ$ is a cusp form belonging to the same
eigenspace as $\psi_j$, for $\psi_jJ(\sigma_\f{a}(z))=
\psi_j(\sigma_{\f{b}_1}(z+b_1))$ converges to $0$ as $z$ tends to
$\infty$.  Thus
$J$ can be diagonalized on each eigenspace of
$\Delta$; that is, we may choose an orthonormal base $\{\psi_j\}$ in such
a way that
$$
\psi_j(-\overline{z})=\epsilon_j\psi_j(z),\quad \epsilon_j=\pm1.
\eqno(17.6)
$$
Also, we observe that
$$
\eqalignno{
J\sigma_{1/c}J\sigma_{1/c}^{-1}&=
\left(\matrix{\sqrt{d}&-f/\sqrt{d}\cr\cr
-c\sqrt{d}&(1+cf)/\sqrt{d}}\right)\left(
\matrix{(1+cf)/\sqrt{d}&-f/\sqrt{d}\cr\cr
-c\sqrt{d}&\sqrt{d}}\right)\cr
&=\left(\matrix{1+2cf&-2f\cr\cr-2c(1+cf)
&1+2cf}\right)\in \varGamma_0(cd).&(17.7)
}
$$
This implies that
$$
\psi_j(J\sigma_{1/c}J(z))=\psi_j(\sigma_{1/c}
(z))\Longleftrightarrow \epsilon_j\psi_j
(\sigma_{1/c}(-\overline{z}))=\psi_j(\sigma_{1/c}(z));\eqno(17.8)
$$
namely
$$
\varrho_j(-n,1/c)=\epsilon_j
\varrho_j(n,1/c).\eqno(17.9)
$$
In particular, we have
$$
L_j^{-}(s;1/c)=\epsilon_jL_j^{+}(s;1/c).
\eqno(17.10)
$$
\medskip
Next, we consider the action of the Fricke operator 
$F_q$. We put $F=F_{cd}$.
Then each $\psi_jF$ is $\varGamma_0(cd)$-invariant, and
is a cusp form such that $\Delta\psi_jF=\nu_j\psi_jF$;
in fact it is a unit vector as 
$$
\int_{\varGamma\backslash{\cal H}}|\psi_jF(z)|^2d\mu(z)
=\int_{F\varGamma\backslash{\cal H}}|\psi_j(z)|^2d\mu(z)=
\int_{\varGamma\backslash{\cal H}}|\psi_j(z)|^2d\mu(z)=1,\eqno(17.11)
$$
for $F\varGamma\backslash{\cal H}$ is a fundamental domain of
$\varGamma=\varGamma_0(cd)$; moreover, $\psi_jF(\sigma_\f{a}(z))
=\psi_j(\sigma_{\f{b}_2}(z+b_2))$ converges to $0$ as $z$ tends to
$\infty$. Since $FJ=JF$, we may assume, besides $(17.6)$, that
$$
\psi_jF=\varpi_j\psi_j,\quad \varpi_j=\pm1.\eqno(17.12)
$$
Further, we observe
$$
\eqalignno{
&\sigma_{1/c}
F\sigma_{1/c}^{-1}F\cr
=&{1\over cd}\left(\matrix{\sqrt{d}&f/\sqrt{d}\cr\cr
c\sqrt{d}&(1+cf)/\sqrt{d}}\right)
\left(\matrix{&-1\cr\cr cd&}\right)
\left(\matrix{(1+cf)/\sqrt{d}&-f/\sqrt{d}\cr\cr
-c\sqrt{d}&\sqrt{d}}\right)\left(\matrix{&-1\cr\cr cd&}\right)\cr
=&{1\over cd}\left(\matrix{cf\sqrt{d}&-\sqrt{d}\cr\cr
c(1+cf)\sqrt{d}&-c\sqrt{d}}\right)
\left(\matrix{-cf\sqrt{d}&-(1+cf)/\sqrt{d}\cr\cr
cd\sqrt{d}&c\sqrt{d}}\right)\cr
=&{1\over cd}\left(\matrix{-(cf)^2d-cd^2&-cf
-(cf)^2-cd\cr\cr c^2f(1+cf)d-(cd)^2
&-c(1+cf)^2-c^2d}\right)\cr
=&\left(\matrix{-cf^2-d&-1
-(1+cf)/d\cr\cr-
cd(1+f(1+cf)/d)&-c-(1+cf)^2/d}\right)
\in \varGamma_0(cd).&(17.13)
}
$$
Hence we have
$$
\eqalignno{
\psi_j(\sigma_{1/c}F(z))
&=\psi_j(F^{-1}\sigma_{1/c}(z))\cr
&=\psi_j(F\sigma_{1/c}(z))=\varpi_j\psi_j
(\sigma_{1/c}(z));&(17.14)
}
$$
that is, we have
$$
\psi_j\left(\sigma_{1/c}\left(-1/cdz\right)\right)=
\varpi_j\psi_j\left(\sigma_{1/c}\left(z\right)\right).
\eqno(17.15)
$$
\bigskip\noindent
{\bf 18.} We may now prove the functional equation for
$L_j(s;1/c)=L^+_j(s;1/c)$; note that
we have $(17.10)$. We have to discuss two cases separately according as
$\epsilon_j=+1$ or $-1$.
\medskip
The case $\epsilon_j=+1$: We have, by $(17.9)$,
$$
\eqalignno{
&\int_0^\infty\psi_j\left(\sigma_{1/c}
\left({iy\over\sqrt{cd}}\right)\right)y^{s-{3\over2}}dy\cr
&=2
\int_0^\infty \left({y\over\sqrt{cd}}\right)^{1\over2}
\sum_{n>0}\varrho_j(n,1/c)
K_{i\kappa_j}\left(2\pi{n\over\sqrt{cd}}y\right)y^{s-{3\over2}}dy\cr
&=2^{s-1}(cd)^{-{1\over4}}\left({2\pi\over\sqrt{cd}}\right)^{-s}
\Gamma\left(\txt{1\over2}(s+i\kappa_j)\right)
\Gamma\left(\txt{1\over2}(s-i\kappa_j)\right)
L_j(s;1/c).&(18.1)
}
$$
On the other hand, by $(17.15)$,
$$
\eqalignno{
&\int_0^\infty\psi_j\left(\sigma_{1/c}
\left({iy\over\sqrt{cd}}\right)\right)y^{s-{3\over2}}dy\cr
&=\int_1^\infty\psi_j\left(\sigma_{1/c}
\left({iy\over\sqrt{cd}}\right)\right)y^{s-{3\over2}}dy+
\int_1^\infty\psi_j\left(\sigma_{1/c}
\left({i\over\sqrt{cd}y}\right)\right)y^{1-s-{3\over2}}dy\cr
&=\int_1^\infty\left\{\psi_j\left(\sigma_{1/c}
\left({iy\over\sqrt{cd}}\right)\right)y^{s-{3\over2}}+\varpi_j
\psi_j\left(\sigma_{1/c}
\left({iy\over\sqrt{cd}}
\right)\right)y^{1-s-{3\over2}}\right\}dy,&(18.2)
}
$$
which is entire in $s$, for $\psi_j\sigma_{1/c}$ decays
exponentially as $y$ tends to $+\infty$. Namely, the function
$L_j(s;1/c)$ is entire, and we have
$$
\eqalignno{
&\left({\pi\over\sqrt{cd}}\right)^{-s}
\Gamma\left(\txt{1\over2}(s+i\kappa_j)\right)
\Gamma\left(\txt{1\over2}(s-i\kappa_j)\right)
L_j(s;1/c)\cr
&=\varpi_j\left({\pi\over\sqrt{cd}}\right)^{1-s}
\Gamma\left(\txt{1\over2}(1-s+i\kappa_j)\right)
\Gamma\left(\txt{1\over2}(1-s-i\kappa_j)\right)
L_j(1-s;1/c),&(18.3)
}
$$
By the duplication formula for $\Gamma$-function, one may transform
this relation into
$$
\eqalignno{
L_j(s;1/c)
&={\varpi_j\over\pi}\left({2\pi\over\sqrt{cd}}\right)^{2s-1}
\Gamma(1-s+i\kappa_j)\Gamma(1-s-i\kappa_j)\cr
&\times\left(\cosh\pi\kappa_j-\cos\pi s\right)
L_j(1-s;1/c).&(18.4)
}
$$
\medskip
The case $\epsilon_j=-1$: We have
$$
\psi_j(\sigma_{1/c}(z))=2i\sqrt{y}\sum_{n>0}
\varrho_j(n,1/c)K_{i\kappa_j}(2\pi ny)
\sin(2\pi nx),\eqno(18.5)
$$
We put $f_j(z)=\partial_x\psi_j(\sigma_{1/c}(z-\overline{c}/d))$. We have
$$
f_j(z)=4\pi i\sqrt{y}\sum_{n>0}n
\varrho_j(n,1/c)K_{i\kappa_j}(2\pi ny)
\cos(2\pi nx),\eqno(18.6)
$$
which implies that as $x\to0$
$$
\psi_j(\sigma_{1/c}(z))=f_j(iy)x+O(x^2)\eqno(18.7)
$$
as well as
$$
\eqalignno{
\psi_j\left(\sigma_{1/c}(-1/cdz)\right)
&=\psi_j\left(\sigma_{1/c}(i/cdy-x/cdy^2+O(x^2))\right)\cr
&=-(x/cdy^2)f_j(i/cdy)+O(x^2);&(18.8)
}
$$
that is,
$$
f_j(i/cdy)=-\varpi_jcdy^2f_j(iy).\eqno(18.9)
$$
Hence, 
$$
\eqalignno{
&\int_0^\infty f_j\left(
{iy\over\sqrt{cd}}\right)y^{s-{1\over2}}dy\cr
&=\pi i(cd)^{-{1\over4}}\left({\pi\over\sqrt{cd}}\right)^{-s-1}
\Gamma\left(\txt{1\over2}(1+s+i\kappa_j)\right)
\Gamma\left(\txt{1\over2}(1+s-i\kappa_j)\right)
L_j(s;1/c);&(18.10)
}
$$
and
$$
\int_0^\infty f_j\left(
{iy\over\sqrt{cd}}\right)y^{s-{1\over2}}dy=
\int_1^\infty\left\{f_j\left(
{iy\over\sqrt{cd}}\right)y^{s-{1\over2}}-\varpi_jf_j\left(
{iy\over\sqrt{cd}}\right)y^{1-s-{1\over2}}\right\}dy.\eqno(18.11)
$$
Namely, we have that
$$
\eqalignno{
L_j(s;1/c)
&=-{\varpi_j\over\pi}\left({2\pi\over\sqrt{cd}}\right)^{2s-1}
\Gamma(1-s+i\kappa_j)\Gamma(1-s-i\kappa_j)\cr
&\times\left(\cosh\pi\kappa_j+\cos\pi s\right)
L_j(1-s;1/c).&(18.12)
}
$$
\medskip
\noindent
{\bf Lemma 6.} {\it The function 
$L_j(s;1/c)$
is entire, and it holds that for any $s$
$$
\eqalignno{
L_j(s;1/c)
&={\varpi_j\over\pi}\left({2\pi\over\sqrt{cd}}\right)^{2s-1}
\Gamma(1-s+i\kappa_j)\Gamma(1-s-i\kappa_j)\cr
&\times\left(\epsilon_j\cosh\pi\kappa_j-\cos\pi s\right)
L_j(1-s;1/c).&(18.13)
}
$$
We have also
$$
L_j(s;1/c)\ll
(\kappa_j+|s|+1))^{c_0}\exp\left(\txt{1\over2}\pi\kappa_j\right),
\eqno(18.14)
$$
where the constant $c_0$ depends at most on $\Re s$, and the implied
constant on $\Re s$.
\/}
\medskip
\noindent
The second assertion follows via a convexity argument.
\medskip
We may omit the discussion on $L_{j,k}$, as it is analogous to $L_j$.
\bigskip\noindent
{\bf 19.} We turn to $D_j(s,\alpha)$. There are
at least two possible ways for us to take here. One is to exploit the
theory of Hecke operators in order to
relate $D_j$ with a product of two values of Hecke
$L$-functions analogously as we did in the case of $M_2(g;1)$ in [11].
However, the cusp form $\psi_j$ cannot generally be assumed to be
such that the corresponding Hecke series is fully
decomposed into an Euler product. This is because those
$\varrho_j(n,\infty)$ with $n|(cd)^\infty$ are not well related to
eigenvalues of Hecke operators, and thus the corresponding part of
$D_j(s,\alpha)$ causes difficulties in the continuation as well as the
estimation procedures, which is a serious drawback of the method
as far as our present purpose is concerned. One may
appeal to the notion of new forms 
whose Hecke series admits a full Euler product; yet it does not seem
to resolve our difficulties. Hence, we shall take the second method
which is in fact  a special instance of applications of
Rankin's unfolding method (see [11, pp.\ 181--182]). This causes, however,
still a technical difficulty, for it requires us to have an explicit
description of the scattering matrix of $\varGamma_0(q)$ and all
Fourier coefficients of Eisenstein series at each cusp
 (see $(24.1)$ below). This task is highly involved. The note [1]
contains, in fact, a discussion of the arithmetical nature of those
Fourier coefficients and the result appears to be essentially adequate
for our purpose, if we let our reasoning in the later sections be 
somewhat inexplicit; note that the same can be done by 
extending $(15.14)$ to a full localization. 
Under such a circumstance, it may be appropriate for us to make here a
compromise by introducing the assumption that $A$ is defined by a sum over
square-free integers, as underlined in the introduction.
Since we have $(14.7)$, this eases our task considerably, yet it
does not seem to restrict the scope of our method. In the future, we shall
work out a fuller account of $M_2(d;A)$.
\bigskip\noindent
{\bf 20.} Thus, we shall hereafter assume that
$$
\hbox{$q=cd$ is square-free.}\eqno(20.1)
$$
\medskip
By Lemma 5 in Section 12, we have now
$$
\left\{\hbox{inequivalent cusps of $\varGamma_0(q)$}\right\}\equiv
\left\{{1\over w}: w|q\right\};\eqno(20.2)
$$
and we have $(14.7)$ for any combination of cusps. In particular, 
for those Hecke congruence groups that are relevant in the sequel,
$(9.6)$ and thus Lemmas 2--4 have been verified, without the discussion in
Section 15.
\medskip
To make Lemmas 3--4 more explicit, let us 
compute the Fourier coefficients of Eisenstein series at each cusp. 
Thus, by the assertion $(14.7)$, 
$$
\eqalignno{
&E(\sigma_{1/w_2}(z),1/w_1;s)\cr
&=\delta_{w_1,w_2}y^s+\sqrt{\pi}y^{1-s}{\Gamma\left(s-{1\over2}\right)\over
\Gamma(s)}\sum_{((v_1,w_2)(w_1,v_2),r)=1}{\varphi((v_1,v_2)(w_1,w_2)r)
\over((w_1,w_2)r\sqrt{v_1v_2})^{2s}}\cr
&+2\sqrt{y}{\pi^s\over\Gamma(s)}\sum_{n\ne0}\exp(2\pi inx)
K_{s-{1\over2}}(2\pi|n|y)|n|^{s-{1\over2}}
\sum_{((v_1,w_2)(w_1,v_2),r)=1}{c_{(v_1,v_2)(w_1,w_2)r}(n)
\over((w_1,w_2)r\sqrt{v_1v_2})^{2s}},\qquad&(20.3)
}
$$
where the last numerator is a Ramanujan sum. We have
$$
\eqalignno{
&\sum_{((v_1,w_2)(w_1,v_2),r)=1}{\varphi((v_1,v_2)(w_1,w_2)r)
\over((w_1,w_2)r\sqrt{v_1v_2})^{2s}}\cr
&={1\over(w_1,w_2)^{2s}(v_1,v_2)^s}
\left\{\sum_{r|((v_1,v_2)(w_1,w_2))^\infty}
{\varphi((v_1,v_2)(w_1,w_2)r)
\over r^{2s}}\right\}
\left\{\sum_{(r,q)=1}{\varphi(r)\over r^{2s}}\right\}\cr
&={1\over(w_1,w_2)^{2s}(v_1,v_2)^s}\prod_{p|(v_1,v_2)}
\left(\sum_{j=0}^\infty{\varphi(p^{j+1})\over p^{2js}}\right)
\prod_{p|(w_1,w_2)}
\left(\sum_{j=0}^\infty{\varphi(p^{j+1})\over p^{2js}}\right)
\left\{\sum_{(r,q)=1}{\varphi(r)\over r^{2s}}\right\}\cr
&={\zeta(2s-1)\over\zeta(2s)}\prod_{p|(v_1,v_2)(w_1,w_2)}
\left({p-1\over p^{2s}-1}\right)
\prod_{p|(v_1,w_2)(w_1,v_2)}
\left({p^s-p^{1-s}\over p^{2s}-1}\right).&(20.4)
}
$$
\medskip
Next,
$$
\eqalignno{
&\sum_{((v_1,w_2)(w_1,v_2),r)=1}{c_{(v_1,v_2)(w_1,w_2)r}(n)
\over((w_1,w_2)r\sqrt{v_1v_2})^{2s}}\cr
&={1\over(w_1,w_2)^{2s}(v_1v_2)^s}
\left\{\sum_{r|((v_1,v_2)(w_1,w_2))^\infty}
{c_{(v_1,v_2)(w_1,w_2)r}(n)
\over r^{2s}}\right\}
\left\{\sum_{(r,q)=1}{c_r(n)\over r^{2s}}\right\}.&(20.5)
}
$$
We have
$$
\eqalignno{
&\sum_{r|((v_1,v_2)(w_1,w_2))^\infty}
{c_{(v_1,v_2)(w_1,w_2)r}(n)
\over r^{2s}}=\prod_{p|(v_1,v_2)(w_1,w_2)}\left\{
\sum_{j=0}^\infty{c_{p^{j+1}}(n)\over p^{2js}}\right\}\cr
&=\prod_{p|(v_1,v_2)(w_1,w_2)}p^{2s}\left\{
\sum_{j=0}^\infty{c_{p^{j}}(n)\over p^{2js}}-1\right\}\cr
&=((v_1,v_2)(w_1,w_2))^{2s}
\prod_{p|(v_1,v_2)(w_1,w_2)}\left\{\sigma_{1-2s}(n_p)\left(1-{1\over
p^{2s}}\right)-1
\right\}&(20.6)
}
$$
and
$$
\sum_{(r,q)=1}{c_r(n)\over r^{2s}}={\sigma_{1-2s}(n,\chi_q)
\over L(2s,\chi_q)},\eqno(20.7)
$$
where $n_p=(n,p^\infty)$ and $\chi_q$ is the principal character
$\bmod q$. Thus,
$$
\eqalignno{
&\sum_{((v_1,w_2)(w_1,v_2),r)=1}{c_{(v_1,v_2)(w_1,w_2)r}(n)
\over((w_1,w_2)r\sqrt{v_1v_2})^{2s}}\cr
&={\sigma_{1-2s}(n,\chi_q)
\over L(2s,\chi_q)}\left({(v_1,v_2)\over[v_1,v_2]}\right)^s
\prod_{p|(v_1,v_2)(w_1,w_2)}\left\{\sigma_{1-2s}(n_p)\left(1-{1\over
p^{2s}}\right)-1
\right\}.&(20.8)
}
$$
\medskip
Collecting these assertions, we obtain in particular that
\medskip
\noindent
{\bf Lemma 7.} {\it The function $s(1-s)\Gamma(s)L(2s,\chi_{cd})
E(\sigma_{1/w_2}(z),1/w_1;s)$ is regular 
for all $s$, and it is
$\ll y^{\Re s}+y^{1-\Re s}$
as $y=\Re z$ tends to infinity, as far as $s$ remains bounded.
\/}
\bigskip\noindent
{\bf 21.} Lemma 2 holds safely for $\varGamma=\varGamma_0(q)$, 
$\mu(q)\ne0$, and Lemmas 3 and 4 become as follows (see [10]):
\medskip
\noindent
{\bf Lemma 8.} {\it
Let $h(r)$ be even, regular and of fast decay on the strip $|\r{Im}
r|<{1\over2}+\eta$ with an $\eta>0$.
Then it holds that for any $m,n>0$ and $w_1|q,w_2|q$
$$
\eqalignno{
&\sum_{j=1}^\infty{\overline{\varrho_j(m,1/w_1)}
\varrho_j(\pm n,1/w_2)
\over\cosh\pi\kappa_j}h(\kappa_j)\cr
+&{1\over\pi}\sum_{q=vw}\int_{-\infty}^\infty \left({n\over m}\right)^{ir}
{\sigma_{2ir}(m;\chi_q)\sigma_{-2ir}(n;\chi_q)\over
|L(1+2ir,\chi_q)|^2}
\left({(v,v_1)\over[v,v_1]}\right)^{{1\over2}-ir}
\left({(v,v_2)\over[v,v_2]}\right)^{{1\over2}+ir}\cr
&\qquad\times\prod_{p|(v,v_1)(w,w_1)}\left\{\sigma_{2ir}(m_p)
\left(1-{1\over p^{1-2ir}}\right)-1\right\}\cr
&\qquad\times\prod_{p|(v,v_2)(w,w_2)}
\left\{\sigma_{-2ir}(n_p)\left(1-{1\over
p^{1+2ir}}\right)-1\right\}h(r)dr\cr
=&{1\over\pi^2}\delta_{w_1,w_2}\delta_{m,\pm n}\exp(2\pi imb_{w_1,w_2})
\int_{-\infty}^\infty r\tanh(\pi r)h(r)dr\cr
+&\sum_{(r,\, (v_1,w_2)(w_1,v_2))=1}{1\over (w_1,w_2)r\sqrt{v_1v_2}}\cr
&\times S(\overline{(v_1,w_2)}m,\pm\overline{(w_1,v_2)}n;
(v_1,v_2)(w_1,w_2)r)h_\pm\!\left({4\pi\sqrt{mn}\over
(w_1,w_2)r\sqrt{v_1v_2}}\right),&(21.1)  
}
$$
with $h_\pm$ as in $(11.4)$.
\/}
\medskip
\noindent
{\bf Lemma 9.} {\it
Let $\varphi$ be smooth and of fast decay over the positive real axis.
Then we have, for any $m,n>0$ and $w_1|q,w_2|q$,
$$
\eqalignno{
&\sum_{(r,\, (v_1,w_2)(w_1,v_2))=1}
{ S(\overline{(v_1,w_2)}m,\pm\overline{(w_1,v_2)}n;
(v_1,v_2)(w_1,w_2)r)\over (w_1,w_2)r\sqrt{v_1v_2}}
\varphi\!\left({4\pi\sqrt{mn}\over
(w_1,w_2)r\sqrt{v_1v_2}}\right)\cr
=&\sum_{j=1}^\infty
{\overline{\varrho_j(m,1/w_1)}\varrho_j(\pm n,1/w_2)
\over\cosh\pi\kappa_j}\hat{\varphi}_\pm(\kappa_j)\cr
+&{1\pm1\over4\pi(4\pi\sqrt{mn})^{2k-1}}\sum_{k=1}^\infty
\Gamma(2k)\hat\varphi_+\!\left(\left(\txt{1\over2}-2k\right)i\right)
\sum_{j=1}^{\vartheta(k)}\overline{\varrho_{j,k}(m,1/w_1)}
\varrho_{j,k}(n,1/w_2)\cr
+&{1\over\pi}\sum_{q=vw}\int_{-\infty}^\infty \left({n\over m}\right)^{ir}
{\sigma_{2ir}(m;\chi_q)\sigma_{-2ir}(n;\chi_q)\over
|L(1+2ir,\chi_q)|^2}
\left({(v,v_1)\over[v,v_1]}\right)^{{1\over2}-ir}
\left({(v,v_2)\over[v,v_2]}\right)^{{1\over2}+ir}\cr
&\qquad\times\prod_{p|(v,v_1)(w,w_1)}\left\{\sigma_{2ir}(m_p)
\left(1-{1\over p^{1-2ir}}\right)-1\right\}\cr
&\qquad\times\prod_{p|(v,v_2)(w,w_2)}
\left\{\sigma_{-2ir}(n_p)\left(1-{1\over
p^{1+2ir}}\right)-1\right\}\hat{\varphi}_\pm(r)dr,&(21.2)
}
$$
where $\hat{\varphi}_\pm$ are as in $(11.6)$.
\/}
\medskip
We specialize the last assertion as in $(14.8)$, and have,
in place of $(16.4)$,
$$
\eqalignno{
&{(2\pi)^{u-w+1}\over2c^ud^{{1\over2}(u+v-w+z)}}
J_+^*(u,v,w,z;g;d/c)\cr
=&\sum_\pm\sum_{j=1}^\infty
{[g]_\pm(\kappa_j;u,v,w,z)\over\cosh\pi\kappa_j}\cr
&\hskip1cm\times\left\{\sum_n{\overline{\varrho_j(n,\infty)}
\sigma_{w+z-1}(n)
\over n^{{1\over2}(u+v+w+z-1)}}\right\}\left\{
\sum_n{\varrho_j(\pm n,1/c)
\over n^{{1\over2}(u-v-w+z+1)}}\right\}\cr
+&2\sum_{k=1}^\infty{(2k-1)!\over (4\pi)^{2k}}
[g]_+\!\left(\left(\txt{1\over2}-2k\right)i;u,v,w,z\right)\cr
&\hskip1cm\times\left\{\sum_n{\overline{\varrho_{j,k}(n,\infty)}
\sigma_{w+z-1}(n)
\over n^{k-{1\over2}}n^{{1\over2}(u+v+w+z-1)}}\right\}
\left\{
\sum_n{\varrho_{j,k}(n,1/c)
\over n^{k-{1\over2}}n^{{1\over2}(u-v-w+z+1)}}\right\}\cr
+&{1\over\pi d^{{1\over2}+ir}}
\sum_\pm\sum_{cd=c_1d_1}{1\over d_1}\int_{-\infty}^\infty 
{(d_1,d)^{1+2ir}\over |L(1+2ir,\chi_{cd})|^2}
[g]_\pm(r;u,v,w,z)\cr
&\hskip1cm\times\sum_n{\sigma_{2ir}(n;\chi_{cd})\sigma_{w+z-1}(n)
\over n^{{1\over2}(u+v+w+z-1)+ir}}\prod_{p|c_1}\left\{\sigma_{2ir}(n_p)
\left(1-{1\over p^{1-2ir}}\right)-1\right\}\cr
&\hskip1cm\times\sum_n{\sigma_{-2ir}(n;\chi_{cd})\over
n^{{1\over2}(u-v-w+z+1)-ir}}\prod_{p|(c_1,c)(d_1,d)}
\left\{\sigma_{-2ir}(n_p)\left(1-{1\over
p^{1+2ir}}\right)-1\right\}dr,&(21.3)
}
$$
with $q=cd$, where we have used the fact that
$\sigma_{1/cd}\in\varGamma_0(cd)$ and thus $\varrho_j(n,1/cd)=
\varrho_j(n,\infty)$, $\varrho_{j,k}(n,1/cd)=
\varrho_{j,k}(n,\infty)$.
\bigskip\noindent
{\bf 22.} We now deal with the function $D_j(s,\alpha)$. As remarked
in Section 19, we shall employ the unfolding method.
\medskip
To this end we introduce the scattering matrix $\B{S}$ of
$\varGamma_0(cd)$. We thus write
$(20.3)$ as
$$
E(\sigma_{1/w_2}(z),1/w_1;s)=\delta_{w_1,w_2}y^s
+\varphi(s;w_1,w_2)y^{1-s}+\cdots.\eqno(22.1)
$$
We put
$$
\B{S}(s)=\Bigg(\varphi(s;w_1,w_2)\Bigg)_{w_1,w_2|cd}.\eqno(22.2)
$$
and
$$
\B{E}(s)=\left(\matrix{\vdots\cr E(z,1/w;s)\cr
\vdots}\right)_{w|cd},\eqno(22.3)
$$
so that
$$
\B{E}(s)=\left(\matrix{1&&\cr&\ddots&\cr&&1}\right)y^s+\B{S}(s)y^{1-s}
+\cdots,\eqno(22.4)
$$
where the error terms decays exponentially and is
$O(y^{{1\over2}-\varepsilon})$ as $y$ tends to infinity and to $0$,
respectively. 
\medskip
We have the functional equation
$$
\B{E}(z,s)=\B{S}(s)\B{E}(z,1-s),\eqno(22.5)
$$
provided both sides are finite.
To confirm this, we let $\Re s$, $\Im s$ be sufficiently large.
Then $(22.4)$ implies in particular that 
$\B{E}(z,1-s)-\B{S}(1-s)\B{E}(z,s)$
is in an obvious vector extension of $L^2(\varGamma_0(q)\backslash{\cal
H})$. However, this vector function, if not trivial, has the eigenvalue
$s(1-s)$ against  $\Delta$ the hyperbolic Laplacian. Since $\Delta$ is
self-adjoint, its eigenvalues $s(1-s)$ should be real, which is 
a contradiction, and hence $(22.5)$ holds for 
all complex $s$ by analytic continuation as far as $\B{E}(z,s)$ is finite.
Consequently, we have got also
$$
\B{S}(s)\B{S}(1-s)=1.\eqno(22.6)
$$
\bigskip
\noindent
{\bf 23.}
We shall assume $\epsilon_j=1$ till the end of Section 24. 
\medskip
Let $E(z,s)$ be the Eisenstein series for
$\varGamma_0(1)$, and put $E^*(z,s)=\pi^{-s}\Gamma(s)\zeta(2s)
E(z,s)$, so that
$$
E^*(z,s)=E^*(z,1-s)\eqno(23.1)
$$
and
$$
\eqalignno{
E^*(z,s)&=\pi^{-s}\Gamma(s)\zeta(2s)y^s+\pi^{s-1}
\Gamma\left(1-s\right)\zeta(2(1-s))y^{1-s}\cr
&+2\sqrt{y}\sum_{n\ne0}|n|^{s-{1\over2}}\sigma_{1-2s}(n)
K_{s-{1\over2}}(2\pi|n|y)\exp(2\pi nx),&(23.2)
}
$$
which shows that $s(1-s)E^*(z,s)$ is regular for all $s$.
We have, on a suitable assumption on $s,\alpha$ to secure 
convergence, that
$$
\eqalignno{
&\int_{\varGamma_0(cd)\backslash{\cal H}}\overline{\psi_j(z)}E^*
\!\left(z,\txt{1\over2}(1-\alpha)\right)
E\left(z,\infty;s-\txt{1\over2}\alpha\right)d\mu(z)\cr
=&\int_{\varGamma_\infty\backslash{\cal H}}\overline{\psi_j(z)}E^*
\!\left(z,\txt{1\over2}(1-\alpha)\right)y^{s-{1\over2}\alpha}d\mu(z)\cr
=&4\sum_{n>0}n^{-{1\over2}\alpha}\sigma_\alpha(n)
\overline{\varrho_j(n,\infty)}\int_0^\infty K_{{1\over2}\alpha}
(2\pi ny)K_{i\kappa_j}(2\pi ny)y^{s-{1\over2}\alpha-1}dy\cr
=&{\Gamma(s,\alpha;\kappa_j)\over
2\pi^{s-{1\over2}\alpha}
\Gamma\left(s-{1\over2}\alpha\right)}D_j(s,\alpha),
&(23.3)
}
$$
with
$$
\Gamma(s,\alpha;\kappa)=\Gamma\left(\txt{1\over2}(s+i\kappa)\right)
\Gamma\left(\txt{1\over2}(s-i\kappa)\right)
\Gamma\left(\txt{1\over2}(s-\alpha+i\kappa)\right)
\Gamma\left(\txt{1\over2}(s-\alpha-i\kappa)\right)
.\eqno(23.4)
$$
\medskip
On noting this, we consider also the relation
$$
\eqalignno{
\int_{\varGamma_0(cd)\backslash{\cal H}}&\overline{\psi_j(z)}E^*
\!\left(z,\txt{1\over2}(1-\alpha)\right)
E\left(z,1/w;s-\txt{1\over2}\alpha\right)d\mu(z)\cr
&={\Gamma(s,\alpha;\kappa_j)\over
2\pi^{s-{1\over2}\alpha}
\Gamma\left(s-{1\over2}\alpha\right)}
D_j(s,\alpha;1/w),
&(23.5)
}
$$
where
$$
D_j(s,\alpha;1/w)=\sum_{n>0}\overline{\varrho_j(n,1/w)}
\sigma_\alpha(n,1/w)n^{-s},\eqno(23.6)
$$
with $\sigma_\alpha(n,1/w)$ an analogue of $\sigma_\alpha(n)$. 
\medskip
By
$(13.7)$, we have $E^*(\sigma_{1/w}(z),s)=E^*(\tau_v(z),s)=E^*(vz,s)$,
and thus
$$
\eqalignno{
E^*(\sigma_{1/w}(z),s)&=\pi^{-s}\Gamma(s)\zeta(2s)(vy)^s
+\pi^{s-1}
\Gamma\left(1-s\right)\zeta(2(1-s))(vy)^{1-s}\cr
&+2\sqrt{y}\sum_{n\ne0}|n|^{s-{1\over2}}\sigma_{1-2s}(n,1/w)
K_{s-{1\over2}}(2\pi|n|y)\exp(2\pi nx),&(23.7)
}
$$
That is, we have
$$
\sigma_\alpha(n,1/w)=
v^{{1\over2}(\alpha+1)}\sigma_\alpha(n/v),\eqno(23.8)
$$
which vanishes if $v\nmid n$.
\medskip
Put
$$
\B{D}_j(s,\alpha)=
\left(\matrix{\vdots\cr
D_j(s,\alpha;1/w)\cr
\vdots}\right)_{w|cd}.\eqno(23.9)
$$
Then we have, by $(22.5)$ and $(23.1)$,
$$
{\Gamma(s,\alpha;\kappa_j)\over
\pi^{s-{1\over2}\alpha}
\Gamma\left(s-{1\over2}\alpha\right)}
\B{D}_j(s,\alpha)=
{\Gamma(1-s,-\alpha;\kappa_j)\over
\pi^{1-s+{1\over2}\alpha}
\Gamma\left(1-s+{1\over2}\alpha\right)}
\B{S}\left(s-\txt{1\over2}\alpha\right)
\B{D}_j(1-s,-\alpha).\eqno(23.10)
$$
In particular, we get the functional equation
$$
\eqalignno{
D_j(s,\alpha)&=\pi^{2s-\alpha-1}
{\Gamma(1-s,-\alpha;\kappa_j)\over
\Gamma(s,\alpha;\kappa_j)}{
\Gamma\left(s-{1\over2}\alpha\right)\over
\Gamma\left(1-s+{1\over2}\alpha\right)}\cr
&\times\sum_{w|cd}
\varphi\left(s-\txt{1\over2}\alpha;
\infty,1/w\right)D_j(1-s,-\alpha;1/w).&(23.11)
}
$$
\bigskip
\noindent
{\bf 24.} We decompose the left side of $(23.5)$ as
$$
\eqalignno{
\sum_{w_1|cd}&\int_0^1\!\int_{y_0}^\infty
\overline{\psi_j(\sigma_{1/w_1}(z))}E^*
\!\left(\sigma_{1/w_1}(z),\txt{1\over2}(1-\alpha)\right)
E\left(\sigma_{1/w_1}(z),1/w;s-\txt{1\over2}\alpha\right)d\mu(z)\cr
+&\int_{(\varGamma_0(cd)\backslash{\cal H})_{y_0}
}\overline{\psi_j(z)}E^*
\!\left(z,\txt{1\over2}(1-\alpha)\right)
E\left(z,1/w;s-\txt{1\over2}\alpha\right)d\mu(z),&(24.1)
}
$$
where $y_0$ is chosen so that the remainder domain
$(\varGamma_0(cd)\backslash{\cal H})_{y_0}$ is a compact set in
${\cal H}$. We then apply Lemma 7 and $(23.7)$ to each term of
$(24.1)$. We obtain the crucial assertion
\medskip
\noindent
{\bf Lemma 10.} {\it The functions
$$
\left(1-\alpha^2\right)\left(s-\txt{1\over2}\alpha\right)
\left(1-s+\txt{1\over2}\alpha\right)\Gamma(2s-\alpha)
L(2s-\alpha,\chi_{cd})D_j(s,\alpha;1/w)\eqno(24.2)
$$
of the complex variables $s$ and $\alpha$ 
are all entire over $\B{C}^2$.
}
\medskip
\noindent
In fact, it suffices to note
that the multiple of $(24.1)$ by the factor 
$\left(1-\alpha^2\right)\left(s-\txt{1\over2}\alpha\right)
\left(1-s+\txt{1\over2}\alpha\right)\break\Gamma(2s-\alpha)
L(2s-\alpha,\chi_{cd})$ is regular in $s$ and $\alpha$ by Lemma 7.
\medskip
On the other hand, we have, by $(20.4)$,
$$
\eqalignno{
&L(2s-\alpha,\chi_{cd})
\varphi\left(s-\txt{1\over2}\alpha;\infty,1/w\right)\cr
&={1\over\pi}\varphi(w)\left(\pi\over cd\right)^{2s-\alpha}
\zeta(2(1-s)+\alpha)
{\Gamma\left(1-s+{1\over2}\alpha\right)\over
\Gamma\left(s-{1\over2}\alpha\right)}
\prod_{p|v}
\left({p^{s-{1\over2}\alpha}-p^{1-s+{1\over2}\alpha}}\right).&(24.3)
}
$$
Inserting this into $(23.11)$, we get
$$
\eqalignno{
L(2s-\alpha,\chi_{cd})D_j(s,\alpha)
&={1\over\pi^2}\left({\pi^2\over
cd}\right)^{2s-\alpha}\zeta(2(1-s)+\alpha)
{\Gamma(1-s,-\alpha;\kappa_j)\over\Gamma(s,\alpha;\kappa_j)}\cr
&\times\sum_{w|cd}\varphi(w)
\prod_{p|v}
\left({p^{s-{1\over2}\alpha}-p^{1-s+{1\over2}\alpha}}\right)
D_j(1-s,-\alpha;1/w).&(24.4)
}
$$
We then let $\Re s$ be negative and so large that
both $\zeta(2(1-s)+\alpha)$ and $D_j(1-s,-\alpha;1/w)$
are absolutely convergent. In this way we obtain, via 
Lemma 2, Stirling's formula, and the convexity argument,
\medskip
\noindent
{\bf Lemma 11.} {\it Provided that $\Re s$ and $\alpha$ are bounded,
we have
$$
\eqalignno{
\left(1-\alpha^2\right)\left(s-\txt{1\over2}\alpha\right)&
\left(1-s+\txt{1\over2}\alpha\right)
L(2s-\alpha,\chi_{cd})D_j(s,\alpha;1/w)\cr
&\ll
(\kappa_j+|s|+1)^\tau\exp\left(\txt{1\over2}\pi\kappa_j\right),
&(24.5)
}
$$
where $\tau$ depends only on $\Re s$ and $\Re\alpha$, and the
implied constant additionally on $cd$ too.
}
\bigskip
\noindent
{\bf 25.} We still need to deal with the case $\epsilon_j=-1$. Here we
shall have to overcome an additional technical difficulties, because
Eisenstein series of non-zero weights naturally come up in our
argument (see [11, Section 3.2]).
\medskip
We introduce 
$$
\psi_j^-(z)=y(\partial_x-i\partial_y)\psi_j(z),\eqno(25.1)
$$ 
with our present vector $\psi_j$ such that $\psi_jJ=-\psi_j$. We have
$$
\psi_j^-(\gamma(z))=\psi_j^-(z)(\jmath(\gamma,z)/|\jmath(\gamma,z)|)^2,
\quad\gamma\in\varGamma.
\eqno(25.2)
$$
In fact, writing
$\xi=\Re\lambda(z)$, $\eta=\Im\lambda(z)$ for a regular function
$\lambda$, we have
$(\partial_x-i\partial_y)[H(\lambda(z))]=
\{(\partial H/\partial\xi)(\partial\xi/\partial x-
i\partial\xi/\partial y)+
(\partial H/\partial\eta)(\partial\eta/\partial x-
i\partial\eta/\partial y)\}=\{(\partial H/\partial\xi
-i\partial H/
\partial\eta)(\partial\xi/\partial x-i\partial\xi/\partial y)\}
=[(\partial_\xi-i\partial_\eta)H]
(d\lambda/dz)$ by the Cauchy--Riemann equation applied to $\lambda$.
We put $\lambda=\gamma$, $H=\psi_j$, and get
$y(\partial_x-i\partial_y)\psi_j(z)
=(y/\eta)(d\gamma/dz)\eta(\partial_\xi-i\partial_\eta)\psi_j(\xi+i\eta)$,
which confirms $(25.2)$.
\medskip
To offset the automorphic factor in $(25.2)$, we introduce
$$
E_-(z,1/w;s)=\sum_{\gamma\in\varGamma_{1/w}\backslash\varGamma}
\left(\Im \sigma_{1/w}^{-1}\gamma(z)\right)^s
\left(\jmath(\sigma_{1/w}^{-1}\gamma ,z)/
|\jmath(\sigma_{1/w}^{-1}\gamma
,z)|\right)^{-2}.\eqno(25.3)
$$
We should note the relation
$$
y(\partial_x-i\partial_y)
[E(z,1/w;s)]=-is E_-(z,1/w;s),\eqno(25.4)
$$
which can be confirmed by setting $\lambda=\sigma_{1/w}^{-1}\gamma$,
$H=y^s$ in the above; and more precisely
$$
y(\partial_x-i\partial_y)
[E(\sigma_{1/w_1}(z),1/w;s)]
=-is E_-(\sigma_{1/w_1}(z),1/w;s)
\left(\jmath(\sigma_{1/w_1},z)/
|\jmath(\sigma_{1/w_1}
,z)|\right)^{-2}.\eqno(25.5)
$$
In particular, we have the functional
equation
$$
s\B{E}_-(z,s)=(1-s)\B{S}(s)\B{E}_-(z,1-s),\eqno(25.6)
$$
with
$$
\B{E}_-(s)=\left(\matrix{\vdots\cr E_-(z,1/w;s)\cr
\vdots}\right)_{w|cd}.\eqno(25.7)
$$
Also, $(25.5)$ implies that
$$
\Gamma(s+1)L(2s,\chi_{cd})E_-(\sigma_{1/w_2}(z),1/w_1;s)
\ll y^{\Re s}+y^{1-\Re s},
\eqno(25.8)
$$
as $y$ tends to infinity while $s$ remains bounded, which means that
the left side is regular for all $s$, too. This is a counterpart of
Lemma 7. 
\medskip
In the region of absolute convergence, we have, by $(25.2)$,
$$
\eqalignno{
&\int_{\varGamma\backslash{\cal H}}\overline{\psi_j^-(z)}
E^*\!\left(z,\txt{1\over2}(1-\alpha)\right)
E_-\!\left(z,1/w;s-\txt{1\over2}\alpha\right)d\mu(z)\cr
=&\sum_{\gamma\in\varGamma_{1/w}\backslash\varGamma}
\int_{\sigma_{1/w}^{-1}\gamma(\varGamma\backslash{\cal
H})}\overline{\psi_j^-(\sigma_{1/w}(z))}
\left({\jmath(\gamma^{-1},\sigma_{1/w}(z)\over
|\jmath(\gamma^{-1},\sigma_{1/w}(z)|}\right)^{-2}
E^*\!\left(\sigma_{1/w}(z),\txt{1\over2}(1-\alpha)\right)\cr
&\hskip 2cm\times
y^{s-{1\over2}\alpha}\left({\jmath(\sigma_{1/w}^{-1}\gamma,
\gamma^{-1}\sigma_{1/w}(z))\over|\jmath(\sigma_{1/w}^{-1}\gamma,
\gamma^{-1}\sigma_{1/w}(z))|}\right)^{-2}d\mu(z)\cr
=&\sum_{\gamma\in\varGamma_{1/w}\backslash\varGamma}
\int_{\sigma_{1/w}^{-1}\gamma(\varGamma\backslash{\cal
H})}\overline{\psi_j^-(\sigma_{1/w}(z))}
\left({\jmath(\sigma_{1/w},z)\over
|\jmath(\sigma_{1/w}, z)|}\right)^2
E^*\!\left(\sigma_{1/w}(z),\txt{1\over2}(1-\alpha)\right)
y^{s-{1\over2}\alpha}d\mu(z)\cr
=&\int_{0}^\infty\!\int_0^1
\overline{(\partial_x-i\partial_y)[\psi_j(\sigma_{1/w}(z))]}
E^*\!\left(\sigma_{1/w}(z),\txt{1\over2}(1-\alpha)\right)
y^{s-{1\over2}\alpha-1}dxdy,
&(25.9)
}
$$
since 
$$
\psi_j^-(\sigma_{1/w}(z))=y(\partial_x-i\partial_y)
[\psi_j(\sigma_{1/w}(z))](\jmath(\sigma_{1/w},z)/
|\jmath(\sigma_{1/w},z)|)^2.\eqno(25.10)
$$ 
We observe then that
$E^*\!\left(\sigma_{1/w}(z),\txt{1\over2}(1-\alpha)\right)$ is even
in $x$ as $(23.7)$ implies, and
$\partial_y[\psi_j(\sigma_{1/w}(z))]$ is odd by $(18.5)$. Hence $(25.9)$
becomes
$$
\eqalignno{
&\int_{\varGamma\backslash{\cal H}}\overline{\psi_j^-(z)}
E^*\!\left(z,\txt{1\over2}(1-\alpha)\right)
E_-\!\left(z,1/w;s-\txt{1\over2}\alpha\right)d\mu(z)\cr
=&-i{\Gamma(s+1,\alpha;\kappa_j)
\over\pi^{s-{1\over2}\alpha}\Gamma\left(s+1-{1\over2}\alpha\right)}
D_j(s,\alpha;1/w),&(25.11)
}
$$
provided absolute convergence holds throughout.
\medskip
We decompose the left side of $(25.11)$ in just the same way as we did
in $(24.1)$, and see, via $(25.8)$, that
$$
\left(1-\alpha^2\right)\Gamma
\left(s+1-\txt{1\over2}\alpha\right)L(2s-\alpha,\chi_{cd})
D_j(s,\alpha;1/w)\eqno(25.12)
$$
are all regular in $s$ and $\alpha$. Also, $(25.11)$ gives, via
$(25.6)$, 
$$
{\Gamma(s+1,\alpha;\kappa_j)\over
\pi^{s-{1\over2}\alpha}
\Gamma\left(s-{1\over2}\alpha\right)}\B{D}_j(s,\alpha)
={\Gamma(2-s,-\alpha;\kappa_j)\over
\pi^{1-s+{1\over2}\alpha}
\Gamma\left(1-s+{1\over2}\alpha\right)}\B{S}\left(s-\txt{1\over2}
\alpha\right)\B{D}_j(1-s,-\alpha),\eqno(25.13)
$$
and in particular
$$
\eqalignno{
D_j(s,\alpha)&=\pi^{2s-\alpha-1}
{\Gamma(2-s,-\alpha;\kappa_j)\over
\Gamma(s+1,\alpha;\kappa_j)}{
\Gamma\left(s-{1\over2}\alpha\right)\over
\Gamma\left(1-s+{1\over2}\alpha\right)}\cr
&\times\sum_{w|cd}
\varphi\left(s-\txt{1\over2}\alpha;
\infty,1/w\right)D_j(1-s,-\alpha;1/w).&(25.14)
}
$$
Hence, by $(24.3)$, we have
$$
\eqalignno{
L(2s-\alpha,\chi_{cd})D_j(s,\alpha)
&={1\over\pi^2}\left({\pi^2\over
cd}\right)^{2s-\alpha}\zeta(2(1-s)+\alpha)
{\Gamma(2-s,-\alpha;\kappa_j)\over\Gamma(s+1,\alpha;\kappa_j)}\cr
&\times\sum_{w|cd}\varphi(w)
\prod_{p|v}
\left({p^{s-{1\over2}\alpha}-p^{1-s+{1\over2}\alpha}}\right)
D_j(1-s,-\alpha;1/w).&(25.15)
}
$$
\medskip
With this, we obtain
\medskip
\noindent
{\bf Lemma 12.} {\it With $\epsilon_j=-1$ as well, the assertions of
Lemmas 10 and 11 hold.\/}
\medskip
This ends our treatment of $L_j$ and $D_j$. 
We omit the discussion of $L_{j,k}$, $D_{j,k}$, for they are
analogous.
\bigskip
\noindent
{\bf 26.} Now we may return to $(21.3)$. Here we shall deal with the
first term on the right, the contribution of real analytic
cusp forms. Its contribution to $I(u,v,w,z;g;b/a)$ is, via $(2.2)$,
$(2.3)$, $(17.10)$, equal to
$$
\eqalignno{
&{2\over a^vb^u(2\pi)^{u-w+1}}\sum_{c|a,d|b}
c^{u+v}d^{{1\over2}(3u+v-w+z)}\sum_j 
R_j\left(\txt{1\over2}(u+v+w+z-1),w+z-1\right)\cr
&\times
L_j\left(\txt{1\over2}(u-v-w+z+1);1/c\right)
{([g]_++\epsilon_j[g]_-)(\kappa_j;u,v,w,z)
\over\cosh\pi\kappa_j}.&(26.1) }
$$
with
$$
R_j(s,\alpha)=\zeta(2s-\alpha)D_j(s,\alpha).\eqno(26.2)
$$
By Lemmas 9--12, we see readily that the expression $(26.1)$ is
meromorphic over $\B{C}^4$, and especially in the vicinity of
$p_{1\over2}$ it is regular; the necessary facts about $[g]_\pm$ is
to be given shortly. Hence
its value at
$p_{1\over2}$ equals
$$
\eqalignno{
&{1\over \pi\sqrt{ab}}\sum_{c|a,d|b}
cd\sum_j 
R_j\left(\txt{1\over2},0\right)
L_j\left(\txt{1\over2};1/c\right)
{([g]_++\epsilon_j[g]_-)(\kappa_j;p_{1\over2})
\over\cosh\pi\kappa_j}.&(26.3)
}
$$
We have another contribution of real analytic cusp forms that comes
from $J_-$, which is, however, exactly the same as $(26.3)$.
\medskip
Let us make the last factor in $(26.3)$ explicit. Thus,
comparing $(6.4)$ with [11, $(4.3.13)$--$(4.3.14)$], we see that the
exchange of variables $u$ and $z$ is to be applied 
to [11, Sections 4.6--4.7] to get corresponding identities. More
precisely, we have, under $(3.4)$ and $(4.1)$,
$$
\eqalignno{
[g]_+(r;&\,u,v,w,z)={1\over{4\pi{i}}}\cos(\txt{1\over2}\pi(v-z))
\int_{(\eta_1)}\sin(\txt{1\over2}\pi(u+v+w+z-2s))\cr
&\times\Gamma(\txt{1\over2}(u+v+w+z-1)+ir-s)
\Gamma(\txt{1\over2}(u+v+w+z-1)-ir-s)\cr
&\times\Gamma(s+1-w-z)\Gamma(s+1-v-w) g^*(s,w)ds,&(26.4)\cr
[g]_-(r;&u,v,w,z)= -{1\over{4\pi{i}}}\cosh(\pi r)
\int_{(\eta_1)}\cos(\pi(w+\txt{1\over2}(v+z)-s))\cr
&\times\Gamma(\txt{1\over2}(u+v+w+z-1)+ir-s)
\Gamma(\txt{1\over2}(u+v+w+z-1)-ir-s)\cr
&\times\Gamma(s+1-w-z)\Gamma(s+1-v-w)g^*(s,w)ds,&(26.5)\cr 
}
$$
corresponding to [11, $(4.4.12)$] and [ibid, $(4.4.15)$], respectively. We
then put
$$
\eqalignno{
\Phi_+(\xi;&\,u,v,w,z;g)=-i(2\pi)^{w-v-2}\cos(\txt{1\over2}\pi(v-z))\cr
&\times\int_{-i\infty}^{i\infty}\sin(\txt{1\over2}\pi(u+v+w+z-2s))\cr
&\times\Gamma(\txt{1\over2}(u+v+w+z-1)+\xi-s)
\Gamma(\txt{1\over2}(u+v+w+z-1)-\xi-s)\cr
&\times\Gamma(s+1-w-z)\Gamma(s+1-v-w)g^*(s,w)ds;&(26.6)\cr
\Phi_-(\xi;&\,u,v,w,z;g)=i(2\pi)^{w-v-2}\cos(\pi\xi)
\int_{-i\infty}^{i\infty}\cos(\pi(w+\txt{1\over2}(v+z)-s))\cr
&\times\Gamma(\txt{1\over2}(u+v+w+z-1)+\xi-s)
\Gamma(\txt{1\over2}(u+v+w+z-1)-\xi-s)\cr
&\times\Gamma(s+1-w-z)\Gamma(s+1-v-w)g^*(s,w)ds;&(26.7)\cr
}
$$
and
$$
\eqalignno{
\Xi(\xi;u,v,w,z;g)=&
{1\over{2\pi{i}}}\int_{-\infty i}^{\infty i}
{{\Gamma(\xi+{1\over2}(u+v+w+z-1)-s)}\over{\Gamma(\xi+
{1\over 2}(3-u-v-w-z)+s)}}\cr
&\times\Gamma(s+1-w-z)
\Gamma(s+1-v-w)g^*(s,w)ds.&(26.8)\cr
}
$$
The paths in $(26.6)$ and $(26.7)$ are such that
the poles of the first two gamma-factors and those of the
other three factors in each integrand are separated to the
right and the left, respectively, by the path, and $\xi,u,v,w,z$ are
assumed to be such that the path can be drawn. 
The path in $(26.8)$ separates the poles of $\Gamma(\xi+{1\over
2}(u+v+w+z-1)-s)$ and those of
$\Gamma(s+1-w-z)\Gamma(s+1-v-w)g^*(s,w)$ to the left and the right
of the path, respectively. We have the relations
$$\eqalignno{
\Phi_{+}(\xi;u,v,w,z;g)&=-{{(2\pi)^{w-u}\cos\big({1\over
2}\pi(v-z))}\over{4\sin(\pi\xi)}}\cr
&\times\{\Xi(\xi;u,v,w,z;g)-\Xi(-\xi;u,v,w,z;g)\},
&(26.9)\cr
\Phi_{-}(\xi;u,v,w,z;g)&={(2\pi)^{w-u}\over{4\sin(\pi\xi)}}
\{\sin(\pi(\txt{1\over2}(u-w)+\xi))\Xi(\xi;u,v,w,z;g)\cr
&-\sin(\pi({\txt{1\over
2}}(u-w)-\xi))\Xi(-\xi;u,v,w,z;g)\},&(26.10)\cr
}
$$
provided the left sides are well-defined. 
\medskip
Under $(4.1)$, we can obviously take
$(\eta_1)$ as the contours in the last three integrals; and
we have, for $r\in\B{R}$,
$$
\eqalignno{
[g]_+(r;u,v,w,z)&={1\over2}(2\pi)^{1+u-w}
\Phi_{+}(ir;u,v,w,z;g),\cr
[g]_-(r;u,v,w,z)&={1\over2}(2\pi)^{1+u-w}
\Phi_{-}(ir;u,v,w,z;g).
&(26.11)
}
$$
In particular, we have, after continuation,
$$
\eqalignno{
[g]_+(r;p_{1\over2})&=-{\pi\over4\sin(\pi ir)}
\left(\Xi\left(ir;p_{1\over2};g\right)-
\Xi\left(-ir;p_{1\over2};g\right)\right),\cr
[g]_-(r;p_{1\over2})&={\pi\over4}
\left(\Xi\left(ir;p_{1\over2};g\right)+
\Xi\left(-ir;p_{1\over2};g\right)\right),
&(26.12)
}
$$
and
$$
([g]_++\epsilon_j[g]_-)\left(r;p_{1\over2}\right)={\pi\over2}
\Re\left\{\left(\epsilon_j+{i\over\sinh\pi
r}\right)\Xi\left(ir;p_{1\over2};g\right)\right\},\eqno(26.13)
$$
since $(3.2)$ and $(26.8)$ imply $\overline
{\Xi\left(ir;p_{1\over2};g\right)}=\Xi\left(-ir;p_{1\over2};g\right)$.
\medskip
From this, we get immediately
\medskip
\noindent
{\bf Lemma 13.} {\it Provided the polynomial $A$ is supported by the
set of square-free integers,
the contribution of real analytic cusp forms to $M_2(g;A)$ is
equal to
$$
\sum_{c,\,d} \e{A}(c,d)
\e{C}(c,d;g),\eqno(26.14)
$$
where
$$
\e{A}(c,d)=\sum_{(ac,bd)=1}{\alpha_{acl}
\overline{\alpha_{bdl}}\over abl},\eqno(26.15)
$$
and
$$
\eqalignno{
\e{C}(c,d;g)=\sum_{\scr{j}\atop\scr{\kappa_j^2+\txt{1\over4}\in{\rm
Sp}(\varGamma_0(cd))}}&{1\over\cosh\pi\kappa_j} 
R_j\left(\txt{1\over2},0\right)
L_j\left(\txt{1\over2};1/c\right)\cr
&\times\Re\left\{\left(\epsilon_j+{i\over\sinh\pi\kappa_j}\right)
\Xi\left(i\kappa_j;p_{1\over2};g\right)\right\}.&(26.16)
}
$$
}
\medskip
\noindent
The fact that the parity symbol $\epsilon_j$ appears in this way will
turn out to be crucial in our later discussion of a certain non-vanishing
assertion (Sections 31--36).
\medskip
The contribution of holomorphic cusp forms is analogous,
and we may skip it.
\bigskip
\noindent
{\bf 27.} We turn to the contribution of continuous spectrum; and we
see from $(21.3)$ that we need first to consider the sum
$$
\eqalignno{
&\sum_n{\sigma_{-2ir}(n;\chi_{cd})\over
n^s}\prod_{p|(c_1,c)(d_1,d)}
\left\{\sigma_{-2ir}(n_p)\left(1-{1\over
p^{1+2ir}}\right)-1\right\}\cr
=&\sum_{l|(c_1,c)(d_1,d)}\mu((c_1,c)(d_1,d)/l)
\prod_{p|l}\left(1-{1\over
p^{1+2ir}}\right)\sum_n\sigma_{-2ir}(n;\chi_{cd})
\sigma_{-2ir}(n_l)n^{-s},
&(27.1)
}
$$
with $n_l=(n,l^\infty)$. We have
$$
\eqalignno{
\sum_n\sigma_{-2ir}(n;\chi_{cd})
\sigma_{-2ir}(n_l)n^{-s}&=\left\{\sum_{(n,l)=1}
{\sigma_{-2ir}(n;\chi_{cd})\over n^s}\right\}
\left\{\sum_{n|l^\infty}
{\sigma_{-2ir}(n)\over n^s}\right\}\cr
&=L(s,\chi_l)L(s+2ir,\chi_{cd})
\prod_{p|l}\left(1-{1\over p^s}\right)^{-1}
\left(1-{1\over p^{s+2ir}}\right)^{-1}\cr
&=\zeta(s)L(s+2ir,\chi_{cd})\prod_{p|l}
\left(1-{1\over p^{s+2ir}}\right)^{-1}.&(27.2)
}
$$
Thus
$$
\eqalignno{
&\sum_n{\sigma_{-2ir}(n;\chi_{cd})\over
n^s}\prod_{p|(c_1,c)(d_1,d)}
\left\{\sigma_{-2ir}(n_p)\left(1-{1\over
p^{1+2ir}}\right)-1\right\}\cr
=&\zeta(s)L(s+2ir,\chi_{cd})
\sum_{l|(c_1,c)(d_1,d)}\mu((c_1,c)(d_1,d)/l)
\prod_{p|l}\left(1-{1\over p^{s+2ir}}\right)^{-1}
\left(1-{1\over p^{1+2ir}}\right)\cr
=&\zeta(s)L(s+2ir,\chi_{cd})\prod_{p|(c_1,c)(d_1,d)}
\left\{\left(1-{1\over p^{s+2ir}}\right)^{-1}
\left(1-{1\over p^{1+2ir}}\right)-1\right\}.
&(27.3)
}
$$
\medskip
Next, we need to treat 
$$
\eqalignno{
&\sum_n{\sigma_{2ir}(n;\chi_{cd})\sigma_\alpha(n)
\over n^s}\prod_{p|c_1}\left\{\sigma_{2ir}(n_p)
\left(1-{1\over p^{1-2ir}}\right)-1\right\}\cr
=&\sum_{l|c_1}\mu(c_1/l)\prod_{p|l}
\left(1-{1\over p^{1-2ir}}\right)\sum_n
{\sigma_{2ir}(n;\chi_{cd})\sigma_\alpha(n)\sigma_{2ir}(n_l)
\over n^s}.
&(27.4)
}
$$
We have
$$
\sum_n{\sigma_{2ir}(n;\chi_{cd})\sigma_\alpha(n)\sigma_{2ir}(n_l)
\over n^s}=\left\{\sum_{(n,l)=1}{\sigma_{2ir}(n;\chi_{cd})
\sigma_\alpha(n)
\over n^s}\right\}
\left\{\sum_{n|l^\infty}{\sigma_{2ir}(n)\sigma_\alpha(n)
\over n^s}\right\}.\eqno(27.5)
$$
Analogously to a famous formula of Ramanujan, we have
$$
\eqalignno{
&\sum_{(n,l)=1}{\sigma_{2ir}(n;\chi_{cd})\sigma_\alpha(n)
\over n^s}\cr
=&\sum_{(n,cd)=1}{\sigma_{2ir}(n)\sigma_\alpha(n)
\over n^s}\sum_{n|(cd/l)^\infty}{\sigma_\alpha(n)
\over n^s}\cr
=&{L(s,\chi_q)L(s-2ir,\chi_q)
L(s-\alpha,\chi_{cd})L(s-2ir-\alpha,\chi_{cd})
\over L(2s-2ir-\alpha,\chi_{cd})}\prod_{p|cd/l}
\left(1-{1\over p^s}\right)^{-1}\left(1-{1\over
p^{s-\alpha}}\right)^{-1},\cr
&\sum_{n|l^\infty}{\sigma_{2ir}(n)\sigma_\alpha(n)
\over n^s}\cr
=&\prod_{p|l}{\displaystyle{1-{1\over p^{2s-2ir-\alpha}}}\over
\displaystyle{\left(1-{1\over p^s}\right)
\left(1-{1\over p^{s-2ir}}\right)
\left(1-{1\over p^{s-\alpha}}\right)
\left(1-{1\over p^{s-2ir-\alpha}}\right)}}.
 &(27.6)
}
$$
Thus,
$$
\eqalignno{
\sum_n{\sigma_{2ir}(n;\chi_{cd})\sigma_\alpha(n)\sigma_{2ir}(n_l)
\over n^s}&={\zeta(s)L(s-2ir,\chi_{cd})
\zeta(s-\alpha)L(s-2ir-\alpha,\chi_{cd})
\over L(2s-2ir-\alpha,\chi_{cd})}\cr
&\times\prod_{p|l}{\displaystyle{1-{1\over p^{2s-2ir-\alpha}}}\over
\displaystyle{\left(1-{1\over p^{s-2ir}}\right)
\left(1-{1\over p^{s-2ir-\alpha}}\right)}}.
&(27.7)
}
$$
Hence,
$$
\eqalignno{
\sum_n&{\sigma_{2ir}(n;\chi_{cd})\sigma_\alpha(n)
\over n^s}\prod_{p|c_1}\left\{\sigma_{2ir}(n_p)
\left(1-{1\over p^{1-2ir}}\right)-1\right\}\cr
=&{\zeta(s)L(s-2ir,\chi_q)
\zeta(s-\alpha)L(s-2ir-\alpha,\chi_{cd})
\over L(2s-2ir-\alpha,\chi_{cd})}\cr
&\times\prod_{p|c_1}\left\{
{\displaystyle{\left(1-{1\over p^{1-2ir}}\right)
\left(1-{1\over
p^{2s-2ir-\alpha}}\right)}\over
\displaystyle{\left(1-{1\over p^{s-2ir}}\right)
\left(1-{1\over p^{s-2ir-\alpha}}\right)}}-1\right\}\cr
=&{\zeta(s)\zeta(s-2ir)\zeta(s-\alpha)\zeta(s-2ir-\alpha)
\over \zeta(2s-2ir-\alpha)}\prod_{p|cd}
\left(1-{1\over p^{2s-2ir-\alpha}}\right)^{-1}\cr
&\times\prod_{p|d_1}
\left(1-{1\over p^{s-2ir}}\right)
\left(1-{1\over p^{s-2ir-\alpha}}\right)\cr
&\times\prod_{p|c_1}
\left\{\left(1-{1\over p^{1-2ir}}\right)
\left(1-{1\over p^{2s-2ir-\alpha}}\right)
-\left(1-{1\over p^{s-2ir}}\right)
\left(1-{1\over p^{s-2ir-\alpha}}\right)\right\}.
&(27.8)
}
$$
\medskip
\noindent
{\bf 28.} Under the conditions $(3.4)$, $(4.1)$ and by $(21.3)$,
$(27.3)$, $(27.8)$, the contribution of the continuous spectrum to $I$
via $J_+^*$ is equal to
$$
4{(2\pi)^{w-u-2}\over a^v b^u}\int_{-\infty}^\infty
{Y_{a,b}(ir;u,v,w,z)
Z(ir;u,v,w,z)\over\zeta(1+2ir)\zeta(1-2ir)}
([g]_-+[g]_-)(r;u,v,w,z)dr
\eqno(28.1)
$$
where
$$
\eqalignno{
Z(\xi;u,v,w,z)&=\zeta\left(\txt{1\over2}(u+v+w+z-1)+\xi\right)
\zeta\left(\txt{1\over2}(u+v+w+z-1)-\xi\right)\cr
&\times\zeta\left(\txt{1\over2}(u+v-w-z+1)+\xi\right)
\zeta\left(\txt{1\over2}(u+v-w-z+1)-\xi\right)\cr
&\times\zeta\left(\txt{1\over2}(u-v-w+z+1)+\xi\right)
\zeta\left(\txt{1\over2}(u-v-w+z+1)-\xi\right)&(28.2)
}
$$
and
$$
Y_{a,b}(\xi;u,v,w,z)=\sum_{c|a, d|b}c^{u+v}
d^{{1\over2}(3u+v-w+z-1)-\xi}X_{cd}(\xi;u,v,w,z),\eqno(28.3)
$$
with
$$
\eqalignno{
&X_{cd}(\xi;u,v,w,z)
=\prod_{p|cd}\left\{
\left(1-{1\over p^{1+2\xi}}\right)
\left(1-{1\over p^{1-2\xi}}\right)
\left(1-{1\over p^{u+v}}\right)\right\}^{-1}\cr
&\times\sum_{cd=c_1d_1}
{1\over d_1}\left({(d_1,d)\over(c_1,c)}\right)^{{1\over2}+\xi}\cr
&\times\prod_{p|(d_1,c)(c_1,d)}\left(1-{1\over
p^{{1\over2}(u-v-w+z+1)+\xi}}\right)
\prod_{p|(c_1,c)(d_1,d)}
\left({1\over p^{{1\over2}(u-v-w+z)}}-{1\over
p^{{1\over2}+\xi}}\right)\cr 
&\times\prod_{p|d_1}\left(1-{1\over p^{{1\over2}(u+v+w+z-1)-\xi}}\right)
\left(1-{1\over p^{{1\over2}(u+v-w-z+1)-\xi}}\right)\cr
&\times\prod_{p|c_1}\Bigg\{
\left(1-{1\over p^{1-2\xi}}\right)\left(1-{1\over p^{u+v}}\right)\cr
&\hskip 2cm-
\left(1-{1\over p^{{1\over2}(u+v+w+z-1)-\xi}}\right)
\left(1-{1\over p^{{1\over2}(u+v-w-z+1)-\xi}}\right)
\Bigg\}.&(28.4)
}
$$
One may carry out the last sum and transform $X_{cd}$ and thus
$Y_{ab}$ into a more closed expression that is a product over prime
divisors of $ab$; however, for our aim it does not seem particularly
expedient to do so, and we leave $(28.3)$ as it is.
\medskip
To continue $(28.1)$ to a neighbourhood of $p_{1\over2}$, 
we need to shift the contour rightward and leftward appropriately as is
done in [11, Section 4.7], and there appears a residual contribution,
which will be treated in detail later. Here
we shall compute, at $p_{1\over2}$,  the integral thus continued.
\medskip
By $(28.4)$, we have, for $r\in\B{R}$,
$$
\eqalignno{
&X_{cd}(ir;p_{1\over2})\cr
=&\prod_{p|cd}{\displaystyle{1-{1\over
p^{{1\over2}+ir}}}\over
\displaystyle{\left|1-{1\over p^{1+2ir}}\right|^2
\left(1-{1\over p}\right)}}\sum_{cd=c_1d_1}
{1\over d_1}\left({(d_1,d)\over(c_1,c)}\right)^{{1\over2}+ir} 
\prod_{p|d_1}\left(1-{1\over p^{{1\over2}-ir}}\right)^2\cr
&\hskip1cm\times\prod_{p|c_1}\Bigg\{
\left(1-{1\over p^{1-2ir}}\right)\left(1-{1\over p}\right)
-\left(1-{1\over p^{{1\over2}-ir}}\right)^2\Bigg\}\cr
=&\prod_{p|cd}{\displaystyle{1-{1\over
p^{{1\over2}+ir}}}\over
\displaystyle{\left|1-{1\over p^{1+2ir}}\right|^2
\left(1-{1\over p}\right)}}\prod_{p|cd}\Bigg\{
{(p,d)^{{1\over2}+ir}\over
p}
\left(1-{1\over p^{{1\over2}-ir}}\right)^2\cr
&\hskip1cm+{1\over(p,c)^{{1\over2}+ir}}\Bigg\{
\left(1-{1\over p^{1-2ir}}\right)\left(1-{1\over p}\right)
-\left(1-{1\over p^{{1\over2}-ir}}\right)^2\Bigg\}\cr
=&c^{-{1\over2}-ir}\prod_{p|cd}{\displaystyle{1-{1\over
p^{{1\over2}+ir}}}\over
\displaystyle{\left|1-{1\over p^{1+2ir}}\right|^2
\left(1-{1\over p}\right)}}\Bigg\{
\left(1-{1\over p^{1-2ir}}\right)\left(1-{1\over p}\right)
-\left(1-{1\over p^{{1\over2}-ir}}\right)^3\Bigg\}\cr
=&c^{-{1\over2}-ir}\prod_{p|cd}
\left|1+{1\over p^{{1\over2}+ir}}\right|^{-2}
\left(1-{1\over p}\right)^{-1}\Bigg\{
\left(1+{1\over p^{{1\over2}-ir}}\right)\left(1-{1\over p}\right)
-\left(1-{1\over p^{{1\over2}-ir}}\right)^2\Bigg\}.\quad&(28.5)
}
$$
This implies that
$$
\eqalignno{
&Y_{a,b}(ir;p_{1\over2})\cr
=&\sum_{c|a,d|b}(cd)^{{1\over2}-ir}
\prod_{p|cd}
\left|1+{1\over p^{{1\over2}+ir}}\right|^{-2}
\left(1-{1\over p}\right)^{-1}\Bigg\{
\left(1+{1\over p^{{1\over2}-ir}}\right)\left(1-{1\over p}\right)
-\left(1-{1\over p^{{1\over2}-ir}}\right)^2\Bigg\}\cr
=&{ab\over\varphi(ab)}\prod_{p|ab}\left(4
\left|1+{1\over p^{{1\over2}+ir}}\right|^{-2}-{1\over p}\right).
&(28.6)
}
$$
\medskip
We have obtained
\medskip
\noindent
{\bf Lemma 14.} {\it Provided the polynomial $A$ is supported
by the set of square-free integers, the
contribution of continuous spectrum to
$M_2(g;A)$ is equal to
$$
\eqalignno{
{1\over\pi}&\sum_{(a,b)=1}
{\alpha_{al}\overline{\alpha_{bl}}\over\varphi(ab)l}
\int_{-\infty}^\infty
{\left|\zeta\left({1\over2}+ir\right)\right|^6\over|\zeta(1+2ir)|^2}\cr
&\times\prod_{p|ab}\left(4
\left|1+{1\over p^{{1\over2}+ir}}\right|^{-2}-{1\over p}\right)
\Re\left\{\left(1+{i\over\sinh\pi r}\right)
\Xi\left(ir;p_{1\over2};g\right)\right\}dr.&(28.7)
}
$$
}
\bigskip
\noindent
{\bf 29.} We shall give the continuation procedure of $(28.1)$ to a
neighbourhood of
$p_{1\over2}$. This is, however, analogous to that
pertaining to the pure fourth moment $M_2(g;1)$ that is developed in
[11, Sections 4.6--4.7]; and we can be brief.
\medskip
By $(26.9)$--$(26.11)$, we transform $(28.1)$ into
$$
\eqalignno{
i{(2\pi)^{w-u-1}\over a^vb^u}\int_{(0)}&
{{Z(\xi;u,v,w,z)}
\over{\sin(\pi\xi)\zeta(1+2\xi)\zeta(1-2\xi)}}
\left\{Y_{a,b}(\xi;u,v,w,z)+Y_{a,b}(-\xi;u,v,w,z)\right\}\cr
\times&\left\{\cos(\txt{1\over2}\pi(v-z))
-\sin(\pi({\txt{1\over2}}(u-w)+\xi))\right\}
\Xi(\xi;u,v,w,z;g)d\xi;&(29.1)\cr
}
$$
and applying the functional equation for $\zeta$ to $\zeta(1-2\xi)$, this
becomes
$$
\eqalignno{
2i{(2\pi)^{w-u-2}\over a^vb^u}\int_{(0)}&
{(2\pi)^{2\xi}\Gamma(1-2\xi)
Z(\xi;u,v,w,z)\over \zeta(2\xi)
\zeta(1+2\xi)}\left\{Y_{a,b}(\xi;u,v,w,z)+Y_{a,b}(-\xi;u,v,w,z)\right\}
\cr
&\times\left \{\cos({\txt{1\over2}}\pi(v-z))
-\sin(\pi({\txt{1\over2}}(u-w)+\xi))\right\}
\Xi(\xi;u,v,w,z;g)d\xi&(29.2)\cr
}
$$ 
(see [11, $(4.6.14)$--$(4.6.15)$]). We shift the last contour to
the far right, and we obtain a meromorphic continuation to a
domain containing the point $p_{1\over2}$; then, restricting
ourselves to the vicinity of $p_{1\over2}$, we shift the contour
back to the imaginary axis. The resulting integral has been considered
already in the last section.
\medskip
The residual contribution of the last procedure takes place when
$$
\eqalignno{
\xi_1=\;&{1\over2}(u+v+w+z-3),\quad\xi_2=\;{1\over2}(u-v-w+z-1),\cr
\xi_3=\;&{1\over2}(3-u-v-w-z),\quad\xi_4=\;{1\over2}(u+v-w-z-1).&(29.3)
}
$$
(see [11,$(4.6.16)$] and the bottom lines of [ibid, p.\ 173]). It should be
stressed that this assertion depends on the fact that the singularities,
save for those belonging to $Z(\xi;u,v,w,z)$,
that we encounter in this procedure are independent of the
location of $(u,v,w,z)$; especially those of $Y_{a,b}(\pm\xi;u,v,w,z)$
come only from the first product on the right of $(28.4)$ and
are independent of $(u,v,w,z)$.
\medskip
\noindent
{\csc Remark 4}.  However, one should note that 
the set of poles of $Y_{a,b}(\xi;u,v,w,z)$ as a function of $\xi$ cluster
at the point $\xi=\pm{1\over2}$ if $a,b$ are allowed to vary
arbitrarily. Thus, if the length of the polynomial $A$ increases
indefinitely, then the nature of the main term of
$M_2(g;A)$ should become subtler.
\bigskip
\noindent
{\bf 30.} With this, we have essentially finished spectrally
decomposing $M_2(g;A)$. Although we have not yet computed the main
term explicitly, the above is already quite adequate to analyze
the error term in the asymptotic formula for the unweighted mean
$$
\int_0^T\left|\zeta\left(\txt{1\over2}+it\right)\right|^4
\left|A\left(\txt{1\over2}+it\right)\right|^2dt.\eqno(30.1)
$$
With this in mind, we shall investigate the location of poles of
the Mellin transform $Z_2(s;A)$, focusing our attention to the 
contribution of real analytic cusp forms, for the relevant part 
of $Z_2(s;A)$ seems to be the most interesting. 
\medskip
Having the assertion of Lemma 13, the argument of [11, Section 5.3] 
works with $Z_2(s;A)$ as well without any essential change. We find,
on the assumption on eigenvalues $\kappa_j^2+{1\over4}$ made
in the introduction, that
\medskip
\noindent
{\bf Lemma 15.} {\it The function $Z_2(s;A)$ is meromorphic over the
entire complex plane. It has a pole of the fifth order at $s=1$; and all
other poles are in the half plane $\Re s\le{1\over2}$. More precisely,
$Z_2(s;A)$ has a pole at ${1\over2}+i\kappa$, $\kappa>0$,
if and only if it holds that
$$
\sum_{c,d}\e{A}(c,d)\sum_{\scr{\kappa_j=\kappa}\atop\scr{\kappa_j^2
+\txt{1\over4}\in{\rm
Sp}(\varGamma_0(cd))}}R_j\left(\txt{1\over2},0\right)
L_j\left(\txt{1\over2};1/c\right)\left(\epsilon_j-{i\over
\sinh\pi\kappa}\right)\ne0.\eqno(30.2)
$$
}
\medskip
We are going to show that if $A$ is fixed besides a natural condition
on its coefficients, then $(30.2)$ holds for infinitely many $\kappa$. 
To this end we shall establish in the sequel that there are infinitely
many $\kappa$ such that
$$
\e{R}(\kappa;A)=\sum_{c,d}\e{A}(c,d)\sum_{\scr{\kappa_j
=\kappa}\atop\scr{\kappa_j^2+\txt{1\over4}\in{\rm
Sp}(\varGamma_0(cd))}}\epsilon_jR_j\left(\txt{1\over2},0\right)
L_j\left(\txt{1\over2};1/c\right)\ne0.\eqno(30.3)
$$
\medskip
\noindent
{\csc Remark 5}. As to the possible poles coming from the contribution of
the continuous spectrum, one may follow the discussion in [11, p.\ 211].
In view of $(28.7)$, we may have poles at 
$$
(2l+1){\pi i\over \log p},\quad l\in\B{Z},\eqno(30.4)
$$
where $p|ab$ with $\alpha_a\alpha_b\ne0$. Thus it can be asserted,
somewhat informally, that
as the length of $A$ tends to infinity the imaginary axis is
gradually filled up with poles of $Z_2(s;A)$.
\bigskip
\noindent
{\bf 31.} To deal with $\e{R}(\kappa; A)$,
we adopt the argument of [11, Section 3.3]. Thus,
on noting the definitions $(17.1)$ and $(26.2)$, we consider 
more generally the sum
$$
\eqalignno{
\e{D}(u,v;h)=&\zeta(u+v)\sum_j
\varrho_j(-f;1/c)D_j(u,u-v){h(\kappa_j)
\over\cosh\pi\kappa_j}\cr
=&\zeta(u+v)\e{D}_1(u,v;h),&(31.1)
}
$$
with an integer $f>0$, where the sum is extended over 
$\kappa_j^2+\txt{1\over4}\in{\rm Sp}(\varGamma_0(cd))$ with a
fixed pair $c, d$, $\mu(cd)\ne0$; also the weight $h$ is assumed
to be an even, entire function such that
$$
h\left(\pm\txt{1\over2}i\right)=0\eqno(31.2)
$$
and
$$
h(r)\ll\exp(-c_0|r|^2),\eqno(31.3)
$$
with a certain $c_0>0$, in any fixed horizontal strip. By Lemmas
10--12, $\e{D}(u,v;h)$
is meromorphic over $\B{C}^2$, and regular in the vicinity of
$\left({1\over2},{1\over2}\right)$; in particular, we have
$$
\e{D}\left(\txt{1\over2},\txt{1\over2};h\right)=
\sum_j\epsilon_j\varrho_j(f;1/c)R_j\left(\txt{1\over2},0\right)
{h(\kappa_j)\over\cosh\pi\kappa_j}.
\eqno(31.4)
$$
\medskip
In the region of absolute convergence, we have, by definition,
$$
\e{D}_1(u,v;h)=\sum_m {m}^{-u}
\sigma_{u-v}(m)\sum_j
\varrho_j(-f;1/c)\varrho_j(m;\infty)
{h(\kappa_j)\over\cosh\pi\kappa_j}.\eqno(31.5)
$$
We apply $(21.1)$ to the inner sum, getting
$$
\e{D}_1(u,v;h)=\e{D}_2(u,v;h)+\e{D}_3(u,v;h)\eqno(31.6)
$$
where
$$
\e{D}_2(u,v;h)={1\over c\sqrt{d}}\sum_m {m}^{-u}\sigma_{u-v}(m)
\sum_{\scr{l}\atop\scr{(l,d)=1}}
{1\over l}S(m,-\overline{d}f;cl)
\psi\left({{4\pi}\over{cl\sqrt{d}}}\sqrt{mf}\right),
\eqno(31.7)
$$
with
$$
\psi(x)={4\over\pi^2}\int_{-\infty}^\infty{r}
\sinh(\pi{r}){K}_{2ir}(x)h(r)dr,\eqno(31.8)
$$
and
$$
\eqalignno{
L(u+v&,\chi_{cd})\e{D}_3(u,v;h)=-{1\over\pi
}\sum_{c_1d_1=cd}{1\over d_1}
\int_{-\infty}^\infty\left({(d_1,d)\over\sqrt{d}}\right)^{1+2ir}
{f^{ir}\sigma_{-2ir}(f;\chi_{cd})\over|L(1+2ir,\chi_{cd})|^2}\cr
&\times\zeta(u+ir)\zeta(u-ir)\zeta(v+ir)\zeta(v-ir)
\cr
&\times \prod_{p|(c_1,c)(d_1,d)}\left(\sigma_{-2ir}(f_p)
\left(1-{1\over p^{1+2ir}}\right)-1\right)
\prod_{p|d_1}
\left(1-{1\over p^{u-ir}}\right)
\left(1-{1\over p^{v-ir}}\right)\cr
&\times\prod_{p|c_1}
\left\{\left(1-{1\over p^{1-2ir}}\right)
\left(1-{1\over p^{u+v}}\right)
-\left(1-{1\over p^{u-ir}}\right)
\left(1-{1\over p^{v-ir}}\right)\right\}h(r)dr,
\quad&(31.9)
}
$$
in which we have used $(27.8)$ with $s=u+ir$, $\alpha=u-v$.
\medskip
\noindent
{\bf 32.} To transform $\e{D}_2$ we use the formula
$$
\psi(x)={1\over\pi^2}\int_{(\alpha)}{{\hat{h}(s)}\over{\cos\pi{s}}}
\Big({x\over2}\Big)^{-2s}ds,\qquad 0<\alpha<\txt{1\over2}\,,
\eqno(32.1)
$$
where
$$
\hat{h}(s)=\int_{-\infty}^\infty{r}h(r)
{{\Gamma(s+ir)}\over{\Gamma(1-s+ir)}}dr
\eqno(32.2)
$$
(see [11, p.\ 113]).
Moving the last path far down, we see that $\hat{h}$ is entire. Also we
have
$$
\hat{h}(\pm\txt{1\over2})=0,
\eqno(32.3)
$$
and $(32.1)$ is replaced by
$$
\psi(x)={1\over\pi^2}\int_{(\alpha)}{{\hat{h}(s)}
\over{\cos\pi{s}}}\Big({x\over2}
\Big)^{-2s}ds,\qquad -\txt{3\over2}<\alpha<\txt{3\over2}.\eqno(32.4)
$$
The integrand decays exponentially, which facilitate our discussion
greatly. We stress that the presence of the factor
$\epsilon_j$ in $(30.3)$ has induced this effect.   
\medskip
Thus in $(31.7)$ we have
$$
\eqalignno{
\sum_{\scr{l}\atop\scr{(l,d)=1}}&
{1\over l}S(m,-\overline{d}f;cl)
\psi\left({{4\pi}\over{cl\sqrt{d}}}\sqrt{mf}\right),\cr
&={1\over\pi^2}\sum_{\scr{l}\atop\scr{(l,d)=1}}{1\over{l}}
S(m,-\overline{d}f;cl)\int_{(\alpha)}
{{\hat{h}(s)}\over{\cos\pi{s}}}
\left({{2\pi}\over{cl\sqrt{d}}}\sqrt{mf}\right)^{-2s}ds,&(32.5)\cr
}
$$
with
$$
-\txt{3\over2}<\alpha<-\txt{1\over4}.\eqno(32.6)
$$
The right side of $(32.5)$ converges absolutely. 
Then we assume that 
$$
\Re{u},\,\Re{v}>1-\alpha.
\eqno(32.7)
$$
On this we insert $(32.5)$ into $(31.7)$, and get
$$
\e{D}_2(u,v;h)={1\over\pi^2c\sqrt{d}}
\sum_{\scr{l}\atop\scr{(l,d)=1}}{1\over l}P(u,v;l),
\eqno(32.8)
$$
where
$$
\eqalignno{
&P(u,v;l)=\int_{(\alpha)}\Big({{2\pi}\over{cl\sqrt{d}}}
\sqrt{f}\Big)^{-2s}{{\hat{h}(s)}\over{\cos(\pi{s})}}\cr
&\times \left\{\sum_{\scr{a=1}\atop\scr{(a,cl)=1}}^{cl}
\exp(-2\pi i\overline{d}fa/cl)
\sum_{m=1}^\infty\sigma_{u-v}(m)
\exp(2\pi im\overline{a}/cl)m^{-u-s}\right\}ds, &(32.9)
}
$$
with $a\overline{a}\equiv1\bmod{cl}$. 
\medskip
We introduce further a sub-region of $(32.7)$:
$$
1-\alpha<\Re(u),\Re(v)<-\beta, \quad
-\txt{3\over2}<\beta<\alpha-1,
\quad-\txt{1\over2}<\alpha<-\txt{1\over4}.
\eqno(32.10)
$$
Then we move the path in $(32.9)$ to $(\beta)$. On the assumption
$u\ne v$, we have, by Estermann's functional equation (see [11, Lemma
3.7]),
$$
\eqalignno{
P(u,v;l)&=-2\pi{i}c_{cl}(f)(cl)^{1-u-v}\Big\{(2\pi\sqrt{f/d})^{2(u-1)}
\hat{h}(1-u)\zeta(1-u+v)/\cos\pi{u}\cr
&+(2\pi\sqrt{f/d})^{2(v-1)}
\hat{h}(1-v)\zeta(1-v+u)/\cos\pi{v}\Big\}\cr
&+2(2\pi)^{u+v-2}(cl)^{1-u-v}
\Big\{\sum_{m=1}^\infty{m}^{u-1}\sigma_{v-u}(m)c_{cl}(dm+f)
\Psi_+(u,v;dm/f;h)\cr
&+\sum_{m=1}^\infty{m}^{u-1}\sigma_{v-u}(m)c_{cl}(dm-f)
\Psi_-(u,v;dm/f;h)\Big\},\qquad&(32.11)\cr
}
$$
where
$$
\Psi_+(u,v;x;h)=-\int_{(\beta)}\Gamma(1-u-s)\Gamma(1-v-s)
\cos\left(\pi\left(s+\txt{1\over2}(u+v)\right)\right)
{{\hat{h}(s)}\over{\cos\pi{s}}}x^sds
\eqno(32.12)
$$
and
$$
\Psi_-(u,v;x;h)=\cos\left(\txt{1\over2}\pi(u-v)\right)\int_{(\beta)}
\Gamma(1-u-s)\Gamma(1-v-s)
{{\hat{h}(s)}\over{\cos(\pi{s})}}x^sds.\eqno(32.13)
$$
\bigskip
\noindent
{\bf 33.} We insert $(32.11)$ into $(32.8)$. We get under $(32.10)$ that
$$
L(u+v,\chi_{cd})\e{D}_2(u,v;h)=\left\{\e{D}_2^1+\e{D}_2^2+\e{D}_2^3
+\e{D}_2^4\right\}(u,v;h),\eqno(33.1)
$$
where
$$
\eqalignno{
\e{D}_2^1&={2\over\pi
i\sqrt{d}}\left\{\left(\!2\pi\sqrt{f\over
d}\right)^{2(u-1)}\!{{\hat{h}(1-u)}\over{\cos\pi{u}}} \zeta(1-u+v)
+\left(\!2\pi\sqrt{f\over
d}\right)^{2(v-1)}\!{{\hat{h}(1-v)}\over{\cos\pi{v}}}\zeta(1-v+u)
\right\}\cr
&\hskip
1cm\times\sigma_{1-u-v}(f,\chi_{cd})
\prod_{p|c}\left(\!\sigma_{1-u-v}(f_p)
\left(1-{1\over p^{u+v}}\right)-1\right)\cr
\e{D}_2^2&=8{(2\pi)^{u+v-4}\over\sqrt{d}}\sum_m m^{u-1}
\sigma_{v-u}(m)\sigma_{1-u-v}(dm+f;\chi_{cd})\Psi_+(u,v;dm/f;h)\cr
&\hskip 1cm\times\prod_{p|c}\left(\!\sigma_{1-u-v}((dm+f)_p)
\left(1-{1\over p^{u+v}}\right)-1\right),\cr
\e{D}_2^3&=8{(2\pi)^{u+v-4}\over\sqrt{d}}
\sum_{\scr{m}\atop\scr{dm\ne f}} m^{u-1}
\sigma_{v-u}(m)\sigma_{1-u-v}(dm-f;\chi_{cd})\Psi_-(u,v;dm/f;h)\cr
&\hskip 1cm\times\prod_{p|c}\left(\!\sigma_{1-u-v}((dm-f)_p)
\left(1-{1\over p^{u+v}}\right)-1\right),\cr
\e{D}_2^4&=8(2\pi)^{u+v-4}{\varphi(c)\over c^{u+v}\sqrt{d}}
L(u+v-1,\chi_d)(f/d)^{u-1}\sigma_{v-u}(f/d)\Psi_-(u,v;1;h),&(33.2)
}
$$
in which $\e{D}_2^4$ appears only when $d|f$. 
\medskip
The expansion $(33.1)$ with $(33.2)$ has been proved under the assumption
that $u\ne v$ and $(32.10)$ holds. However, the former can be dropped now;
and also $\e{D}_2^2$ and $\e{D}_2^3$ converge absolutely if 
$1+\beta<\Re u,\,\Re v<-\beta$. In particular, 
$L(u+v,\chi_{cd})\e{D}_2(u,v;h)$ is regular
at $\left({1\over2},{1\over2}\right)$, and there $(33.1)$ holds.
\medskip
Further, shifting the path in $(31.9)$ upward and downward
appropriately, we have
the following continuation of $\e{D}_3$ to the domain
$\Re{u},\,\Re{v}<1 $:
$$
L(u+v,\chi_{cd})\e{D}_3(u,v;h)=\left\{\e{D}_3^1+\e{D}_3^2+\e{D}_3^3
\right\}(u,v;h),\eqno(33.3)
$$
where
$$
\eqalignno{
\e{D}_3^1&=-{1\over\pi
}\sum_{c_1d_1=cd}{1\over d_1}
\int_{-\infty}^\infty\left({(d_1,d)\over\sqrt{d}}\right)^{1+2ir}
{f^{ir}\sigma_{-2ir}(f;\chi_{cd})\over|L(1+2ir,\chi_{cd})|^2}\cr
&\times\zeta(u+ir)\zeta(u-ir)\zeta(v+ir)\zeta(v-ir)\cr
&\times \prod_{p|(c_1,c)(d_1,d)}\left(\sigma_{-2ir}(f_p)
\left(1-{1\over p^{1+2ir}}\right)-1\right)
\prod_{p|d_1}
\left(1-{1\over p^{u-ir}}\right)
\left(1-{1\over p^{v-ir}}\right)\cr
&\times\prod_{p|c_1}
\left\{\left(1-{1\over p^{1-2ir}}\right)
\left(1-{1\over p^{u+v}}\right)
-\left(1-{1\over p^{u-ir}}\right)
\left(1-{1\over p^{v-ir}}\right)\right\}h(r)dr,\cr
\e{D}_3^2&=-2f^{1-u}\sigma_{2(u-1)}(f;\chi_{cd}){\zeta(u+v-1)
\zeta(v-u+1)\over L(3-2u,\chi_{cd})}h(i(u-1))\cr
&\times\sum_{c_1d_1=cd}{\varphi(c_1)\over c_1^{u+v}d_1}
\left({(d_1,d)\over\sqrt{d}}\right)^{3-2u}
\prod_{p|(c_1,c)(d_1,d)}\left(\sigma_{2(u-1)}(f_p)
\left(1-{1\over p^{3-2u}}\right)-1\right)\cr
&\times\prod_{p|d_1}\left(1-{1\over p^{u+v-1}}\right)\cr
&-2f^{1-v}\sigma_{2(v-1)}(f;\chi_{cd}){\zeta(u+v-1)
\zeta(u-v+1)\over L(3-2v,\chi_{cd})}h(i(v-1))\cr
&\times\sum_{c_1d_1=cd}{\varphi(c_1)\over c_1^{u+v}d_1}
\left({(d_1,d)\over\sqrt{d}}\right)^{3-2v}
\prod_{p|(c_1,c)(d_1,d)}\left(\sigma_{2(v-1)}(f_p)
\left(1-{1\over p^{3-2v}}\right)-1\right)\cr
&\times\prod_{p|d_1}\left(1-{1\over p^{u+v-1}}\right),\cr
\e{D}_3^3&=-2f^{u-1}\sigma_{2(1-u)}(f;\chi_{cd}){\zeta(u+v-1)
\zeta(v-u+1)\over L(3-2u,\chi_{cd})}h(i(u-1))\prod_{p|cd}
\left(1-{1\over p^{2u-1}}\right)^{-1}\cr
&\times\sum_{c_1d_1=cd}{\varphi(d_1)\over d_1^2}
\left({(d_1,d)\over\sqrt{d}}\right)^{2u-1}
\prod_{p|(c_1,c)(d_1,d)}\left(\sigma_{2(1-u)}(f_p)
\left(1-{1\over p^{2u-1}}\right)-1\right)\cr
\times\prod_{p|d_1}&\left(1-{1\over p^{v-u+1}}\right)\prod_{p|c_1}
\left\{\left(1-{1\over p^{3-2u}}\right)
\left(1-{1\over p^{u+v}}\right)-\left(1-{1\over p}\right)
\left(1-{1\over p^{v-u+1}}\right)\right\},\cr
&-2f^{v-1}\sigma_{2(1-v)}(f;\chi_{cd}){\zeta(u+v-1)
\zeta(u-v+1)\over L(3-2v,\chi_{cd})}h(i(v-1))\prod_{p|cd}
\left(1-{1\over p^{2v-1}}\right)^{-1}\cr
&\times\sum_{c_1d_1=cd}{\varphi(d_1)\over d_1^2}
\left({(d_1,d)\over\sqrt{d}}\right)^{2v-1}
\prod_{p|(c_1,c)(d_1,d)}\left(\sigma_{2(1-v)}(f_p)
\left(1-{1\over p^{2v-1}}\right)-1\right)\cr
\times\prod_{p|d_1}&\left(1-{1\over p^{u-v+1}}\right)
\prod_{p|c_1}
\left\{\left(1-{1\over p^{3-2v}}\right)
\left(1-{1\over p^{u+v}}\right)
-\left(1-{1\over p}\right)
\left(1-{1\over p^{u-v+1}}\right)\right\}.&(33.4)
}
$$
We see readily that $\e{D}_3^1$ and $\e{D}_3^2$ are regular
at $\left({1\over2},{1\over2}\right)$. As to $\e{D}_3^3$, the factors 
$\prod_{p|cd}(1-p^{1-2u})^{-1}$ and $\prod_{p|cd}(1-p^{1-2v})^{-1}$
diverge at the point unless
$cd=1$; however, $\e{D}_3^3$ itself must be regular there, for
$L(u+v,\chi_{cd})\e{D}_1$,
$L(u+v,\chi_{cd})\e{D}_2$ are regular, and thus 
$L(u+v,\chi_{cd})\e{D}_3$ as well.
\medskip
Hence, from $(31.1)$, $(31.6)$,
$(33.1)$, and $(33.3)$, we obtain 
$$
\e{D}\left(\txt{1\over2},\txt{1\over2};h\right)=
{cd\over\varphi(cd)}\left\{\e{D}_2^1+\e{D}_2^2+\e{D}_2^3
+\e{D}_2^4+\e{D}_3^1+\e{D}_3^2+\e{D}_3^3\right\}
\left(\txt{1\over2},\txt{1\over2};h\right).\eqno(33.5)
$$ 
\bigskip
\noindent
{\bf 34.} The last equation gives
\medskip
\noindent
{\bf Lemma 16.} {\it
We have, with the weight $h$ as above,
$$
\sum_j\epsilon_j\varrho_j(f;1/c)R_j\left(\txt{1\over2},0\right)
{h(\kappa_j)\over\cosh\pi\kappa_j}=\sum_{a=1}^7\e{H}_a(f;h),\eqno(34.1)
$$
where
$$
\eqalignno{
\e{H}_1&={2cd\over\pi^3i\varphi(cd)}\Big\{(c_E-\log(2\pi\sqrt{f/d}))
(\hat{h})'(\txt{1\over2})
+\txt{1\over4}(\hat{h})''(\txt{1\over2})\Big\}
\tau(f,\chi_{cd})f^{-{1\over2}}\prod_{p|c}
\left(\tau(f_p)\left(1-{1\over p}\right)-1\right),\cr
\e{H}_2&={c\sqrt{d}\over\pi^3\varphi(cd)
}\sum_m m^{-{1\over2}}
\tau(m)\tau(dm+f;\chi_{cd})\Psi_+\left(dm/f;h\right)
\prod_{p|c}\left(\!\tau((dm+f)_p)
\left(1-{1\over p}\right)-1\right),\cr
\e{H}_3&={c\sqrt{d}\over\pi^3\varphi(cd)}
\sum_{\scr{m}\atop\scr{dm\ne f}}m^{-{1\over2}}
\tau(m)\tau(dm-f;\chi_{cd})\Psi_-(dm/f;h)
\prod_{p|c}\left(\!\tau((dm-f)_p)
\left(1-{1\over p}\right)-1\right),\cr
\e{H}_4&=-{\delta_{d,1}\over2\pi^3}
f^{-{1\over2}}\tau(f)\Psi_-(1;h),\cr
\e{H}_a&={cd\over\varphi(cd)}\e{D}_3^{a-4}
\left(\txt{1\over2},\txt{1\over2};h\right),\quad 5\le a\le 7.&(34.2)
}
$$
Here $\tau$ is the divisor
function, $\Psi_\pm(x;h)=\Psi_\pm\left(\txt{1\over2},
\txt{1\over2};x;h\right)$, and $\e{H}_4$ vanishes unless $d=1$.
}
\medskip
\noindent
This is a counterpart of [11, Lemma 3.8], and follows immediately from
$(31.4)$, $(33.2)$ and $(33.5)$. We have left
$\e{H}_a(f;h)$, $5\le a\le7$, without computing it
explicitly, because it seems better to avoid
the highly complicated computation of
$\e{D}_3^3\left(\txt{1\over2},\txt{1\over2};h\right)$ caused by
the two products over $p|cd$ mentioned above; and in fact
those $\e{H}_a(f;h)$ will readily turn out to be negligible in our
application of $(34.1)$ to be given in the next section. 
\medskip
From [11, pp.\ 119--121], we quote the following:
$$
(\hat{h})'(\txt{1\over2})=2\int_{-\infty}^\infty{r}
h(r){{\Gamma'}\over\Gamma}
(\txt{1\over2}+ir)dr,\quad
(\hat{h})''(\txt{1\over2})=4\int_{-\infty}^\infty{r}h(r)
\Big\{{{\Gamma'}\over\Gamma}(\txt{1\over2}+ir)\Big\}^2dr,
\eqno(34.3)
$$
$$
\Psi_+(x;h)=2\pi\int_0^1\big\{y(1-y)(1+y/x)\big\}^{-{1\over2}}
\int_{-\infty}^\infty{r}
h(r)\tanh(\pi{r})\Big\{{{y(1-y)}\over{x+y}}\Big\}^{ir}drdy.
\eqno(34.4)
$$
For $x>1$
$$
\Psi_-(x;h)=2\pi{i}\int_0^1\big\{y(1-y)(1-y/x)\big\}^{-{1\over2}}
\int_{-\infty}^\infty{{{r}h(r)}\over{\cosh(\pi{r})}}
\Big\{{{y(1-y)}\over{x-y}}\Big\}^{ir}drdy.
\eqno(34.5)
$$
For $x=1$
$$
\Psi_-(1;h)=2\pi^2\int_{-\infty}^\infty{r}h(r)
{{\sinh(\pi{r})}\over{(\cosh(\pi{r}))^2}}dr.
\eqno(34.6)
$$
For $0<x<1$
$$\eqalignno{
\Psi_-(x;h)=\int_0^\infty&\Big\{\int_{(\beta)}x^s(y(y+1))^{s-1}
{{\Gamma({1\over2}-s)^2}\over{\Gamma(1-2s)\cos(\pi{s})}}ds\Big\}\cr
&\times\Big\{
\int_{-\infty}^\infty{r}h(r)\Big({y\over{1+y}}\Big)^{ir}dr\Big\}dy,
\quad -\txt{3\over2}<\beta<\txt{1\over2},\;\beta\ne-\txt{1\over2}.
&(34.7)
}
$$
\bigskip
\noindent
{\bf 35.} We shall continue our discussion, adopting the argument given
in [11, pp.\ 124--130]. Thus we first state the following approximation for
$L_j\left(\txt{1\over2};1/c\right)$:  Let $K$ tend to
infinity, and assume that
$$
|\kappa_j-K|\le{G}\log{K}\eqno(35.1)
$$
with
$$
K^{{1\over2}+\delta}<G<K^{1-\delta},\quad 0<\delta<\txt{1\over2}.
\eqno(35.2)
$$
Then we have, for any $N\ge1$ and $\lambda=C\log{K}$ 
with a sufficiently large $C>0$,
$$
\eqalignno{
&L_j\left(\txt{1\over2};1/c\right)=
\sum_{f\le3K\sqrt{cd}}\varrho_j(f;1/c)f^{-{1\over2}}
\exp(-(f/(K\sqrt{cd}))^\lambda)\cr
-&\sum_{\nu=0}^{N_1}\sum_{f\le{3K\sqrt{cd}}}\varrho_j(f;1/c)
f^{-{1\over2}}U_\nu(f/(K\sqrt{cd}))(1-(\kappa_j/K)^2)^\nu
+O(K^{-{1\over5}N}+K^{-{1\over2}C}),\quad&(35.3)\cr
}
$$
with the implied constant depending only on 
$\delta$, $C$, and $N$. Here $N_1=[3N/\delta]$ and
$$
U_\nu(x)={1\over{2\pi{i}\lambda}}\int_{(-\lambda^{-1})}
(4\pi^2 x)^wu_\nu(w)\Gamma(w/\lambda)dw,
\eqno(35.4)
$$
where $u_\nu(w)$ is a polynomial of degree $\le2N_1$, 
whose coefficients are independent of $\kappa_j$ and bounded 
by a constant depending only on $\delta$ and $N$.
\medskip
\noindent
In fact, this assertion is a counterpart of [11, Lemma 3.9] and
the proof is analogous; the necessary change is only in that
we now use $(23.13)$ instead under the assumption $\varpi_j\epsilon_j
=+1$ as $L\left({1\over2};1/c\right)=0$ if $\varpi_j\epsilon_j=-1$.
\medskip
With this, we now set, in $(34.1)$,
$$
h(r)=\left(r^2+\txt{1\over4}\right)\left\{
\exp(-((r-K)/G)^2)+\exp(-((r+K)/G)^2)\right\}.\eqno(35.5)
$$
We have
$$
\eqalignno{
(\hat{h})'\left(\txt{1\over2}\right)&=2i\pi^{3\over2}K^3G+O(KG^3),\cr
(\hat{h})''\left(\txt{1\over2}\right)&=8i\pi^{3\over2}K^3G\log K
+O(KG^3\log K).&(35.6)
}
$$
(see [11, p.\ 129]).
\medskip
We have, by $(34.1)$ and $(35.4)$,
$$
\eqalignno{
\sum_j\epsilon_jR_j\left(\txt{1\over2},0\right)
L_j&\left(\txt{1\over2};1/c\right){h(\kappa_j)\over\cosh\pi\kappa_j}=
\sum_{f\le3K\sqrt{cd}}f^{-{1\over2}}
\exp(-(f/(K\sqrt{cd})^\lambda)\sum_{a=1}^7
\e{H}_a(f;h)\cr
-&\sum_{\nu\le N_1}\sum_{f\le3K\sqrt{cd}}f^{-{1\over2}}
U_\nu(f/(K\sqrt{cd}))\sum_{a=1}^7
\e{H}_a(f;h_\nu)+O(1)
,&(35.7)
}
$$
where the five terms correspond to those on the right side of $(34.1)$,
respectively, with the present $h$ and $h_\nu(r)=h(r)(1-(r/K)^2)^\nu$. 
Since we have imposed $(35.1)$--$(35.2)$, those terms with $\nu\ge1$
can actually be ignored, and it suffices to consider instead
$$
\eqalignno{
\sum_{f\le3K\sqrt{cd}}&f^{-{1\over2}}
\exp(-(f/(K\sqrt{cd})^\lambda)\sum_{a=1}^7
\e{H}_a(f;h)\cr
-&\sum_{f\le3K\sqrt{cd}}f^{-{1\over2}}
U_0(f/(K\sqrt{cd}))\sum_{a=1}^7\e{H}_a(f;h).&(35.8)
}
$$
The discussion in [11, pp.\ 128--129] works just fine with our present
situation as well; and the contribution of $\e{H}_a$, $a=2,3,4$, turns out
to be negligible. 
\medskip
\noindent
{\csc Remark 6}. However, if the uniformity in the Stufe $cd$ is
required, then this part of our argument should become subtle. 
\medskip
As to $\e{H}_1$, its contribution to $(35.8)$ is equal to
$$
{4cd\over\pi^{3\over2}\varphi(cd)}
K^3G\left(\e{K}_1+\e{K}_2\right)
+O(KG^3(\log K)^2),\eqno(35.9)
$$
where we have used $(35.6)$, and 
$$
\eqalignno{
\e{K}_1&=\sum_f
\left(c_E-\log(2\pi\sqrt{f/d})+\log K\right)
\exp(-(f/(K\sqrt{cd})^\lambda)\cr
&\qquad\times {\tau(f;\chi_{cd})\over f}
\prod_{p|c}\left(\tau(f_p)\left(1-{1\over p}\right)-1\right),\cr
\e{K}_2&=-\sum_f
\left(c_E-\log(2\pi\sqrt{f/d})+\log K\right)
U_0(f/(K\sqrt{cd}))\cr
&\qquad\times {\tau(f;\chi_{cd})\over f}
\prod_{p|c}\left(\tau(f_p)\left(1-{1\over p}\right)-1\right).
&(35.10)
}
$$
To compute $\e{K}_1$, $\e{K}_2$, let us put
$$
z(s)=\sum_f {\tau(f;\chi_{cd})\over
f^{s+1}}\prod_{p|c}\left(\tau(f_p)\left(1-{1\over p}\right)-1\right).
\eqno(35.11)
$$
Then
$$
\eqalignno{
\e{K}_1&={1\over 2\pi i\lambda}\int_{(1)} 
\left\{(\log(K\sqrt{d}/2\pi)+c_E)
z(s)+\txt{1\over2}z'(s)\right\}(K\sqrt{cd})^s\Gamma(s/\lambda)ds,\cr
\e{K}_2&=-{1\over 2\pi i\lambda}\int_{(-1)} 
\left\{(\log(K\sqrt{d}/2\pi)+c_E)
z(-s)+\txt{1\over2}z'(-s)\right\}\cr
&\qquad\times(4\pi^2/K\sqrt{cd})^{s}u_0(w)\Gamma(s/\lambda)ds.&(35.12)
}
$$
The latter can be replace by
$$
-{1\over 2\pi i\lambda}\int_{(-1)} 
\left\{(\log(K\sqrt{d}/2\pi)+c_E)
z(-s)+\txt{1\over2}z'(-s)\right\}
(4\pi^2/K\sqrt{cd})^{s}\Gamma(s/\lambda)ds\eqno(35.13)
$$
with an admissible error (see [11, p.\ 127] for a description of $u_0$).
\medskip
We have
$$
z(s)=\zeta(s+1)^2\prod_{p|d}\left(1-{1\over p^{s+1}}\right)
\prod_{p|c}\left(1-{1\over p^{s+1}}-{2\over p^{s+2}}+
{1\over p^{2s+2}}+{1\over p^{2s+3}}\right).
\eqno(35.14)
$$
Hence, we get 
$$
\e{K}_1,\e{K}_2
\sim\txt{1\over3}(\log K)^3{\varphi(cd)\over cd}\prod_{p|c}
\left(1-{1\over p^2}\right).\eqno(35.15)
$$
\bigskip
\noindent
{\bf 36.} It remains for us to deal with $\e{H}_a$, $5\le a\le7$.
\medskip
We have obviously
$$
\e{D}_3^1\left(\txt{1\over2},\txt{1\over2};h\right)
\ll\tau(f)\int_{-\infty}^\infty
{\left|\zeta\left({1\over2}+ir\right)\right|^4\over|\zeta(1+2ir)|^2}
h(r)dr\ll \tau(f)K^3(\log K)^5,\eqno(36.1)
$$
which can of course be replaced by a better bound, but for our
purpose this is sufficient. We see that the contribution of
$\e{H}_5$ to $(35.8)$ is $\ll K^{7\over2}(\log K)^7$, which is
negligible in view of $(35.2)$, $(35.9)$, and $(35.15)$. 
\medskip
As to $\e{H}_6$ and $\e{H}_7$, we shall treat the latter only, for
the former is analogous and certainly easier than the latter.
As we have remarked already, $\e{D}_3^3(u,v;h)$ is regular in the
vicinity of $\left({1\over2},{1\over2}\right)$. Thus we have
$$
\e{D}_3^3\left(\txt{1\over2},\txt{1\over2};h\right)=
-{1\over(2\pi)^2}\int_{C_2}\!\int_{C_1}
{\e{D}_3^3(u,v;h)\over \left(u-{1\over2}\right)
\left(v-{1\over2}\right)}dudv,\eqno(36.2)
$$
where 
$$
C_1: \, \left|u-\txt{1\over2}\right|={1\over B(2+\log cd)},\quad
C_2: \, \left|v-\txt{1\over2}\right|={1\over 2B(2+\log cd)},\eqno(36.3)
$$
with a sufficiently large constant $B$. This integrand is, by the explicit
formula for $\e{D}_3^3(u,v;h)$ in $(33.4)$,
$$
\ll \exp\left(-\txt{1\over2}(K/G)^2\right),\eqno(36.4)
$$
and $\e{H}_7$ is negligible.
\medskip
Hence we have obtained
\medskip
\noindent
{\bf Lemma 17.} {\it Let $h$ be as in $(35.5)$ with 
$(35.1)$--$(35.2)$. Then we have,  for any fixed $c,d$
with $\mu(cd)\ne0$, 
$$
\sum_{\scr{\kappa_j}\atop\scr{\kappa_j^2+\txt{1\over4}\in{\rm
Sp}(\varGamma_0(cd))}}\epsilon_jR_j\left(\txt{1\over2},0\right)
L_j\left(\txt{1\over2};1/c\right){h(\kappa_j)\over\cosh\pi\kappa_j}
\sim{8\over3\pi^{3\over2}}GK^3(\log K)^3
\prod_{p|c}\left(1-{1\over p^2}\right).\eqno(36.5)
$$
}
\medskip
\noindent
In particular, if $A$ is fixed, we have
$$
\sum_\kappa\e{R}(\kappa;A){h(\kappa)\over\cosh\pi\kappa}\sim
{8\over3\pi^{3\over2}}GK^3(\log K)^3\sum_{c,d}\e{A}(c,d)
\prod_{p|c}\left(1-{1\over p^2}\right),\eqno(36.6)
$$
where ${1\over4}+\kappa^2\in\bigcup_{c,d}\r{Sp}(\varGamma_0(cd))$
with $\mu(cd)\ne0$.
\bigskip
Therefore we have established
\medskip
\noindent
{\bf Theorem.}
{\it Provided $\alpha_n>0$ for square-free $n$ and $=0$ otherwise, 
the function $Z_2(s;A)$ has infinitely many simple poles on the line
$\Re s={1\over2}$. 
}
\bigskip
\noindent
This restriction on the support of $\alpha_n$ will be lifted in our
forthcoming work.
\medskip
Our result suggests that the Mellin
transform
$$
Z_3(s;1)=\int_1^\infty\left|\zeta\left(\txt{1\over2}+it\right)\right|^6
t^{-s}dt\eqno(36.7)
$$
should have the line $\Re s={1\over2}$ as a natural boundary, for
$|\zeta|^6=|\zeta|^4|\zeta|^2$
and $|\zeta|^2$ may be replaced by a finite expression similar to
$|A|^2$ via the approximate functional equation.
The same was speculated also by a few people other than us, but it
appears that our theorem is so far the sole explicit evidence supporting
this conjectural assertion. At any event, in view of
of {\csc Remark 5} above, it appears reasonable for us to maintain that
$Z_3(s;1)$ does not continue beyond the imaginary axis.
\bigskip
This entails 
\medskip
\noindent
{\bf Problems:} 
\medskip
\noindent
(1) Is the set $\bigcup_{q\ge1}\r{Sp}(\varGamma_0(q))$ dense in the
positive real axis?
\smallskip
\noindent
(2) Is the set of $\kappa$ satisfying $(30.3)$ dense in the positive
real axis?
\smallskip
\noindent
(3) Is the set of $\kappa$ satisfying $(30.3)$ dense in 
any half line?
\smallskip
\noindent
(4) Is the set of $\kappa$ satisfying $(30.3)$ dense in any interval
whose left end point is the origin?
\bigskip
\noindent
Obviously (1) is to be solved first and (2) must be far more 
difficult than (1). The third, weaker than (2), appears highly plausible
in the light of Lemma 17; on the other hand our method does not seem to
extend without new twists so as to include the situation of (4), i.e., the
detection of low lying poles.
\bigskip
\noindent
{\csc Addendum.} Recently C.P. Hughes and M.P. Young (arXiv:0709.2345
[math.NT]) obtained an asymptotic formula for the mean value $(30.1)$
where the length of $A$ is less than $T^\eta$ with any
fixed $\eta<1/11$. They did not employ the spectral theory of Kloosterman
sums. Our method should give a better result than theirs, 
if it is combined with works by N. Watt on this mean value.
\vskip 1cm
\centerline{\bf References}
\bigskip
\item{[1]} R. Bruggeman. Kloosterman sums (v.5). 
September 2007.
\item{[2]} R. Bruggeman and Y. Motohashi. A new approach to the
spectral theory of the fourth moment of the Riemann zeta-function.
J. reine angew.\ Math., {\bf 579} (2005), 75--114.
\item{[3]} J.-M. Deshouillers and H. Iwaniec.
Kloosterman sums and Fourier coefficients of cusp forms. 
Invent.\ math., {\bf 70} (1982), 219--288.
\item{[4]} ---. Power mean-values  for
Dirichlet's polynomials and the Riemann zeta-function.\ I. 
\quad  Mathematika, {\bf 29} (1982), 202--212.
\item{[5]} A. Good. The convolution method  for Dirichlet series.
Contemp.\ Math., {\bf 53} (1986), 207--214.
\item{[6]} D.A. Hejhal. {\it The Selberg Trace Formula for 
${\rm PSL}(2,\B{R})$\/}.\ II. Lect.\ Notes in Math.,
vol.\ {\bf 1001}, Springer-Verlag, Berlin, 1983.
\item{[7]} Y. Motohashi. An explicit formula for the fourth power mean of
the Riemann zeta-function. Acta Math., {\bf 170} (1993), 181--220.
\item{[8]} ---. A relation between the Riemann zeta-function
and the hyperbolic Laplacian. Ann.\ Scuola Norm.\ Sup.\ Pisa.\ (4)
{\bf 22} (1995), 299--313.
\item{[9]} ---. The Riemann zeta-function and Hecke
congruence subgroups. RIMS Kyoto Univ.\ Kokyuroku, {\bf 958} (1996),
166--177.
\item{[10]} ---. The mean square of Dedekind zeta-functions
of quadratic number fields. In:  {\it Sieve Methods, Exponential
Sums and their Applications in Number Theory\/} 
(G.R.H. Greaves, G. Harman, M.N. Huxley, editors) Cambridge University
Press, Cambridge 1997, pp.\ 309--323.
\item{[11]} ---. {\it Spectral Theory of the Riemann Zeta-Function\/}.
Cambridge University Press,
Cambridge 1997.
\item{[12]} ---. A note on the mean value of the zeta and
$L$-functions.\ XV. Proc.\ Japan Acad.,  
{\bf 83A} (2007), 73--78.

\vskip 1cm
\noindent
Department of Mathematics,
Nihon University
\par\noindent
Surugadai, Tokyo 101-8308, JAPAN
\smallskip\noindent
ymoto@math.cst.nihon-u.ac.jp
\par
\noindent
www.math.cst.nihon-u.ac.jp/$\!\sim$ymoto/

\bye